\documentclass[]{article}

\usepackage{subfiles}
\usepackage{graphicx}
\graphicspath{ {./images/} }

\usepackage[margin=1in]{geometry}
\usepackage{amsmath}
\usepackage{amsfonts}
\usepackage{amssymb}
\usepackage{amsthm}
\usepackage{physics}
\usepackage{ulem}
\usepackage{tikz}
\usetikzlibrary{calc}
\usetikzlibrary{arrows.meta}
\usepackage{hyperref}
\usepackage{marginnote}
\usepackage{subcaption}
\usepackage{mathtools}
\usepackage[backend=bibtex,
style=alphabetic,sorting=nyt,maxbibnames=10,maxalphanames=4,isbn=false,url=false,doi=false]{biblatex}
\DeclareFieldFormat*{title}{\mkbibitalic{#1}}
\DeclareFieldFormat*{journaltitle}{#1\isdot}
\DeclareFieldFormat*{booktitle}{#1\isdot}

\theoremstyle{plain}% default
\newtheorem{thm}{Theorem}[section]
\newtheorem{lem}[thm]{Lemma}
\newtheorem{prop}[thm]{Proposition}
\newtheorem{cor}[thm]{Corollary}

\newtheorem{manualtheoreminner}{Theorem}
\newenvironment{manualtheorem}[1]{%
	\renewcommand\themanualtheoreminner{#1}%
	\manualtheoreminner
}{\endmanualtheoreminner}

\theoremstyle{definition}
\newtheorem{defn}[thm]{Definition}

\theoremstyle{remark}
\newtheorem{rem}[thm]{Remark}
\newtheorem{exmp}[thm]{Example}

\DeclareMathOperator{\ad}{ad}
\DeclareMathOperator{\Ad}{Ad}
\DeclareMathOperator{\conj}{conj}
\renewcommand{\trace}{\mathrm{trace}}

\DeclareMathOperator{\Stab}{Stab}

\DeclareMathOperator{\End}{End}

\DeclareMathOperator{\id}{id}
\DeclareMathOperator{\mfa}{\mathfrak{a}}
\DeclareMathOperator{\mfp}{\mathfrak{p}}
\DeclareMathOperator{\mfk}{\mathfrak{k}}
\DeclareMathOperator{\mfg}{\mathfrak{g}}
\DeclareMathOperator{\mfn}{\mathfrak{n}}

\DeclareMathOperator{\mfe}{\mathfrak{e}}

\DeclareMathOperator{\orb}{\mathrm{orb}}
\DeclareMathOperator{\X}{\mathbb{X}}
\DeclareMathOperator{\Y}{\mathbb{Y}}

\DeclareMathOperator{\Flag}{Flag}
\DeclareMathOperator{\Flagt}{Flag(\tau_{mod})}

\DeclareMathOperator{\st}{st}
\DeclareMathOperator{\ost}{ost}
\DeclareMathOperator{\midp}{mid}
\DeclareMathOperator{\interior}{int}
\DeclareMathOperator{\Hess}{Hess}
\DeclareMathOperator{\R}{\mathbb{R}}
\DeclareMathOperator{\HH}{\mathbb{H}}
\DeclareMathOperator{\im}{im}
\DeclareMathOperator{\length}{length}
\DeclareMathOperator{\sl2r}{\mathfrak{sl}(2,\mathbb{R})}

\DeclareMathOperator{\eva}{ev}

\DeclareMathOperator{\SL}{SL}

\DeclareMathOperator{\SO}{SO}

\renewcommand{\d}[1]{\ensuremath{\mathrm{d} {#1}}}

\newcommand{\pnorm}[1]{\ensuremath{\left\lvert #1 \right\rvert_{B_p}}}
\newcommand{\qnorm}[1]{\ensuremath{\left\lvert #1 \right\rvert_{B_q}}}

\newcommand{\deriv}[1]{
	\ensuremath{\frac{\mathrm{d}}{\mathrm{d} {#1}} }}
\newcommand{\derivtwo}[1]{
	\ensuremath{\frac{\mathrm{d}^2}{\mathrm{d} {#1}^2} }}
\newcommand{\atzero}[1]{
	\ensuremath{\mid_{{#1}=0} }}

\newcommand{\lied}{\mathcal{L}}
\renewcommand{\grad}{\mathrm{grad}}
\newcommand{\evp}{\mathrm{ev}_p}

\title{A quantified local-to-global principle for Morse quasigeodesics}
\author{J. Maxwell Riestenberg}

\begin{document}

\maketitle

\begin{abstract}
	In \cite{KLP14}, Kapovich, Leeb and Porti gave several new characterizations of Anosov representations $\Gamma \to G$, including one where geodesics in the word hyperbolic group $\Gamma$ map to ``Morse quasigeodesics" in the associated symmetric space $G/K$. In analogy with the negative curvature setting, they prove a local-to-global principle for Morse quasigeodesics and describe an algorithm which can verify the Anosov property of a given representation in finite time. However, some parts of their proof involve non-constructive compactness and limiting arguments, so their theorem does not explicitly quantify the size of the local neighborhoods one needs to examine to guarantee global Morse behavior. In this paper, we supplement their work with estimates in the symmetric space to obtain the first explicit criteria for their local-to-global principle. This makes their algorithm for verifying the Anosov property effective. As an application, we demonstrate how to compute explicit perturbation neighborhoods of Anosov representations with two examples.
\end{abstract}

\section{Introduction}

Anosov representations were introduced by Labourie and defined in general by Guichard and Wienhard \cite{Lab06,GW12}. An Anosov representation is a homomorphism from a word hyperbolic group $\Gamma$ to a semisimple Lie group $G$ satisfying a strong dynamical condition. These representations have come to be widely studied as an interesting source of infinite covolume discrete subgroups of higher rank semisimple Lie groups, see the surveys \cite{Kas18,KL18c}. This paper is concerned with certifying the Anosov property of a given representation. For some well-studied examples of Anosov representations, such as Hitchin representations and maximal representations of surface groups, the Anosov property can be certified via coarse topological invariants \cite{BIW03}. However, in the most general setting, deciding whether a given representation is Anosov is difficult. Building on the work of Kapovich, Leeb and Porti in \cite{KLP14}, we give here the first explicit, finite criteria that certify the Anosov property for a general representation.

One important property of Anosov representations is stability: any sufficiently small perturbation of an Anosov representation remains Anosov. It can happen that a connected component of the representation space consists entirely of Anosov representations, such as the Hitchin component, or the components consisting of maximal representations of surface groups, see also \cite{Wie18,GW18}. In these cases, the Anosov condition is closed: every deformation of such a representation remains Anosov. However the Anosov condition is not closed in general. For instance, given an Anosov representation of a free group, or the representations of surface groups studied by Barbot in \cite{Bar10}, it is unclear how large to expect Anosov neighborhoods to be. As an application of our main result, we demonstrate how to construct explicit perturbation neighborhoods of a given Anosov representation with two examples, see Theorem \ref{main theorem demonstration 1} and Theorem \ref{main theorem demonstration 2}.

Anosov representations have come to be viewed as the appropriate generalization to higher rank semisimple Lie groups of convex cocompact actions on rank $1$ symmetric spaces. Indeed, when $G$ has real rank $1$, a representation of a finitely generated group is Anosov if and only if it has finite kernel and the image is convex cocompact, i.e.\ acts cocompactly on a nonempty convex subset of the associated negatively curved symmetric space. A finitely generated group of isometries of a negatively curved symmetric space is convex cocompact if and only if it is undistorted, i.e.\ any orbit map is a quasi-isometric embedding. By the Morse Lemma in hyperbolic geometry, geodesics in $\Gamma$ then map within uniformly bounded neighborhoods of geodesics in the symmetric space. Moreover, the Morse Lemma implies a local-to-global principle for quasigeodesics, allowing one to establish finite criteria for a finitely generated group to be undistorted. One can then exhaust the group by balls in the Cayley graph and if any such ball passes a finite check then the group is undistorted. This is a semidecidable algorithm to verify undistortion: if the group is undistorted, this algorithm will eventually terminate and certify so; otherwise, it will run on forever. 

The naive generalization of convex cocompactness to higher rank turns out to be too restrictive. For example, the work of Kleiner and Leeb and independently Quint implies that a Zariski dense, discrete subgroup of a higher rank simple Lie group which acts cocompactly on a convex subset of the associated symmetric space is a uniform lattice  \cite{KL06,Q05}. On the other hand, the naive generalization of undistortion to higher rank turns out to be too loose: In his thesis, Guichard described an example of an undistorted subgroup in $\SL(2,\R) \times \SL(2,\R)$ which is unstable, in the sense that representations arbitrarily close to the inclusion fail to have discrete image \cite{Gui04}, see also \cite{GGKW17}. In \cite{KLP14}, Kapovich, Leeb and Porti describe an example of a discrete undistorted subgroup of $\SL(2,\R) \times \SL(2,\R)$ which is finitely generated but not finitely presentable, using work Baumslag and Roseblade \cite{BR84}. The Anosov property strikes a balance between these two naive generalizations to give a large class of representations that still exhibit good behavior. We will be concerned with a newer characterization that directly strengthens the undistortion condition. 

In \cite{KLP14}, Kapovich, Leeb and Porti gave several new characterizations of Anosov representations generalizing some of the many characterizations of convex cocompact subgroups. We will use their characterization, called Morse actions, that strengthens the undistortion condition by requiring geodesics in $\Gamma$ to map to \textit{Morse quasigeodesics}, described below. They prove a suitable generalization of the local-to-global principle for Morse quasigeodesics in higher rank symmetric spaces, see Theorem \ref{main theorem} below. They then show the Anosov property is semi-decidable by describing an algorithm which can certify the Anosov property of a given representation of a word hyperbolic group in finite time. However, some parts of their proof involve non-constructive compactness and limiting arguments, so their theorem does not explicitly quantify the size of the local neighborhoods one needs to examine to guarantee global Morse behavior. In order to implement their algorithm, one needs a quantified version of the local-to-global principle as we give here.

Roughly speaking, a quasigeodesic is \textit{Morse} if every finite consecutive subsequence is uniformly close to a \textit{diamond}, which plays the role of a geodesic segment in rank $1$. These diamonds are intersections of Weyl cones, see  Sections \ref{visual boundary} and \ref{straight and spaced sequences are morse}, and may also be characterized as unions of Finsler geodesic segments, see \cite{KL18a,KL18b}. An infinite Morse quasiray stays within a uniformly bounded neighborhood of a Weyl cone, which plays the role of a geodesic ray in rank $1$, and a bi-infinite Morse quasigeodesic stays within a uniformly bounded neighborhood of a parallel set, which plays the role of a geodesic line in rank $1$, see Section \ref{section parallel sets and horocycles}. The precise definition of Morse quasigeodesic is given in Section \ref{local to global section}.

The main result of this paper is a quantified version of the following theorem due to Kapovich, Leeb and Porti. We let $\X$ denote a symmetric space of noncompact type.
\begin{thm}[{\cite[Theorem 7.18]{KLP14}}]\label{main theorem}
	For any $\alpha_{new}<\alpha_0,D,c_1,c_2,c_3,c_4$, there exists a scale $L$ so that every $L$-local $(\alpha_0,\tau_{mod},D)$-Morse $(c_1,c_2,c_3,c_4)$-quasigeodesic in $\X$ is an $(\alpha_{new},\tau_{mod},D')$-Morse $(c_1',c_2',c_3',c_4')$-quasigeodesic.
\end{thm}
We reprove Theorem \ref{main theorem} and obtain the first explicit estimate of $L$. This appears in Theorem \ref{local to global}, which depends on Theorem \ref{straight and spaced sequences are morse} and Theorem \ref{midpoint sequences}. The theorem statements involve several auxiliary parameters and inequalities, so they are too cumbersome to give here. In order to apply our quantified version of the local-to-global principle and obtain an explicit scale $L$, one must produce auxiliary parameters satisfying these inequalities; this process is tedious but easy, as we discuss in Section \ref{section applications}. Versions of Theorems \ref{straight and spaced sequences are morse} and \ref{local to global} without explicit conditions are also proved in \cite{KLP14}. 

As a demonstration of our techniques, we compute explicit perturbation neighborhoods of two Anosov representations into $\SL(3,\R)$. To quantify the distance between linear representations we use the Frobenius norm on the generators: for a matrix $A$, let $\abs{A}_{Fr}^2 = \trace(A^T A)$. In both cases we control the orbit map at a basepoint; the Frobenius norm is closely related to distances to that basepoint, see Section \ref{section surface group}. The first example is a neighborhood of Anosov representations of a free group. 

%corrected 3.5.2021
\begin{thm}\label{main theorem demonstration 1}
	Let $\Gamma_1$ be the subgroup of $\SL(3,\R)$ generated by 
$$ g= \begin{bmatrix} e^t & 0 & 0 \\ 0 & 1 & 0 \\ 0& 0 & e^{-t} \end{bmatrix} , \quad  h= \begin{bmatrix} \cosh t & 0 & \sinh t \\ 0 & 1 & 0 \\ \sinh t & 0 &  \cosh t \end{bmatrix} ,$$
with $\tanh t=0.75$. If $\Gamma_1'$ is generated by $g',h'$ where $\max \{ \abs{g-g'}_{Fr},\abs{h-h'}_{Fr} \} \le  10^{-15,309} $, then $\Gamma_1'$ is Anosov.
\end{thm}

The second example is a neighborhood of Anosov representations of a closed surface group. Let $\Gamma_2$ be the subgroup of $\SL(3,\R)$ generated by 
$$ S = \left\{ \begin{bmatrix} \cos \theta & 0 & \sin \theta \\ 0 & 1 & 0 \\ -\sin \theta & 0 & \cos \theta \end{bmatrix}  \begin{bmatrix} \lambda & 0 & 0 \\ 0 & 1 & 0 \\ 0 & 0 & \lambda^{-1} \end{bmatrix} 
\begin{bmatrix} \cos \theta & 0 & -\sin \theta \\ 0 & 1 & 0 \\ \sin \theta & 0 & \cos \theta \end{bmatrix} \bigg\vert \, \theta \in \left\{0,\frac{\pi}{8},\frac{\pi}{4},\frac{3\pi}{8}\right\} \right\} $$
for $\log \lambda = \cosh^{-1} (\cot \frac{\pi}{8})  $. This group is isomorphic to the fundamental group of a closed surface of genus $2$, see Section \ref{section surface group}. In the statement of Theorem \ref{main theorem demonstration 2}, we control the perturbed representation on a larger generating set $S' = \{ \gamma \in \Gamma_2 \mid \sqrt{6} \abs{ \log \gamma }_{Fr} \le 9.5 \}$. The finite set $S'$ contains the standard generating set $S$ and consists of the elements of $\Gamma_2$ which move a basepoint $p$ in the symmetric space associated to $\SL(3,\R)$ by a distance of at most $9.5$. This basepoint is the point stabilized by $\SO(3)$. Using this larger generating set allows us to perturb the initial representation farther.  

%corrected 3.8.2021
\begin{thm}\label{main theorem demonstration 2}
	If $\rho \colon \Gamma_2 \to \SL(3,\R)$ is a representation satisfying $\abs{ \rho(s) - s}_{Fr} \le 10^{-3,698,433}$ for all $s \in S'$, then $\rho$ is Anosov. 
\end{thm}

We briefly sketch the proof of Theorems \ref{main theorem demonstration 1} and \ref{main theorem demonstration 2}. Let $\Gamma$ denote either $\Gamma_1$ or $\Gamma_2$. In either case the group $\Gamma$ acts cocompactly on a closed convex subset of a copy of the hyperbolic plane embedded totally geodesically in the symmetric space associated to $\SL(3,\R)$. We find explicit quasiisometry constants and by the classical Morse Lemma, there exists $R>0$ such that the orbit of any geodesic in $\Gamma$ is within $R$ of a geodesic. We slightly relax the Morse quasiisometric parameters of $\Gamma$ and apply the local-to-global principle Theorem \ref{local to global}. This provides a lower bound on $k$ such that any $2k$-local Morse quasigeodesic is a global Morse quasigeodesic. We control the perturbation of words of length $k$ in terms of the perturbation of the generators, completing the proof. 

We emphasize that our approach is completely general, in the following sense. Let $\rho \colon \Gamma \to G$ be any Anosov representation such that the orbit map at $p \in \X$ has known Morse quasiisometry parameters with respect to a finite symmetric generating set $S$ for $\Gamma$. We may then easily produce explicit parameters $k,\epsilon$ such that: if any other representation $\rho' \colon \Gamma \to G$ satisfies $d(\rho(\gamma)p,\rho'(\gamma)p) \le \epsilon$ for all $\gamma \in \Gamma$ of word length at most $k$, then $\rho'$ is Anosov. Moreover, for linear groups we explicitly bound $d(\rho(\gamma)p,\rho'(\gamma)p)$ in terms of the word length of $\gamma$, the Frobenius norms $ \abs{\rho(s)}_{Fr}$, and $\abs{\rho(s)-\rho'(s)}_{Fr}$, so we obtain a condition on $\rho'$ just in terms of the generators.

% summarize estimates
The bulk of the paper is devoted to a proof of Theorem \ref{main theorem}. We supply a number of estimates in Section \ref{section estimates} related to the geometry of the symmetric space $\X$. An important tool is the $\zeta$-angle $\angle_p^\zeta$, a $\Stab_G(p)$-invariant metric on $\Flagt$ introduced by Kapovich, Leeb and Porti in \cite{KLP14}, see Section \ref{section angles} for the definition. In Lemma \ref{zeta projection} we obtain explicit control on $\angle_p^\zeta(x,y)$ in terms of the Riemannian angle $\angle_p(x,y)$. The proof uses an explicit bound for the Hessian of a Morse function on $\Flagt$, see Proposition \ref{morse function} and Corollary \ref{morse function estimate}. A crucial step in the proof of the local-to-global principle is controlling the distance from the midpoint of a long  regular segment to a nearby diamond. The existence of such a bound is demonstrated in the proof of \cite[Proposition 7.16]{KLP14} via a limiting argument. To achieve explicit control, we consider the lengths of certain curves in $\X$ which are images of curves in $G$ under the orbit map, see Lemma \ref{right translated curves}. In Lemma \ref{cone rotation}, the curve in $G$ is required to lie in a maximal compact subgroup. In Lemma \ref{strongly asymptotic geodesics}, the curve is required to lie in a unipotent horocyclic subgroup. We combine these in Corollary \ref{midpoint projection} to obtain explicit, arbitrary control for the distance of midpoints to nearby Weyl cones (and hence diamonds). Kapovich, Leeb and Porti show that distance from a point $x \in \X$ to the parallel set $P(\tau_-,\tau_+)$ controls the $\zeta$-angle $\angle_x^\zeta(\tau_-,\tau_+)$ and vice versa via a compactness argument \cite[Section 2.4.5]{KLP14}. We give an explicit bound for $\angle_x^\zeta(\tau_-,\tau_+)$ in terms of $d(x,P(\tau_-,\tau_+))$ in Corollary \ref{distance to angles near pi}. This follows from Lemma \ref{simplex displacement}, whose proof relies on controlling the Lie derivative $\mathcal{L}_{X} \grad f_\tau$ where $X$ is a Killing vector field and $f_\tau$ is a Busemann function. Similarly, we obtain an explicit bound for $d(x,P(\tau_-,\tau_+))$ terms of $\angle_x^\zeta(\tau_-,\tau_+)$ in Lemma \ref{angle to distance} by controlling iterated derivatives of Busemann functions. In particular, we obtain an explicit uniform bound for the third derivative of the restriction of a Busemann function to a geodesic. 

% KLP zeta projection is implicit, follows from easy compactness argument
% KLP midpoint projection uses a limiting argument.
% KLP distance vs angle uses compactness arguments

%summarazie theorem 5.1
As in \cite{KLP14}, the proof of Theorem \ref{main theorem} is essentially broken into two parts, Theorem \ref{straight and spaced sequences are morse} and Theorem \ref{midpoint sequences}. Theorem \ref{straight and spaced sequences are morse} guarantees that a sequence $(x_n)$ with sufficiently spaced points forming $\zeta$-angles sufficiently close to $\pi$ is a Morse quasigeodesic. It is a quantified version of Theorem 7.2 in \cite{KLP14} and shares the same outline. One first shows that the property of ``moving away" from a simplex propagates along the sequence, see Section \ref{section straight and spaced sequences}. This implies that we can extract a simplex $\tau_-$ that the sequence $(x_n)$ moves away from (respectively towards) as $n$ increases (respectively decreases), and a simplex $\tau_+$ that the sequence $(x_n)$ moves away from (respectively towards) as $n$ decreases (respectively increases). One then verifies that the simplices $\tau_-,\tau_+$ are opposite and that the projections to the parallel set $P(\tau_-,\tau_+)$ define suitable diamonds, making $(x_n)$ a Morse quasigeodesic.

%summarize theorem 5.4
Theorem \ref{midpoint sequences} is a quantified version of Proposition 7.16 in \cite{KLP14}. It states that sufficiently spaced points on Morse quasigeodesics have straight and spaced midpoint sequences. A crucial ingredient is Corollary \ref{midpoint projection}, which allows us to force the midpoints to be arbitrarily close to the parallel sets in terms of the Morse and spacing parameters. This guarantees that they appear in nested Weyl cones, and makes the $\zeta$-angles arbitrarily straight. %Theorem \ref{midpoint sequences} allows us to apply Theorem \ref{straight and spaced sequences are morse} in the proof of Theorem \ref{local to global}.

%summarize local to global
Armed with Theorem \ref{straight and spaced sequences are morse} and Theorem \ref{midpoint sequences}, the proof of Theorem \ref{local to global} is similar to the proof of Theorem \ref{main theorem} given in \cite{KLP14}. We start with an $L$-local Morse quasigeodesic where $L$ is large enough to satisfy several explicit inequalities. We then replace our Morse quasigeodesic with a coarsification and take the midpoint sequence. Our assumptions together with Theorem \ref{midpoint sequences} shows that this coarse midpoint sequence is sufficiently straight and spaced, see Section \ref{straight and spaced sequences are morse}. An application of Theorem \ref{straight and spaced sequences are morse} shows that the midpoint sequence is a Morse quasigeodesic, and since it is a coarse approximation of the original sequence, the original sequence is also a Morse quasigeodesic, completing the proof. 

The usual proof of the local-to-global principle in hyperbolic geometry depends on the classical Morse Lemma. A higher rank version of the Morse Lemma was proved by Kapovich, Leeb and Porti in \cite{KLP18b}. In particular they prove that the orbit map $\Gamma \to \X$ of a finitely generated group is a coarsely uniformly regular quasiisometric embedding if and only if $\Gamma$ is word hyperbolic and the orbit map is a Morse quasiisometric embedding. It would be interesting to quantify their higher rank Morse Lemma by producing an explicit Morse parameter for (coarsely) uniformly regular quasiisometric embeddings, but we do not do this here. In the special case of the symmetric space associated to $\SL(d,\R)$, another proof of the higher rank Morse Lemma appears in \cite{BPS19}. There, Bochi, Potrie and Sambarino give yet another characterization of Anosov representations in terms of cone-types and dominated splittings.

The organization of the paper is as follows. In Section \ref{section notation} we fix some notation we use throughout the paper. In Section \ref{section background} we review some background of symmetric spaces. Much of this section is classical and may be skipped by experts on symmetric spaces, but we point the reader to our definition of regularity in Definition \ref{regularity} and the definition of $\zeta$-angle in Definition \ref{zeta angle}. The notion of regularity here is slightly different, but equivalent to, that in \cite{KLP14}, see Proposition \ref{klp regularity}. The bulk of the work is in Section \ref{section estimates} where we give several estimates related to the geometry of symmetric spaces. In Section \ref{local to global section}, we supplement the proof of the local-to-global principle in \cite{KLP14} with our estimates from Section \ref{section estimates}, reproving Theorem \ref{main theorem} with explicit bounds. Together with some standard geometric group theory, elementary hyperbolic geometry, and linear algebra in Section \ref{section applications}, this allows us to prove Theorems \ref{main theorem demonstration 1} and \ref{main theorem demonstration 2}. 

\subsection*{Acknowledgements}

I would like to thank my advisor Jeff Danciger for his frequent support and for suggesting this project, as well as Joan Porti for his encouragement and Misha Kapovich for some helpful advice. I also thank Florian Stecker for many interesting conversations about symmetric spaces and Anosov representations. Finally, I thank the UT math department and especially Martin Bobb, Teddy Weisman, Ne\v{z}a \v{Z}ager Korenjak, Casandra Monroe, and Charlie Reid for fostering a healthy atmosphere for math research.

\tableofcontents

	\section{Notation}\label{section notation}
	
	We establish our notational conventions in this paper. When possible, we have tried to keep notation consistent with \cite{KLP14,KLP17,E96}.
	
	\begin{enumerate}
		\item $\X=G/K$ will denote a symmetric space of noncompact type. $G$ is assumed to be the connected component of the isometry group of $\X$, and $K$ is a maximal compact subgroup of $G$, see Section \ref{section background}.
		\item We let $p,q,r,c$ denote points or curves in $\X$. We let $g,h,u,a$ denote elements or curves in $G$. An element or curve in $K$ may be denoted by $k$. 
		\item The Lie algebra of $G$ is denoted $\mfg$. The Lie algebra of $K$ is denoted $\mfk$. When a point $p$ is given, $K$ is the stabilizer of $p$ in $G$. Usually $U,V,W,X,Y,Z$ will denote elements of $\mfg$.
		\item The orbit map $\orb_p \colon G \to \X$, given by  $\orb_p(g)=gp$, has differential $\evp \colon \mfg \to T_p \X$ at the identity, see Section \ref{section background}. 
		\item The Cartan decomposition induced by $p \in \X$ is $\mfg = \mfk \oplus \mfp$. It corresponds to a Cartan involution $\vartheta_p \colon \mfg \to \mfg$, see Section \ref{cartan decomposition}.
		\item The Killing form on $\mfg$ is denoted $B$. Each point $p \in \X$ induces an inner product $B_p$ on $\mfg$ defined by $B_p(X,Y)=-B(\vartheta_pX,Y)$, see Section \ref{cartan decomposition}.
		\item We assume that the Riemmanian metric $\langle\cdot ,\cdot \rangle$ on $\X$ is the one induced by the Killing form, see Equation \ref{metric is Killing form}.
		\item The sectional curvature $\kappa$ of $\X$ has image $[-\kappa_0^2,0]$, see Section \ref{curvature}.
		\item A maximal abelian subspace of $\mfp$ will be denoted $\mfa$. The associated restricted roots are denoted $\Lambda \subset \mfa^\ast$. A choice of simple roots is denoted $\Delta$, see Section \ref{restricted root space decomposition}.
		\item Each maximal abelian subspace $\mfa$ has an action by the Weyl group and decomposition into Euclidean Weyl chambers denoted $V$, see Section \ref{section Weyl group}. 
		\item There is a \textit{vector-valued distance function} $\vec{d}\colon \X \times \X \to V_{mod}$ with image the model Euclidean Weyl chamber, see Equation \ref{vector valued distance}. In \cite{KLP14,KLP17} this map is denoted $\Delta$, and they let $\Delta$ denote the model Euclidean Weyl chamber we call $V_{mod}$. In this paper, $\Delta$ denotes a choice of simple roots.
		\item A spherical Weyl chamber $\sigma$ corresponds to a set of simple roots $\Delta$. For a face $\tau$ of $\sigma$ we have 
		$$ \Delta_\tau = \{ \alpha \in \Delta \mid \alpha(\tau) = 0 \}, \quad \Delta_\tau^+ = \{ \alpha \in \Delta \mid \alpha(\interior \tau)>0 \} ,$$
		see Equations \ref{delta tau}. We have 
		$$ \tau = \sigma \cap \bigcap_{\alpha \in \Delta_\tau} \ker \alpha, \quad \interior_\tau \sigma = \{ X \in \sigma \mid \forall \alpha \in \Delta_\tau^+, \alpha(X) >0 \}, \quad \partial_\tau \sigma = \sigma \cap \bigcup_{\alpha \in \Delta_\tau^+} \ker \alpha .$$
		\item The visual boundary of $\X$ is denoted $\partial \X$, see Section \ref{visual boundary}. We let $\tau, \sigma$ denote a spherical simplex/chamber in $\mfa$ or an ideal simplex/chamber in $\partial \X$. 
		\item There is a \textit{type projection} 
		$ \theta \colon \partial \X \to \sigma_{mod}$
		with image the model ideal Weyl chamber, see Section \ref{visual boundary}.
		\item A face of $\sigma_{mod}$ is called a \textit{model simplex} and denoted $\tau_{mod}$. There is a decomposition $\sigma_{mod}= \interior_{\tau_{mod}} \sigma_{mod} \sqcup \partial_{\tau_{mod}} \sigma_{mod}$, see Section \ref{section regularity ideal}.
		\item We define $(\alpha_0,\tau)$-regular and $(\alpha_0,\tau)$-spanning vectors and geodesics in Section \ref{regularity}. We extend that definition to ideal points in Section \ref{section regularity ideal}.
		\item We write $\Flagt$ for the set of ideal simplices in $\partial \X$ of type $\tau_{mod}$, see Section \ref{visual boundary}. 
		\item We define Weyl cones $V(x,\st(\tau),\alpha_0),V(x,\ost(\tau))$ and Weyl sectors $V(x,\tau)$ in Section \ref{section regularity ideal}.
		\item We describe the generalized Iwasawa decomposition $G=N_\tau A_\tau K$ in Section \ref{iwasawa decomposition}.
		\item A parallel set is denoted $P(\tau_-,\tau_+)$ for opposite simplices $\tau_-,\tau_+ \in \Flagt$. A horocycle is denoted $H(p,\tau)$, see Section \ref{section parallel sets and horocycles}. A diamond is denoted $\diamondsuit(p,q)$ and a truncated diamond is denoted $\diamondsuit_{\alpha_0}(p,q)$, see Section \ref{section straight and spaced sequences}.
		\item For $p \in \X$ and $x,y \in \overline{X} \setminus \{p\}$, $\angle_p(x,y)$ denotes the Riemannian angle at $p$ between $x$ and $y$. For $\eta,\eta' \in \partial \X$, $\angle_{Tits}(\eta,\eta')$ denotes their Tits angle. If $px$ and $py$ are $\tau_{mod}$-regular and $\tau,\tau' \in \Flagt$ then $\angle_p^\zeta(\tau,\tau'),\angle_p^\zeta(\tau,y),\angle_p(\zeta(\tau),\zeta(py))$ denote the $\zeta$-angles, see Section \ref{section angles}.
		\item A $(c_1,c_2,c_3,c_4)$-\textit{quasigeodesic} is a sequence $(x_n)$ (possibly finite, infinite, or biinfinite) in $\X$ such that 
		$$ \frac{1}{c_1} \abs{N}-c_2 \le d(x_n,x_{n+N}) \le \abs{N}c_3+c_4 .$$
		A quasigeodesic is $(\alpha_0,\tau_{mod},D)$-\textit{Morse} if for all $x_n,x_m$ there exists a diamond $\diamondsuit_{\alpha_0}(p,q)$ such that $d(p,x_n),d(q,x_m) \le D$ and for all $n \le i \le m$, $d(x_i,\diamondsuit) \le D$, see Section \ref{local to global section}.		
	\end{enumerate}
	
	\section{Background on symmetric spaces}\label{section background}
	
	We begin with some background on the structure of symmetric spaces of noncompact type. Experts on symmetric spaces can skip this section, but should note that we assume the metric is induced by the Killing form (see Equation \ref{metric is Killing form}), quantify the regularity of geodesics in Definition \ref{regularity}, and define the $\zeta$-angle in Definition \ref{zeta angle}. For detailed references on symmetric spaces see \cite{E96,H01,H79}. \\
	
	A \textit{symmetric space} is a connected Riemannian manifold $\X$ such that for each point $p \in \X$, there exists a \textit{geodesic symmetry} $S_p\colon \X \to \X$, an isometry fixing $p$ whose differential at $p$ is $(\dd{S_p})_p = -\id_{T_p\X}$. A symmetric space is necessarily complete with transitive isometry group. If $\X$ is nonpositively curved, it is simply connected. Simply connected Riemannian manifolds admit a de Rham decomposition into metric factors. If $\X$ is a nonpositively curved symmetric space with no Euclidean de Rham factors, we say $\X$ is a \textit{symmetric space of noncompact type}. Throughout the paper, $\X$ refers to any fixed symmetric space of noncompact type. 
	
	The isometry group of $\X$ is a semisimple Lie group, and we let $G$ be the identity component of the isometry group. For each point $p \in \X$, the stabilizer $K=G_p = \{ g \in G \mid  gp=p \}$ is a maximal compact subgroup of $G$. Hence $\X$ is diffeomorphic to $G/K$ by the orbit-stabilizer theorem for Lie groups and homogeneous spaces. We write $\mfg$ for the Lie algebra of left-invariant vector fields on $G$. 
	
	A \textit{Killing vector field} on a Riemannian manifold is vector field whose induced flow is by isometries. There is a natural linear isomorphism from $\mfg$ to the space of Killing vector fields on $\X$ by defining for $X \in \mfg$ the vector field $X^\ast$ given by 
	\begin{equation}\label{killing field} X^\ast_p \coloneqq \eval{ \dv{t} e^{tX} p }_{t=0} . \end{equation}
	The Lie bracket of two Killing vector fields is again a Killing vector field, but the map $X \mapsto X^\ast$ is a Lie algebra anti-homomorphism: $[X,Y]^\ast = -[X^\ast,Y^\ast]$. 	
	
	\subsection{Cartan decomposition}\label{cartan decomposition}
	
	Each point $p \in \X$ induces a \textit{Cartan decomposition} in the following way. The geodesic symmetry $S_p \colon \X \to \X$ induces an involution of $G$ by
	$$ g \mapsto S_p \circ g \circ S_p .$$
	The differential is a Lie algebra involution $\vartheta_p \colon  \mfg \to \mfg$, so we may write 
	$$ \mfg = \mfk \oplus \mfp $$
	where $\mfk = \{ X \in \mfg \mid \vartheta_pX=X \}$ and $\mfp  = \{ X \in \mfg \mid \vartheta_p X = -X \}$. Since $\vartheta_p$ preserves brackets, we have
	$$ [\mfk,\mfk] \subset \mfk, \quad [\mfk,\mfp] \subset \mfp, \quad [\mfp, \mfp] \subset \mfk .$$	
	We denote the orbit map $g \mapsto gp$ by $\orb_p\colon G \to \X$. The differential $(\dd{\orb_p})_1 \colon \mfg \to T_p \X$ has kernel precisely $\mfk$. Moreover, $\mfk$ is the Lie algebra of $K=G_p$. The restriction $(\dd{\orb_p})_1 \colon \mfp \to T_p\X$ is a vector space isomorphism. For any $X \in \mfg$, $(\dd{\orb_p})_1 X = X^\ast_p \eqqcolon \evp X$, see Equation \ref{killing field}, so we use the less cumbersome notation $\evp=(\dd{\orb_p})_1 \colon \mfg \to T_p\X$ throughout the paper (read as ``evaluation at $p$").
	
	Let $B$ denote the Killing form on $\mfg$ and let $\langle \cdot,\cdot \rangle$ denote the Riemannian metric on $\X$. We will assume that for all $X,Y \in \mfp$,
	\begin{equation}\label{metric is Killing form}
	B(X,Y) = \langle \mathrm{ev}_p X ,\mathrm{ev}_pY \rangle_p ,
	\end{equation}
	i.e.\ that the Riemannian metric on $\X$ is induced by the Killing form. Any other $G$-invariant Riemannian metrics on $\X$ only differs from this one by scaling by a global constant on each de Rham factor of $\X$.
	
	Under the identification of $\mfp$ with $T_p\X$, the Riemannian exponential map $\mfp \to \X$ is given by $X \mapsto e^X p$. In particular, the constant speed geodesics at $p$ are given by $c(t)=e^{tX}p$ for $X \in \mfp$. 
	
	The point $p \in \X$ induces an inner product $B_p$ on $\mfg$ defined by 
	\begin{equation}\label{Bp definition}
	B_p(X,Y)\coloneqq -B(\vartheta_p X,Y) .
	\end{equation}
	On $\mfp$, $B_p$ is just the restriction of the Killing form $B$, and we have required that the identification of $(\mfp,B)$ with $(T_p\X, \langle,\rangle)$ is an isometry. On $\mfk$, $B_p$ is the negative of the restriction of $B$ to $\mfk$. Since $\mfk$ and $\mfp$ are $B$-orthogonal, it follows that $B_p$ is an inner product on $\mfg$. For each $X \in \mfp$, $\ad X$ is symmetric with respect to $B_p$ on $\mfg$, and likewise for each $Y \in \mfk$, $\ad Y$ is skew-symmetric. 
	
	\subsection{Restricted root space decomposition}\label{restricted root space decomposition}
	
	Let $\mfa$ be a maximal abelian subspace of $\mfp$. Via the adjoint action, $\mfa$ is a commuting vector space of diagonalizable linear transformations on $\mfg$. Therefore $\mfg$ admits a common diagonalization called the \textit{restricted root space decomposition}. For each $\alpha \in \mfa^\ast$, define 
	$$ \mfg_\alpha = \{ X \in \mfg \mid \forall A \in \mfa, \ad A (X) = \alpha(A) X \}. $$
	We obtain a collection of \textit{roots} 
	$$ \Lambda = \{ \alpha \in \mfa^\ast \setminus \{0\} \mid \mfg_\alpha \ne 0 \} $$
	corresponding to the nonzero root spaces. The restricted root space decomposition is then 
	$$ \mfg = \mfg_0 \oplus \bigoplus_{\alpha \in \Lambda} \mfg_\alpha .$$
	
	%%% An older version included this footnote.
	%\footnote{Any finite-dimensional vector space $\mfa$ of commuting diagonalizable linear transformations admits a ``common eigenspace decomposition." The ``eigenvalues" are called roots and are naturally a subset of $\mfa^\ast$ and the ``eigenspaces" are called root spaces.}
	
	For each root $\alpha \in \Lambda$, define the coroot $H_\alpha \in \mfa$ by $\alpha(A)=B(H_\alpha,A)$ for all $A \in \mfa$. This induces an inner product, also denoted $B$, on $\mfa^\ast$ by defining $B(\alpha,\beta)\coloneqq B(H_\alpha,H_\beta)$. The set $\Lambda$ forms a root system\footnote{Note that this definition of root system is slightly different from the definition that appears in the study of, say, complex semisimple Lie algebras. There, one assumes that the only multiples of a root $\alpha$ appearing in $\Lambda$ are $\pm \alpha$. This assumption does not hold for restricted roots of symmetric spaces in general; for example, it fails in complex hyperbolic space.} in $(\mfa^\ast,B)$, see \cite[Proposition 2.9.3]{E96}. The restricted root space decomposition is $B_p$-orthogonal. A subset $\Lambda^+$ of the roots is \textit{positive} if for every $\alpha \in \Lambda$, exactly one of $\alpha,-\alpha$ is contained in $\Lambda^+$ and for any $\alpha,\beta \in \Lambda^+$ such that $\alpha +\beta$ is a root, we have $\alpha+\beta \in \Lambda^+$.

%%% An older version included this proposition. It doesn't seem necessary.	
	
%	\begin{prop}[{\cite[Proposition 2.9.3]{E96}}] 
%		$\Lambda$ is a \textit{root system}\footnote{Note that this definition of root system is slightly different from the definition that appears in the study of, say, complex semisimple Lie algebras. There, one assumes that the only multiple of a root $\alpha$ which appears in $\Lambda$ is $\pm \alpha$. This assumption does not hold for restricted roots of symmetric spaces, for example it fails in the symmetric space associated to $SU(2,1)$.} in $(\mfa^\ast,B)$. That is,
%		\begin{enumerate}
%			\item The span of $\Lambda$ is $\mfa^\ast$;
%			\item If $\alpha \in \Lambda$ and a scalar multiple $\lambda \alpha \in \Lambda$, then $\lambda \in \{ \pm \frac{1}{2}, \pm 1, \pm 2 \}$; 
%			\item For each $\alpha \in \Lambda$, reflection in $\alpha^\perp$ permutes $\Lambda$;
%			\item If $\alpha ,\beta \in \Lambda$, then 
%			$ 2 \frac{ B(\alpha,\beta)}{B(\alpha,\alpha)} $
%			is an integer. 
%		\end{enumerate}
%	In addition, the restricted root space decomposition is $B_p$-orthogonal.
%	\end{prop}
	
	The Cartan involution restricts to an isomorphism $\vartheta_p\colon  \mfg_\alpha \to \mfg_{-\alpha} $ for each $\alpha \in \Lambda \cup \{0\}$. Thus we have 
	$$ \mathfrak{p}_\alpha\coloneqq\mathfrak{p} \cap \mathfrak{g}_\alpha \oplus \mathfrak{g}_{-\alpha} = (\id - \vartheta_p)\mathfrak{g}_\alpha = (\id -\vartheta_p)\mathfrak{g}_{-\alpha} .$$
	and 
	$$ \mathfrak{k}_\alpha\coloneqq\mathfrak{k} \cap \mathfrak{g}_\alpha \oplus \mathfrak{g}_{-\alpha} = (\id + \vartheta_p)\mathfrak{g}_\alpha = (\id +\vartheta_p)\mathfrak{g}_{-\alpha} .$$
	Note that $\mathfrak{p}_\alpha = \mathfrak{p}_{-\alpha}$ and likewise  $\mathfrak{k}_\alpha = \mathfrak{k}_{-\alpha}$, so for $\Lambda^+$ a set of positive roots, we have the decomposition
	$$\mathfrak{g} = \mathfrak{a} \oplus \mathfrak{k}_0 \oplus \bigoplus_{\alpha \in \Lambda^+} \mathfrak{p}_\alpha \oplus \bigoplus_{\alpha \in \Lambda^+} \mathfrak{k}_\alpha $$
	which is both $B_p$ orthogonal and $B$ orthogonal. Some authors use the notation $\mathfrak{m}=\mfk_0$. 
	
	\subsection{Curvature and copies of $\HH^2$}\label{curvature}
	
	The curvature tensor $R$ of $\X$ may be defined using the Levi-Civita connection $\nabla$ by
	$$ R(u,v) = \nabla_u \nabla_v - \nabla_v \nabla_u - \nabla_{[u,v]}, $$
	for vector fields $u,v$ on $\X$. In a symmetric space there is a particularly nice formula for the curvature tensor. Our convention is that the sectional curvature spanned by orthonormal unit vectors $u,v$ is 
	$$ \kappa(u \wedge v) = \langle R(u,v)v,u \rangle. $$
	\begin{thm}[{\cite[p242]{Pet06}}]
		Let $X,Y,Z \in \mfp$ and write $X^\ast,Y^\ast,Z^\ast$ for the corresponding Killing vector fields on $\X$. Then 
		$$ (R(X^\ast,Y^\ast)Z^\ast)_p = -\evp[[X,Y],Z] .$$ 
	\end{thm}
	The theorem allows us to work directly with the sectional curvature by using the structure of the Lie algebra. Let $X \in \mfa$, $Y \in \mathfrak{p}$ and assume $X,Y$ are orthogonal unit vectors. For any $Y \in \mathfrak{p}$, we may write $Y=Y_0+\sum_{\alpha \in \Lambda^+ } Y_\alpha$ where $Y_0 \in \mfa$ and each $Y_\alpha \in \mathfrak{p}_\alpha$, and recall that this decomposition is $B$-orthogonal, so we have the lower curvature bound
	\begin{align*} & \kappa(X^\ast_p \wedge Y^\ast_p) = B(-[[X,Y],Y],X) = B([X,Y],[X,Y]) = -B([X,[X,Y]],Y) = \\ & -\sum_{\alpha \in \Lambda^+} B(\alpha(X)^2Y_\alpha,Y) = -\sum_{\alpha,\beta \in \Lambda^+} \alpha(X)^2 B(Y_\alpha,Y_\beta) = -\sum_{\alpha \in \Lambda^+} \alpha(X)^2 B(Y_\alpha,Y_\alpha) \ge - \kappa_0^2 
	\end{align*}
	where $\kappa_0$ is defined to be the maximum of $\{\alpha(X)\mid \alpha \in \Lambda, X \in \mfa, \abs{X}=1 \}$. In general, we have $\kappa_0 \le 1$, as we now explain. Since $\alpha(X)$ is maximized in the direction of the coroot $H_\alpha$, we have
	$$ \kappa_0 = \alpha\left(\frac{H_\alpha}{\abs{ H_\alpha}}\right) = \lvert H_\alpha \rvert $$
	for some $\alpha$. By \cite[2.14.5]{E96}, we have for $A,A' \in \mfa$ that $B(A,A') = \sum_{\beta \in \Lambda} (\dim \mfg_\beta) \beta(A)\beta(A') $, so
	$$  1= B\left(\frac{H_\alpha}{\abs{ H_\alpha}},\frac{H_\alpha}{\abs{ H_\alpha}}\right)  = \sum_{\beta \in \Lambda} \left( \dim \mathfrak{g}_\beta \right)\beta\left(\frac{H_\alpha}{\abs{ H_\alpha}}\right)^2  \ge \alpha\left(\frac{H_\alpha}{\abs{ H_\alpha}}\right)^2 =\kappa_0^2 .$$ 
	In particular, under this normalization where the symmetric space inherits its metric from the Killing form, the sectional curvature is always bounded between $0$ and $-1$. 
	
	\begin{exmp}
		In $\mathfrak{sl}(d,\R)$, each root $\alpha$ has $\abs{H_\alpha}=\frac{1}{\sqrt{d}}$, so we have $\kappa_0 = \frac{1}{\sqrt{d}}$ and the associated symmetric space has lower curvature bound $-\frac{1}{d}$. 
	\end{exmp}
	
	In Section \ref{section applications} we will need to know the curvature of copies of the hyperbolic plane in $\X$. These correspond to copies of $\sl2r$ in $\mfg$. Let $\alpha \in \Lambda$ and $X_\alpha \in \mfg_\alpha$ such that $B_p(X_\alpha,X_\alpha)=\frac{2}{\abs{H_\alpha}^2}$. Set $\tau_\alpha \coloneqq \frac{2}{\abs{H_\alpha}^2}H_\alpha$ so that $\alpha(\tau_\alpha)=2$. Set $Y_\alpha \coloneqq -\vartheta_pX_\alpha \in \mfg_{-\alpha}$. Then 
	$$ [\tau_\alpha,X_\alpha] = 2 X_\alpha, \quad [\tau_\alpha,Y_\alpha] = -2Y_\alpha, \quad \text{and} \quad [X_\alpha,Y_\alpha]=\tau_\alpha, $$
	where the last equality follows from considering $B([X_\alpha,Y_\alpha],A)$ for $A \in \mfa = \R H_\alpha \oplus \ker \alpha$. Then $\vartheta_p (X_\alpha + Y_\alpha) = \vartheta_p X_\alpha - \vartheta_p^2 X_\alpha = -(Y_\alpha + X_\alpha)$, so $X_\alpha + Y_\alpha \in \mfp$ and $\abs{X_\alpha + Y_\alpha}^2 = \pnorm{X_\alpha}^2+ \pnorm{Y_\alpha}^2 = \frac{4}{\abs{H_\alpha}^2}$. So $\frac{\abs{H_\alpha}}{2}(X_\alpha+Y_\alpha)$ and $\frac{H_\alpha}{\abs{H_\alpha}}$ are orthonormal unit vectors in $\mfp$, and 
	$$ \kappa \left( \frac{\abs{H_\alpha}}{2}(X_\alpha+Y_\alpha) \wedge \frac{H_\alpha}{\abs{H_\alpha}} \right) = - \alpha \left( \frac{H_\alpha}{\abs{H_\alpha}}\right)^2 \abs{ \frac{\abs{H_\alpha}}{2}(X_\alpha+Y_\alpha)}^2 = -\frac{\abs{H_\alpha}^4}{\abs{H_\alpha}^2} \frac{\abs{H_\alpha}^2}{4} \frac{4}{\abs{H_\alpha}^2} = -\abs{H_\alpha}^2  $$
	by the formula above.
	
	\begin{exmp}
		In the symmetric space associated to $\mathfrak{sl}(d,\R)$, the root spaces $\mfg_\alpha$ are one-dimensional, so the subalgebra $\mathfrak{sl}(2,\R)_\alpha$ spanned by $X_\alpha,Y_\alpha,\tau_\alpha$ is uniquely determined by $\alpha$ and we denote it by $\sl2r_\alpha$. The image of $\R H_\alpha \oplus \mfp_\alpha$ under the Riemannian exponential map at $p$ is a totally geodesic submanifold $\HH^2_\alpha$ isometric to the hyperbolic plane of curvature $-\frac{1}{d}$. 
	\end{exmp}
	
\subsection{Weyl chambers and the Weyl group}\label{section Weyl group}

In this section we describe Weyl faces as subsets of maximal abelian subspaces $\mathfrak{a} \subset \mathfrak{p}$. In Section \ref{visual boundary} we will define Weyl faces as subsets of the visual boundary $\partial \X$, and explain how the definitions relate. 

Let $\Lambda$ be the roots of a restricted root space decomposition of a maximal abelian subspace $\mathfrak{a}$ of $\mathfrak{p}$. For each $\alpha \in \Lambda \subset \mathfrak{a}^\ast$, the kernel of $\alpha$ is called a \textit{wall}, and a component $C$ of the complement of the union of the walls is called an \textit{open Euclidean Weyl chamber}; $C$ is open in $\mathfrak{a}$. A vector $X \in \mathfrak{a}$ is called \textit{regular} if it lies in an open Euclidean Weyl chamber and \textit{singular} otherwise. The closure $V$ of an open Euclidean Weyl chamber is a \textit{closed Euclidean Weyl chamber}; $V$ is closed in $\mathfrak{p}$.

For a closed Weyl chamber $V$ there is an associated set of \textit{positive roots} 
$$ \Lambda_+ \coloneqq \{ \alpha \in \Lambda \mid \forall v \in V, \alpha(v) \ge 0 \} $$
and \textit{simple roots} $\Delta$, i.e.\ those which cannot be written as a sum of two elements of $\Lambda_+$, see \cite[2.9.6]{E96}. 

We may define 
$$ N_K(\mfa)\coloneqq\{ k \in K \mid \Ad(k)(\mfa)=\mfa \},\quad Z_K(\mfa)\coloneqq\{k \in K \mid \forall A \in \mfa, \Ad(k)(A)=A \}.$$
Since the adjoint action preserves the Killing form, $N_K(\mfa)$ acts by isometries on $\mfa$ with kernel $Z_K(\mfa)$. We call the image of this action the \textit{Weyl group}. For each reflection $r_\alpha$ in a wall, it is possible to find a $k \in K$ whose action on $\mfa$ agrees with $r_\alpha$ \cite[2.9.7]{E96}. It is well-known that the Weyl group acts simply transitively on the set of Weyl chambers, which implies it is generated by the reflections in the walls of a chosen Weyl chamber. It is convenient for us to show this fact in Proposition \ref{morse function}, since the same techniques provide Corollary \ref{morse function estimate}.

\begin{figure}[h]	
	\begin{subfigure}{0.5\textwidth}
		\centering
		\begin{tikzpicture}[scale=0.7]
		\draw (-4,0) -- (4,0);
		\draw (-2,-3.46) -- (2,3.46);
		\draw (-2,3.46) -- (2,-3.46);
		\end{tikzpicture}
		\caption{The walls of a maximal flat in $\SL(3,\R)/\SO(3)$.}
	\end{subfigure}
	\begin{subfigure}{0.5\textwidth}
		\centering
		\begin{tikzpicture}
			\coordinate (fbl) at (0,0);
			\coordinate (fbr) at (4,0);
			\coordinate (ful) at (0,4);
			\coordinate (fur) at (4,4);
			\coordinate (bbl) at (2,1);
			\coordinate (bbr) at (6,1);
			\coordinate (bul) at (2,5);
			\coordinate (bur) at (6,5);
			\draw (fbl) -- (fbr) -- (fur) -- (ful) -- (fbl);
			\draw (ful) -- (bul) -- (bur) -- (bbr) -- (fbr);
			\draw (fur) -- (bur);
			\draw[dashed] (bul) -- (bbl) -- (bbr);
			\draw[dashed] (fbl) -- (bbl);
			\draw (ful) -- (bur) -- (fbr) -- (ful);
			\draw (fbl) -- (fur) -- (bul);
			\draw (fur) -- (bbr);
			\draw[dashed] (fbl) -- (bbr) -- (bul) -- (fbl);
			\draw[dashed] (fbl) -- (bul);
			\draw[dashed] (ful) -- (bbl) --(bur);
			\draw[dashed] (bbl) -- (fbr);
		\end{tikzpicture}
		\caption{The walls of a maximal flat in $\SL(4,\R)/\SO(4)$.}
	\end{subfigure}
\end{figure}

The Riemannian exponential map identifies maximal abelian subspaces in $\mathfrak{p}$ isometrically with maximal flats through $p$. So we can also refer to open/closed Euclidean Weyl chambers in $\X$ as the images of those in some $\mathfrak{a}$ under this identification. For every $X \in \mathfrak{p}$, there exists a maximal abelian subspace $\mathfrak{a}$ containing $X$, and in $\mathfrak{a}$, there exists some closed Euclidean Weyl chamber $V$ containing $X$. 

\subsection{A Morse function on flag manifolds}

In this subsection, we show that the vector-valued distance function $\vec{d}$ on $\X$ (denoted $d_\Delta$ in \cite{KLP14,KLP17}, see Definition \ref{vector valued distance}) is well-defined, and give part of a proof of Theorem \ref{Vmod welldefined}, an important part of the structure theory of symmetric spaces. Along the way we prove the $\vec{d}$-triangle inequality \cite{KLP14,KLP17,KLM09,Par}, and provide an estimate on the Hessian of a certain Morse function defined on flag manifolds embedded in $\mfp$, see Proposition \ref{morse function} and Corollary \ref{morse function estimate}.

We will use the following proposition. For $A \in \mfp$, let $\mfe_A$ be the intersection of all maximal abelian subspaces containing $A$. 
\begin{prop}[{\cite[2.20.18]{E96}}]\label{type uniqueness}
	Let $p$ in $\X$ with Cartan decomposition $\mfg = \mfk \oplus \mfp$ and let $k\in K$ and $A \in \mfp$. If $\Ad(k)(A)=A$ then for all $E \in \mathfrak{e}_A$ we have $\Ad(k)(E)=E$. 
\end{prop}
Note that there is a typo in Eberlein: the word ``maximal" is omitted in the definition of $\mfe_A$. The proof of Proposition \ref{type uniqueness} relies on passing to the compact real form of $\mfg^\mathbb{C}$. 

In this section, a \textit{flag manifold} is the orbit of a vector $Z \in \mfp$ under the adjoint action of $K=\Stab_G(p)$. The following proposition is essentially a standard part of the theory of symmetric spaces, however we will need to extract a specific estimate, recorded in Corollary \ref{morse function estimate}, in order to prove Lemma \ref{zeta projection}. 

\begin{prop}[{Cf. \cite[Lemma 6.3 p211]{H01} and \cite[Proposition 24]{E85}}]\label{morse function}
	Let $X,Z \in \mathfrak{p}$ be unit vectors. Define
	$$ f\colon  K \to \mathbb{R}, \quad f(k)\coloneqq B(X,\Ad(k)Z).$$
	\begin{enumerate}
		\item If $k$ is a critical point for $f$, then $\Ad(k)Z$ commutes with $X$.
		\item If $k$ is a local maximum for $f$, then $\Ad(k)Z$ lies in a common closed Weyl chamber with $X$.
		\item If $X$ is regular then the function $B(X,\cdot)\colon \Ad(K)Z \to \mathbb{R}$ is Morse and has a unique local maximum.
		\item If $X$ is regular then the distance function $d(X,\cdot)\colon \Ad(K)Z \to \mathbb{R}$ has a unique local minimum. 
	\end{enumerate}
\end{prop}

Note that $f$ is the composition of the orbit map $K \to \Ad(K)Z$ with the map $B(X,\cdot)\colon \Ad(K)Z \to \mathbb{R}$.

\begin{proof}
	1. Let $Y \in \mfk$, viewed as a left-invariant vector field on $K$. If $k$ is a critical point for $f$, then 
	\begin{align*}
	0 = & \d f_k(Y)= \eval{ \dv{t} f(ke^{tY}) }_{t=0} = \eval{ \dv{t} B(X,\Ad(ke^{tY})Z)}_{t=0} \\
	& = B(X,\Ad(k)(\ad(Y)(Z))) = B(X,[Y',Z'])=B([Z',X],Y') 
	\end{align*}   
	where we write $Y'=\Ad(k)Y$ and $Z'=\Ad(k)Z$. Since $Y'$ is an arbitrary element of $\mfk$, $[X,Z'] \in \mfk$, and $B$ is negative definite on $\mfk$, we can conclude that $[X,Z']=0$, which is the claim.
	
	2. At a critical point $k$ for $f$, the Hessian of $f$ at $k$ is a symmetric bilinear form on $T_kK$ determined by 
	$$ \Hess(f)(v,v)_k = (f \circ c)''(0) $$
	for any curve $c$ with $c(0)=k$ and $c'(0)=v$. Let $Y \in \mfk$, the left-invariant vector fields on $K$, and choose $c(t)=ke^{tY}$. To compute the Hessian of $f$ we only need to compute 
	\begin{align*}
	\derivtwo{t} f(ke^{tY}) \atzero{t} & = \deriv{t} B(X,\Ad(ke^{tY})(\ad(Y)(Z))) \atzero{t} = B(X,\Ad(k)([Y,[Y,Z]])) = B(X,[Y',[Y',Z']]) \\
	&= B([X,Y'],[Y',Z']) 
	= B([Z',[X,Y']],Y') = B(\ad(Z')\ad(X)(Y'),Y') = B(T Y',Y')
	\end{align*} 
	where we write $T=\ad(Z')\circ \ad(X)$ as a linear transformation on $\mfk$. At a critical point $X$ and $Z'$ commute by part 1, and we can choose a maximal abelian subspace $\mfa$ containing both of them, and then consider the corresponding restricted root space decomposition. For $Y_\alpha \in \mfk_\alpha$, 
	$$ T Y_\alpha = \alpha(Z') \alpha(X) Y_\alpha $$
	so the transformation $T$ has the eigenvalue $\alpha(Z')\alpha(X)$ on its eigenspace $\mfk_\alpha$ and acts as $0$ on $\mfk_0$. Since we assumed $k$ is a local maximum for $f$, we have
	$$ 0 \ge \eval{ \derivtwo{t} f(ke^{tY}) }_{t=0} = B(T Y',Y') $$
	for all $Y \in \mfk$, so for each $\alpha \in \Lambda$, $\alpha(Z')\alpha(X) \ge 0$, and therefore $X$ and $Z'$ lie in a common closed Weyl chamber. 
	
	\begin{figure}[h]
		\centering
		\begin{tikzpicture}
		\draw (-4,0) -- (4,0);
		\draw (-2,-3.46) -- (2,3.46);
		\draw (-2,3.46) -- (2,-3.46);
		\node[circle,scale=0.5,fill=black,label=below:$X$] (x) at (2.6,1.5) {};
		\node[circle,fill=black,scale=0.5,label=right:{$\Ad(k_1)Z$}] (z1) at (2,0.3) {};
		\node[circle,fill=black,scale=0.5,label=right:{$\Ad(k_2)Z$}] (z2) at (2,-0.3) {};
		\node[circle,fill=black,scale=0.5,label=left:{$\Ad(k_3)Z$}] (z3) at (-1.26,1.58) {};
		\node[circle,fill=black,scale=0.5,label=above:{$\quad \Ad(k_4)Z$}] (z4) at (-0.74,1.88) {};
		\node[circle,fill=black,scale=0.5,label=left:{$\Ad(k_5)Z$}] (z5) at (-1.26,-1.58) {};
		\node[circle,fill=black,scale=0.5,label=below:{$\quad \Ad(k_6)Z$}] (z6) at (-0.74,-1.88) {};
		\end{tikzpicture}
		\caption{The intersection $\Ad(K)Z \cap \mfa$}
	\end{figure}
	
	3. We may assume that $Z$ is a critical point of $f$ by precomposing $f$ with a left translation of $K$. The differential $(\dd{\orb_Z})_1\colon \mathfrak{k} \to T_z \Ad(K)Z$ is given by $-\ad Z$ and has kernel $\mfk_Z = Z_{\mfk}(Z) = \{ W \in \mfk \mid [W,Z]=0\}$ with orthogonal complement $\mfk^Z=\bigoplus_{\alpha \in \Lambda : \alpha(Z)>0}\mfk_\alpha $. Then $k$ is a critical point for $f$ if and only if $Z(k)=\Ad(k)Z$ is a critical point for $B(X,\cdot)$. The Hessians satisfy
	$$ \Hess(B(X,\cdot))((\dd{\orb_Z})_k U, (\dd{\orb_Z})_k V)_{\Ad(k)Z} = \Hess(f)(U,V)_k ,$$
	so by the calculation above the critical points are nondegenerate, occur at $\Ad(k)Z$ when $[\Ad(k)Z,X]=0$, and have index the number of positive signs in the collection $\alpha(X)\alpha(\Ad(k)Z)$, (weighted by $\dim \mfk_\alpha$) as $\alpha$ ranges over the roots with $\alpha(Z)>0$. These can only be nonnegative when $\Ad(k)Z$ lies in the closed Weyl chamber containing $X$.
	
	For uniqueness, observe that any two maximizers $Z',Z''$ lie in the closed Weyl chamber containing $X$, and suppose $\Ad(k)(Z')=Z''$. The adjoint action takes walls to walls so $\Ad(k)$ preserves the facet spanned by $Z',Z''$ and hence fixes its soul (i.e.\ its center of mass) \cite[p65]{E96}. By Proposition \ref{type uniqueness}, $\Ad(k)$ fixes each point of the face, and in particular $Z'=Z''$.
	
	4. Since $(\mfp,B)$ is a Euclidean space, 
	$$ d_{\mfp}(X,Y)^2 = B(X-Y,X-Y) = B(X,X)+B(Y,Y) - 2 B(X,Y) $$
	so if $X,Y$ are unit vectors in $\mfp$
	$$ d_{\mfp}(X,Y)^2 = 2(1-B(X,Y))$$
	and the distance function $d_{\mfp}(X,\cdot)$ is minimized when $B(X,\cdot)$ is maximized. Then by part 3, the distance function is uniquely minimized at the unique $\Ad(k)Z$ in the closed Weyl chamber containing $X$. 
\end{proof}

The next two results are part of the standard theory of symmetric spaces. Since we have already proven Proposition \ref{morse function} it is convenient to give the proofs.

\begin{cor}\cite[Section 2.12]{E96}\label{orbits and chambers}
	Every $K$-orbit in the unit sphere $S(\mfp)$ intersects each closed spherical Weyl chamber exactly once.
\end{cor}

\begin{proof}
	Let $X$ be a regular vector in a chosen Weyl chamber. The $K$-orbit of a unit vector $Z$ is compact and therefore the function $d_{\mfp}(X,\cdot)$ has a global minimum on $\Ad(K)Z$. But that function has a unique local minimum which must lie in the chosen closed Weyl chamber. 
\end{proof}

For a point $p \in \X$, maximal abelian subspace $\mfa \subset \mfp$ and closed Euclidean Weyl chamber $V \subset \mfa$, we call $(p,\mfa,V)$ a point-chamber triple.

\begin{thm}\cite[Section 2.12]{E96}\label{Vmod welldefined}
	For any two point-chamber triples $(p,\mathfrak{a},V),(p',\mathfrak{a}',V')$ there exists an isometry $g \in G$ taking $(p,\mathfrak{a},V)$ to $(p',\mathfrak{a}',V')$. If $g$ stabilizes $(p,\mfa,V)$, then it acts trivially on it.
\end{thm}

\begin{proof}
	The group $G$ acts transitively on $X$, so we may assume that $p'=p$ and then show that an element of $K=\Stab_G(p)$ takes $(\mathfrak{a},V)$ to $(\mathfrak{a}',V')$. Choose any regular unit vectors $X\in V$, $Z \in V'$. Then Proposition \ref{morse function} implies there is an element $k \in K$ such that $\Ad(k)Z$ is in the same open Weyl chamber as $X$. Regular vectors lie in unique Weyl chambers in unique maximal abelian subspaces, so $\Ad(k)\mathfrak{a}'=\mathfrak{a}$ and $\Ad(k)V'=V$.
	
	If $g$ fixes $p$ and stabilizes $(\mfa,V)$, then it acts trivially on $V$ by Corollary \ref{orbits and chambers}.
\end{proof}

The above isometry is not necessarily unique. For example, consider hyperbolic space $\mathbb{H}^n, n\ge 3$. There a Euclidean Weyl chamber is just a geodesic ray, which has infinite pointwise stabilizer. However the action on $V$ is unique. 

As a corollary, we may define the \textit{vector-valued distance function}
\begin{equation}\label{vector valued distance}
 \vec{d}\colon  \X \times \X \to (\X \times \X) /G \eqqcolon V_{mod} 
\end{equation}
to have range a model closed Euclidean Weyl chamber. One could think of $V_{mod}$ as some preferred Euclidean Weyl chamber, but it is better to think of it as an abstract Euclidean cone with no reference to a preferred basepoint, flat or Weyl chamber in $\X$. There is an ``opposition involution" $\iota\colon V_{mod} \to V_{mod}$ induced by any geodesic symmetry $S_p$. On a model pointed flat $\mfa_{mod}$, the composition of $-\id$ with the longest element of the Weyl group restricts to $\iota$ on the model positive chamber $V_{mod}$. Note that $\iota \vec{d}(p,q) = \vec{d}(q,p)$.

The triangle inequality implies that for any $p,p',q,q'$ in a metric space,
$$ \abs{d(p,q)-d(p',q')} \le d(p,p')+d(q,q').$$
The next result is the ``vector-valued triangle inequality" for symmetric spaces.
\begin{cor}[{The $\vec{d}$-triangle inequality \cite{KLP17,KLM09,Par}}]\label{vector valued triangle inequality}
	For points $p,p',q,q'$ in $\X$, 
	$$ \lvert \vec{d}(p,q)-\vec{d}(p',q') \rvert \le d(p,p')+d(q,q').$$
\end{cor}

\begin{proof}
	In a moment we will use the proposition to prove that for any $p,q,q'$ in $\X$, 
	\begin{equation}\label{vector valued reverse triangle inequality}
	 \lvert \vec{d}(p,q)-\vec{d}(p,q') \rvert \le d(q,q') ,
	 \end{equation}
	from which the general inequality follows easily:
	\begin{align*}
	\lvert \vec{d}(p,q)-\vec{d}(p',q') & = \lvert \vec{d}(p,q)-\vec{d}(p,q') +\vec{d}(p,q')-\vec{d}(p',q')\rvert \\
	& \le \lvert \vec{d}(p,q)-\vec{d}(p,q') \rvert +\lvert \iota \vec{d}(q',p)-\iota \vec{d}(q',p')\rvert \le d(q,q') + d(p,p'). 
	\end{align*}
	To prove \ref{vector valued reverse triangle inequality}, let $X,Z \in \mfp$ such that $e^Xp=q$ and $e^Zp=q'$. Choose a closed Weyl chamber $V$ containing $X$ and the unique $Z'$ in the $K$-orbit of $Z$ in that Weyl chamber. The map $\vec{d}(p,e^{(\cdot)}p)\colon V\to V_{mod}$ is an isometry. Note that $k\mapsto B(X,\Ad(k)Z)$ is maximized when $k\mapsto B(X,\Ad(k)Z)/\lvert X\rvert \lvert Z\rvert$ is maximized, so by Proposition \ref{morse function}
	\begin{align*}
	\lvert \vec{d}(p,q)-\vec{d}(p,q') \rvert^2 & = \lvert X-Z' \rvert^2 = \lvert X \rvert^2 + \lvert Z' \rvert^2 - 2 \langle X,Z' \rangle \\
	& \le \lvert X \rvert^2 +\lvert Z \rvert^2 - 2 \langle X,Z \rangle = d_{\mfp} (X,Z)^2 \le d(q,q')^2
	\end{align*}
	since the Riemannian exponential map is distance non-decreasing by the nonpositive curvature of $\X$. 
\end{proof}

\subsection{Regularity in maximal abelian subspaces}

A \textit{spherical Weyl chamber} is the intersection of a Euclidean Weyl chamber with the unit sphere $S$ in $\mathfrak{a}$. A spherical Weyl chamber $\sigma$ is a spherical simplex, and each of its faces $\tau$ is called a \textit{Weyl face}. Each Euclidean (resp.\ spherical) Weyl face is the intersection of walls of $\mathfrak{a}$ (resp.\ as well as $S$). The interior of a face $\interior(\tau)$ is obtained by removing its proper faces; the interiors of faces are called \textit{open simplices}. The unit sphere $S$ is a disjoint union of the open simplices. If $\tau$ is the smallest simplex containing a unit vector $X$ in its interior, we say that $\tau$ is \textit{spanned} by $X$ and $X$ is $\tau$-\textit{spanning}. 

We will quantify the regularity of tangent vectors using a parameter $\alpha_0>0$. We will show in Proposition \ref{klp regularity} that our definition of regularity is equivalent to the definition in \cite{KLP14}. A similar definition appears in \cite[Definition (2.6)]{KLP18b}.

\begin{defn}[Regularity]\label{regularity}
	Let $p \in \X$ and $\X$ be a closed spherical Weyl chamber and let $\tau$ be a face of $\sigma$. Consider the corresponding maximal abelian subspace $\mfa$ in $\mfp$ and set of simple roots $\Delta$. We define
	\begin{equation}\label{delta tau}
	\Delta_\tau = \{ \alpha \in \Delta \mid \alpha(\tau) = 0 \}, \quad \Delta_\tau^+ = \{ \alpha \in \Delta \mid \alpha(\interior \tau)>0 \}. 
	\end{equation}	
	A vector $X \in \mfa$ is called \textit{$(\alpha_0,\tau)$-regular} if for each $\alpha \in \Delta_\tau^+, \alpha(X) \ge \alpha_0 \abs{X}$. A geodesic $c$ at $p$ is called \textit{$(\alpha_0,\tau)$-regular} if $c'(0) = \evp X$ for an $(\alpha_0,\tau)$-regular vector $X \in \mfa$.
\end{defn}

It is immediate from the definition that $X$ is $(\alpha_0,\sigma)$-regular for some $\alpha_0 >0$ and $\sigma$ if and only if $X$ is regular. 
We define 
\begin{equation}\label{lambda tau}
\Lambda_\tau \coloneqq \{ \alpha \in \Lambda \mid \alpha(\tau) = 0 \}, \quad \Lambda_\tau^+ \coloneqq \{ \alpha \in \Lambda \mid \alpha(\interior \tau)>0 \} 
\end{equation}
Observe that $X$ is $(\alpha_0,\tau)$-regular if and only if for each root $\alpha \in \Lambda_\tau^+$ we have $\alpha(X)\ge\alpha_0$. 

\begin{rem}\label{roots vs distance}
	The signed distance from a vector $A \in \mfa$ to the wall $\ker \alpha$ is $\alpha(A)/\abs{\alpha} \ge \alpha(A)/\kappa_0$.
\end{rem}

\begin{defn}\label{spanning}
	A unit vector $X$ is $(\alpha_0,\tau)$-\textit{spanning} if it is $\tau$-spanning and $(\alpha_0,\tau)$-regular.
\end{defn}

We may now record a mild extension of Proposition \ref{morse function} which will appear in Lemma \ref{zeta projection}.

%corrected
\begin{cor}\label{morse function estimate}
	Suppose $X\in \mfp$ is an $(\alpha_0,\tau)$-regular unit vector and $Z \in \mfp$ is a $(\zeta_0,\tau)$-spanning unit vector. Then $Z$ is the unique maximum of $B(X,\cdot)\colon \Ad(K)Z \to \R$, and for all $U,V \in T_Z \Ad(K)Z$, 
	$$ \abs{\Hess(B(X,\cdot))(U,V)_Z} \ge \alpha_0\zeta_0 \abs{B_p(U,V)} .$$
\end{cor}

\begin{proof}
	The proof of Proposition \ref{morse function} goes through in this setting, requiring only the following observation: if $X$ is $\tau$-regular and lies in a spherical Weyl chamber $\sigma$, then $\tau$ is a face of $\sigma$. If $U,V \in T_Z \Ad(K)Z$ correspond to $U',V' \in \mfk^\tau$ under the identification $T_Z \Ad(K)Z = \mfk^\tau$, we showed that $\Hess(B(X,\cdot))(U,V)_Z = B( \ad(Z) \ad(X) U',V')$.
\end{proof}

\subsection{The visual boundary $\partial \X$}\label{visual boundary}

We say two unit speed geodesic rays $c_1, c_2$ are \textit{asymptotic} if there exists a constant $D>0$ such that
$$ d(c_1(t),c_2(t)) \le D $$
for all $t \ge 0$. The asymptote relation is an equivalence relation on unit-speed geodesic rays and the set of asymptote classes is called the \textit{visual boundary} of $\X$ and denoted by $\partial \X$. There is a natural topology on $\partial \X$ called the \textit{cone topology}, where for each point $p \in \X$ the map $S(T_p \X) \to \partial \X$ (which takes a unit tangent vector to the geodesic ray with that derivative) is a homeomorphism. In fact the cone topology extends to $\overline{\X}\coloneqq\X \cup \partial \X$, yielding a space homeomorphic to a unit ball of the same dimension as $\X$.

\begin{lem}\label{asymptotic tips}
	If $c_1$ and $c_2$ are asymptotic geodesic rays then for all $t \ge 0$,
	$$ d(c_1(t),c_2(t)) \le d(c_1(0),c_2(0)) .$$
\end{lem}

\begin{proof}
	The left hand side is convex \cite{E96} and bounded above, hence (weakly) decreasing. 
\end{proof}

We have a natural action of $G$ on $\partial \X$: $g[c]= [g \circ c]$. For $\eta \in \partial \X$, we denote the stabilizer 
$$ G_\eta \coloneqq \{ g \in G \mid g\eta =\eta \} $$
and call $G_\eta$ the \textit{parabolic subgroup} fixing $\eta$. (Note that in \cite{GW12} and \cite{GGKW17}, $G$ itself is a parabolic subgroup, but in this paper a parabolic subgroup is automatically a proper subgroup.) When $\eta$ is regular, $G_\eta$ is a \textit{minimal parabolic} subgroup of $G$ (sometimes called a Borel subgroup).

Let $\eta,\eta'$ be ideal points in $\partial \X$, represented by the geodesics $c(t)=e^{tX}p$ and $c'(t)=e^{tY}q$. Then since $G$ is transitive on point-chamber triples, we can find $g \in G$ such that $g q=p$ and $\Ad(g)Y$ lies in a (closed) Euclidean Weyl chamber in common with $X$. In particular, every $G$ orbit in $\partial \X$ intersects every spherical Weyl chamber exactly once. 

Each unit sphere $S(\mfp)$ has the structure of a simplical complex compatible with the action of $G$. By Theorem \ref{Vmod welldefined} this simplicial structure passes to $\partial \X$, which is in fact a thick spherical building whose apartments are the ideal boundaries of maximal flats. In \cite{KLP14,KLP17} the spherical building structure on $\partial \X$ is used to describe the regularity of geodesic rays. We have used the restricted roots to define regularity and will show the notions are equivalent in Proposition \ref{klp regularity}. When we need to distinguish between simplices in $S(\mfp)$ and simplices in $\partial \X$ we call the former \textit{spherical} and the latter \textit{ideal}. Compared to a spherical simplex, an ideal simplex lacks the data of a basepoint $p \in \X$. 

Define the \textit{type map} to be
$$ \theta\colon  \partial \X \to \partial \X /G \eqqcolon \sigma_{mod} $$
with range the model ideal Weyl chamber. The opposition involution $\iota\colon V_{mod} \to V_{mod}$ induces an opposition involution $\iota\colon  \sigma_{mod} \to \sigma_{mod}$, see the discussion after Equation \ref{vector valued distance} in the previous subsection. The faces of $\sigma_{mod}$ are called \textit{model simplices}. For a model simplex $\tau_{mod} \subset \sigma_{mod}$, we define the \textit{flag manifold} $\Flagt$ to be the set of simplices $\tau$ in $\partial \X$ such that $\theta(\tau)=\tau_{mod}$. If ideal points $\eta,\eta'$ span the same simplex $\tau$, then they correspond to the same parabolic subgroup, so we define $G_\tau\coloneqq G_\eta$. A model simplex corresponds to the conjugacy class of a parabolic subgroup of $G$. 
% it is now possible to talk about model simplices

\subsection{Regularity for ideal points}\label{section regularity ideal}

\begin{figure}[h]
	\centering
	\begin{subfigure}{0.3\textwidth}
		\centering
		\begin{tikzpicture}
			\node[circle,scale=0.5,fill=black,label=below:{$\tau_1$}] (t1) at (0,0) {};
			\node[circle,fill=black,scale=0.5,label=above:{$\tau_2$}] (t2) at (2,3.46) {};
			\node[circle,fill=black,scale=0.5,label=below:{$\tau_3$}] (t3) at (4,0) {};
			\draw (t1) -- node[pos=0.2](t12){} (t2) -- (t3) node[pos=0.8](t23){} -- (t1);
			\coordinate (sa) at (1.2,0.6);
			\coordinate (sb) at (2.8,0.6);
			\coordinate (sc) at (2,2.2);
			\draw[dashed,fill=gray!30] (sa) -- (sb) -- (sc) -- cycle;	
		\end{tikzpicture}
		\caption{$(\alpha_0,\sigma_{mod})$-regular}
	\end{subfigure}
	\begin{subfigure}{0.3\textwidth}
	\centering
		\begin{tikzpicture}
			\node[circle,scale=0.5,fill=black,label=below:{$\tau_1$}] (t1) at (0,0) {};
			\node[circle,fill=black,scale=0.5,label=above:{$\tau_2$}] (t2) at (2,3.46) {};
			\node[circle,fill=black,scale=0.5,label=below:{$\tau_3$}] (t3) at (4,0) {};
			\draw (t1) -- (t2) -- (t3)  -- (t1)  node[pos=0.2](t13a){} node[pos=0.8](t13b){};
			\coordinate (t13c) at (2,2.1);
			\fill[gray!30] (t13a.center) -- (t13c)--(t13b.center) -- cycle;		
			\draw[dashed] (t13a.center) -- (t13c)--(t13b.center);	
		\end{tikzpicture}
		\caption{$(\alpha_0,\tau_{13})$-regular}
	\end{subfigure}	
	\begin{subfigure}{0.3\textwidth}
	\centering
	\begin{tikzpicture}
		\node[circle,scale=0.5,fill=black,label=below:{$\tau_1$}] (t1) at (0,0) {};
		\coordinate (t2) at (2,3.46);
		\node[circle,fill=black,scale=0.5,label=below:{$\tau_3$}] (t3) at (4,0) {};
		\draw (t1) -- coordinate[pos=0.2](t12) (t2) -- (t3) coordinate[pos=0.8](t23) -- (t1);
		\fill[gray!30] (t2) -- (t12) -- (t23) -- cycle;	
		\draw[dashed] (t12) -- (t23);	
		\node[circle,fill=black,scale=0.5,label=above:{$\tau_2$}] at (t2) {};	
	\end{tikzpicture}
	\caption{$(\alpha_0,\tau_2)$-regular}
	\end{subfigure}
	\caption{$(\alpha_0,\tau_{mod})$-regularity for various choices of $\tau_{mod}$}
\end{figure}
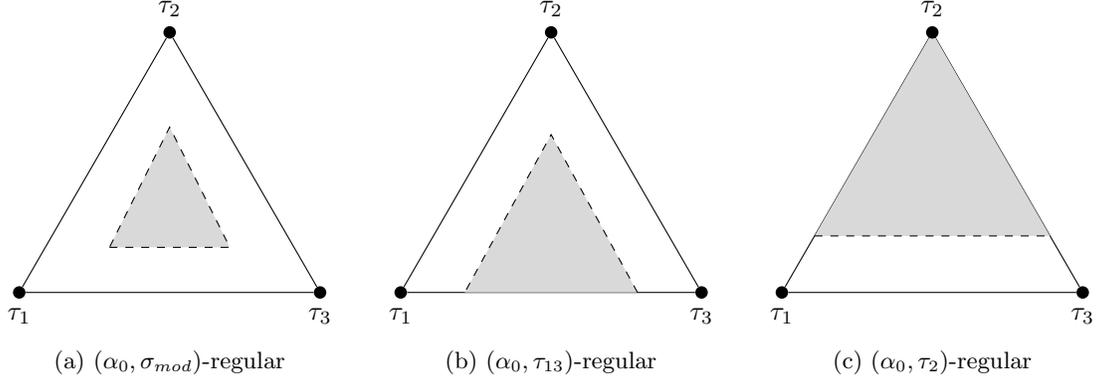

Theorem \ref{Vmod welldefined} implies that ``model roots" are well-defined: if $g \in G$ takes the point-chamber triple $(p,\mfa,V)$ to $(p',\mfa',V')$ and takes the simplex $\tau \subset \partial V$ to $\tau' \subset \partial V'$, it also takes $\Delta_\tau$ to ${\Delta'}_{\tau'}$ and $\Delta_\tau^+$ to ${\Delta'}_{\tau'}^+$, where $\Delta$ is the simple roots in $\mfa^\ast$ corresponding to $V$ and $\Delta'$ is the simple roots in $\mfa'$ corresponding to $V'$.

An ideal point $\eta \in \partial \X$ is called $(\alpha_0,\tau)$\textit{-regular} if every geodesic in its asymptote class is $(\alpha_0,\tau)$-regular. As soon as one representative of an ideal point is $(\alpha_0,\tau)$-regular, every representative is. A vector, geodesic, or ideal point is \textit{$(\alpha_0,\tau_{mod})$-regular} if it is $(\alpha_0,\tau)$-regular for some simplex $\tau$ of type $\tau_{mod}$.

The \textit{open star} of a simplex $\tau$, denoted $\ost(\tau)$, is the union of open simplices $\nu$ whose closures intersect $\tau$. Equivalently, it is the collection of $\tau$-regular points in $\partial X$. For a model simplex, $\interior_{\tau_{mod}}(\sigma_{mod})$ is the collection of $\tau_{mod}$-regular ideal points in $\sigma_{mod}$. Equivalently, it is $\sigma_{mod} \setminus \bigcup_{\alpha \in \Delta_\tau^+} \ker \alpha $.\footnote{ In \cite{KLP14} the notation $\ost(\tau_{mod})$ was used for what is called $\interior_{\tau_{mod}}(\sigma_{mod})$ here and in \cite{KLP17}.}  We have 
$$ \tau = \sigma \cap \bigcap_{\alpha \in \Delta_\tau} \ker \alpha, \quad \interior_\tau \sigma = \{ \eta \in \sigma \mid \forall \alpha \in \Delta_\tau^+, \alpha(\eta) >0 \}, \quad \partial_\tau \sigma = \sigma \cap \bigcup_{\alpha \in \Delta_\tau^+} \ker \alpha .$$ 
There is a decomposition $\sigma_{mod}= \interior_{\tau_{mod}} \sigma_{mod} \sqcup \partial_{\tau_{mod}} \sigma_{mod}$. 

We call the set of $(\alpha_0,\tau)$-regular points the ``$\alpha_0$-star of $\tau$." We define the closed cone on the $\alpha_0$-star of $\tau$
$$  V(p,\st(\tau),\alpha_0) \coloneqq \{ c_{px}(t) \mid t \in \left[ 0,\infty \right), x \text{ is }
(\alpha_0,\tau)\text{-regular} \} $$
the cone on the open star of $\tau$
$$  V(p,\ost(\tau)) \coloneqq \{ c_{px}(t) \mid t \in \left[ 0,\infty \right), x \text{ is }
\tau\text{-regular} \} $$
and the Euclidean Weyl sector 
$$ V(p,\tau) \coloneqq \{ c_{px}(t) \mid t \in \left[ 0,\infty \right), x \text{ is }
\tau\text{-spanning} \} .$$
It follows from Lemma \ref{asymptotic tips} that the Hausdorff distance between $V(p,\st(\tau),\alpha_0)$ and $V(q,\st(\tau),\alpha_0)$ is bounded above by $d(p,q)$, and the same holds for the open cones $V(p,\ost(\tau))$ and $V(q,\ost(\tau))$ and for the Weyl sectors $V(p,\tau), V(q,\tau)$. 

We now describe the notion of regularity used in \cite{KLP14,KLP17} and show it is equivalent to our definition. We always work with respect to a fixed type $\tau_{mod}$.  A subset $\Theta \subset \sigma_{mod}$ is called $\tau_{mod}$\textit{-Weyl convex} if its symmetrization $W_{\tau_{mod}} \Theta \subset a_{mod}$ is a convex subset of the model apartment $a_{mod}$. Here we think of the Weyl group $W$ as acting on the visual boundary $a_{mod}$ of a model flat $\mfa_{mod}$ with distinguished Weyl chamber $\sigma_{mod}$ and $W_{\tau_{mod}}$ is the subgroup of $W$ stabilizing the simplex $\tau_{mod}$. One then quantifies $\tau_{mod}$-regular ideal points by fixing an auxiliary compact $\tau_{mod}$-Weyl convex subset $\Theta$ of $\interior_{\tau_{mod}}(\sigma_{mod}) \subset \sigma_{mod}$.

An ideal point $\eta$ is $\Theta$\textit{-regular} if $\theta(\eta)\in \Theta$. It is easy to see that the notions of $\Theta$-regularity and $(\alpha_0,\tau_{mod})$-regularity are equivalent. 

\begin{prop}\label{klp regularity}
	Let $\Delta_{\tau_{mod}} \subset \Delta$ be the model simple roots corresponding to a simplex $\tau_{mod} \subset \sigma_{mod}$. Then
\begin{enumerate}
	\item If $\Theta$ is a compact subset of $\interior_{\tau_{mod}}(\sigma_{mod})$ then every $\Theta$-regular ideal point is $(\alpha_0,\tau_{mod})$-regular for $\alpha_0 = \min_{\alpha \in \Delta_{\tau_{mod}}^+} \alpha (\Theta)$.
	\item Every $(\alpha_0,\tau_{mod})$-regular ideal point is $\Theta$-regular for $\Theta= \{ \xi \in \sigma_{mod} \mid  \forall \alpha \in \Delta_{\tau_{mod}}^+, \alpha(\zeta) \ge \alpha_0   \} $.	
\end{enumerate}
\end{prop}

\begin{proof}
	We first prove $1$. Since $\Theta$ is a compact subset of $\sigma_{mod} \setminus \bigcup_{\alpha \in \Delta_{\tau_{mod}}^+} \ker \alpha $, the quantity $ \min \{ \alpha(\zeta) \mid \alpha \in \Delta_{\tau_{mod}}^+, \zeta \in \Theta \}$ 
	exists and is positive. 
	
	We now prove $2$. The subset $\Theta= \{\zeta \in \sigma_{mod} \mid \forall \alpha \in \Delta_{\tau_{mod}}^+, \alpha(\zeta) \ge \alpha_0 \}$ has symmetrization $W_{\tau_{mod}}\Theta= \{\xi \in a_{mod} \mid \forall \alpha \in \Delta_{\tau_{mod}}^+, \alpha(\xi) \ge \alpha_0 \}$ which is an intersection of finitely many half-spaces together with the unit sphere, so it is compact and convex. Furthermore $\Theta = \sigma_{mod} \cap W_{\tau_{mod}} \Theta$ is a compact subset of $\interior_{\tau_{mod}}(\sigma_{mod}) \cap \sigma_{mod}$.  
\end{proof}

\subsection{Generalized Iwasawa decomposition}\label{iwasawa decomposition}

Let $p$ be a point in $\X$, $\tau \in \Flagt$ and let $X \in \mfp$ be $\tau$-spanning. Choose a Cartan subspace $X \in \mfa \subset \mfp$, with restricted roots $\Lambda$ and a choice of simple roots $\Delta$ associated to $\sigma \supset \tau$. Recalling the notation in \ref{lambda tau} following Definition \ref{regularity} we define 
\begin{enumerate}
	\item $\mfa_\tau = Z(X) \cap \mfp = \{ Y \in \mfp \mid [X,Y]=0 \}$ and $A_\tau = \exp(\mfa_\tau)$. Note that $\mfa_\tau$ and $A_\tau$ depend on $p$.
	\item The \textit{(nilpotent) horocyclic subalgebra} $\mfn_\tau = \bigoplus_{\alpha \in \Lambda_\tau^+} \mfg_\alpha$ and the \textit{(unipotent) horocylic subgroup} $N_\tau = \exp( \mfn_\tau)$.
	\item The \textit{generalized Iwasawa decomposition} of $\mfg$ is $ \mfg = \mfk \oplus \mfa_\tau \oplus \mfn_\tau .$
	\item The \textit{generalized Iwasawa decomposition} of $G$ is $G = K A_\tau N_\tau = N_\tau A_\tau K$. The indicated decomposition is unique. 
\end{enumerate}
Note that our notation differs from \cite{KLP17}, where $N_\tau$ denotes the full horocyclic subgroup at $\tau$ and $A_\tau$ is the group of translations of the flat factor of the parallel set defined by $p$ and $\tau$, see Section \ref{section parallel sets and horocycles}. In our notation, $N_\tau$ is the unipotent radical of the parabolic subgroup $G_\tau$, see \cite[2.17]{E96}.

\subsection{Antipodal simplices, parallel sets and horocycles}\label{section parallel sets and horocycles}

We say a pair of points $\xi,\eta$ in $\partial \X$ are \textit{antipodal} if there exists a geodesic $c$ with $c(-\infty)=\xi$ and $c(+\infty)=\eta$. Equivalently, $\xi,\eta$ are antipodal if there exists a geodesic symmetry $S_p$ taking $\xi$ to $\eta$.

A pair of simplices $\tau_\pm$ are \textit{antipodal} if there exists some $p \in \X$ such that $S_p \tau_- =\tau _+$, or equivalently if there exists a geodesic $c$ with $c(-\infty) \in \interior (\tau_ -)$ and $c(+\infty) \in \interior (\tau_+)$. If a model simplex $\tau_{mod}$ is $\iota$-invariant then every simplex $\tau$ of type $\tau_{mod}$ has the same type as any of its antipodes.

For antipodal simplices $\tau_\pm$, the \textit{parallel set} $P(\tau_-,\tau_+)$ is the union of (images of) geodesics $c$ with $c(-\infty) \in \tau_-$ and $c(+\infty) \in \tau_+$. Given one such geodesic $c$, we may alternatively define $P(\tau_-,\tau_+)=P(c)$ to be the union of geodesics parallel to $c$, or equivalently to be the union of maximal flats containing $c$. Antipodal $\tau_{mod}$-regular points $\xi,\eta$ lie in the boundary of a unique parallel set $P=P(\tau(\xi),\tau(\eta))$, where $\tau(\xi)$ (resp.\ $\tau(\eta)$) is the unique simplex of type $\tau_{mod}$ in some/every Weyl chamber containing $\xi$ (resp.\ $\eta$). We say that $P(\tau_-,\tau_+)$ \textit{joins} $\tau_-$ and $\tau_+$. The parallel set joining a pair of antipodal Weyl chambers is a maximal flat.

The \textit{horocycle} centered at $\tau \in \Flag(\tau_{mod})$ through $p \in \X$ is denoted $H(p,\tau)$ and is defined to be the orbit $N_\tau \cdot p$. For any $p \in \X$ and $\hat{\tau}$ antipodal to $\tau$, the horocycle $H(p,\tau)$ intersects the parallel set $P(\hat{\tau},\tau)$ in exactly one point. A horocycle is the union of basepoints of strongly asymptotic Weyl sectors/ geodesic rays \cite{KLP14,KLP17}.

\subsection{The $\zeta$-angle and Tits angle}\label{section angles}

We follow \cite{KLP14} in defining the $\zeta$\textit{-angle} between two simplices at a point $p \in \X$. For fixed $p \in X$ and $\zeta$, the $\zeta$-angle provides a metric on $\Flagt$ by viewing it as embedded in the tangent space at $p$ and restricting the angle metric $\angle_p$ to the vectors of type $\zeta$. The $\zeta$-angle also makes sense for $\tau_{mod}$-regular directions by projecting to $\Flagt$. To make this definition, we first fix the auxilary data of a $(\zeta_0,\tau_{mod})$-spanning $\iota$-invariant model ideal point $\zeta=\zeta_{mod} \in \interior (\tau_{mod})$. 

\begin{defn}[{$\zeta$-angle, cf. \cite[Definitions 2.3 and 2.4]{KLP14}}]\label{zeta angle} \hfill
	\begin{enumerate}
		\item For a simplex $\tau \in \Flag(\tau_{mod})$ let $\zeta(\tau)$ denote the unique point in $\interior (\tau)$ of type $\zeta$. 
		\item For a $\tau_{mod}$-regular ideal point $\xi \in \partial \X$, let $\zeta(\xi)=\zeta(\tau(\xi))$ where $\tau(\xi)$ is the simplex spanned by $\xi$. 
		\item Let $p \in \X$, let $\tau,\tau'$ be Weyl chambers in $\partial \X$ and let $x,y \in \overline{\X}$ with $px$ and $py$ $\tau_{mod}$-regular. The $\zeta$\textit{-angle} is given by
		\begin{align*}
		 \angle_p^\zeta(\tau,\tau') & \coloneqq \angle_p(\zeta(\tau),\zeta(\tau')), \\
		 \angle_p^\zeta(\tau,y)  & \coloneqq \angle_p(\zeta(\tau),\zeta(py)), \\
		 \angle_p^\zeta(x,y)  & \coloneqq \angle_p(\zeta(px),\zeta(py)). \\
		\end{align*}
		Note there is a typo in the definition of $\zeta$-angle in \cite[Definition 7.5]{KLP14}.
	\end{enumerate} 
\end{defn}

\begin{figure}[h]
		\centering
		\begin{subfigure}{0.25\textwidth}
		\centering
		\begin{tikzpicture}[scale=0.7]
			\node[circle,scale=0.5,fill=black] (t1) at (0,0) {};
			\node[circle,fill=black,scale=0.5] (t2) at (2,3.46) {};
			\node[circle,fill=black,scale=0.5] (t3) at (4,0) {};
			\draw (t1) -- node[pos=0.2](t12){} (t2) -- (t3) node[pos=0.8](t23){} -- (t1);
			\node[circle,fill=black,scale=0.5,label=right:{$\zeta$}] (z) at (2,1.5) {};
		\end{tikzpicture}
		\caption{$\zeta \in \sigma_{mod}$}	
	\end{subfigure}	
	\begin{subfigure}{0.65\textwidth}
	\centering
	\begin{tikzpicture}
	\coordinate (o) at (0,0);
	\coordinate (r) at (5.5,1);
	\coordinate (l) at (3,1);
	\node[label=right:{$\sigma(X)$}] (u) at (4,3) {};
	\node[label={[xshift=0.2cm, yshift=-0.2cm]:$\zeta(X)$}] (z) at (4,1.8) {};
	\node[label={[xshift=-0.4cm, yshift=-0.2cm]:$X$}] (x) at (5.3,1.2) {};
	\draw (o) -- (r);
	\draw (o) -- (l);
	\draw (o) -- (u.center);
	\draw (u.center) -- (l) -- (r) -- cycle;
	\draw[-{latex[scale=2.5,length=2,width=3]}] (o) -- (4.2,1.8);
	\draw[-{latex[scale=2.5,length=2,width=3]}] (o) -- (x);
	\coordinate (mr) at (-5.5,1);
	\coordinate (ml) at (-3,1);
	\node[label=left:{$\sigma(Y)$}] (mu) at (-4,3) {};
	\node[label={[xshift=-0.55cm, yshift=-0.6cm]:$\zeta(Y)$}] (mz) at (-4,1.8) {};
	\node[label={[xshift=0cm, yshift=-0.3cm]:$Y$}] (y) at (-4.1,2.3) {};
	\draw (o) -- (mr);
	\draw (o) -- (ml);
	\draw (o) -- (mu.center);
	\draw (mu.center) -- (ml) -- (mr) -- cycle;
	\draw[-{latex[scale=2.5,length=2,width=3]}] (o) -- (-4.2,1.8);
	\draw[-{latex[scale=2.5,length=2,width=3]}] (o) -- (y);
	\end{tikzpicture}
	\caption{The $\zeta$-angle between $X$ and $Y$}
	\end{subfigure}
\end{figure}

For $\xi,\eta \in \partial \X$, the \textit{Tits angle} is $$ \angle_{Tits}(\xi,\eta) \coloneqq \sup_{p\in \X} \angle_p(\xi,\eta) .$$
Ideal points $\xi,\eta$ are antipodal if and only if their Tits angle is $\pi$. For $p \in \X$, $\xi,\eta\in \partial \X$, the equality $\angle_p(\xi,\eta) = \angle_{Tits}(\xi,\eta)$ holds if and only if there is a maximal flat $F$ containing $p$ with $\xi,\eta \in \partial F$ and moreover for any $\xi,\eta \in \partial \X$, there exists some maximal flat $F$ with $\xi,\eta \in \partial F$ \cite{E96}.

For simplices $\tau,\tau'$ in $\Flagt$, we may define 
$$ \angle_{Tits}^\zeta(\tau,\tau') \coloneqq \angle_{Tits}(\zeta(\tau),\zeta(\tau') ).$$ 
There are only finitely many possible Tits angles between ideal points of fixed type. Therefore, there exists a bound $\varepsilon(\zeta_{mod})$ such that if $\angle_{Tits}^\zeta(\tau,\tau') > \pi - \varepsilon(\zeta_{mod})$ then $\tau$ and $\tau'$ are antipodal, as observed in \cite[Remark 2.42]{KLP14}. By Remark \ref{roots vs distance}, we have 
$$ \sin(\frac{1}{2} \varepsilon(\zeta_{mod}) ) = \min_{\alpha \in \Lambda_{\tau_{mod}}^+} \ \frac{\alpha(\zeta_{mod})}{\abs{\alpha}} \ge \frac{\zeta_0}{\kappa_0} .$$ 
By the definition of Tits angle, the same holds if the $\zeta$-angle at any point is strictly within $\varepsilon(\zeta_{mod})$ of $\pi$: the inequality
$$ \angle_{Tits}^\zeta(\tau,\tau') \ge \angle_p^\zeta(\tau,\tau') > \pi - \varepsilon(\zeta_{mod}) $$
implies that $\tau$ and $\tau'$ are antipodal. Since $\zeta_0 \le \kappa_0 < 2 \kappa_0$ we have
$$ \sin \frac{1}{2} \frac{\zeta_0^2}{\kappa_0^2} \le \frac{1}{2} \frac{\zeta_0^2}{\kappa_0^2} < \frac{\zeta_0}{\kappa_0} \le \sin \frac{1}{2} \varepsilon(\zeta_{mod}) ,$$
and we obtain the estimate $\frac{\zeta_0^2}{\kappa_0^2} < \varepsilon(\zeta_{mod})$. We record this observation in the following lemma.
\begin{lem}[{Cf.\ \cite[Remark 2.42]{KLP14}}]\label{angle implies antipodal}
	If the inequality $\angle_p^\zeta(\tau_-,\tau_+) \ge \pi - \frac{\zeta_0^2}{\kappa_0^2}$ holds for some $p \in \X$ then $\tau_-$ is antipodal to $\tau_+$. In other words, $\frac{\zeta_0^2}{\kappa_0^2} <\varepsilon(\zeta_{mod})$.
\end{lem}
%Comment: the inequality zeta_0^2/kappa_0^3 < \varepsilon(\zeta_{mod}) does not hold in general.

\section{Estimates}\label{section estimates}

%The exposition of this section has been updated.
This section contains the main contributions of this paper. We prove several explicit estimates in the symmetric space that we will use in Section \ref{local to global section} to give a quantified version of the local-to-global principle for Morse quasigeodesics. Qualitative versions of these estimates appear in \cite{KLP14, KLP17}, but there the proofs rely on topological arguments that do not produce explicit bounds. For example, in subsection \ref{section zeta projection}, Lemma \ref{zeta projection} we consider the natural projection from $(\alpha_0,\tau_{mod})$-regular vectors in $\mfp$ to $\Flagt$. This map is the restriction of a smooth map to a compact submanifold with boundary, so an abstract proof of the existence of a Lipschitz constant is not hard. However, that approach is not suitable for our purposes, so we apply Corollary \ref{morse function estimate} to obtain an explicit local Lipschitz constant. Note that such an estimate cannot be uniform for all $\alpha_0 >0$ and therefore must depend on $\alpha_0$.  

A crucial notion, introduced in \cite{KLP14}, is the $\zeta$-angle, denoted $\angle^\zeta$, see Section \ref{section angles}. Recall that $\zeta=\zeta_{mod}$ is a fixed type in the interior of $\tau_{mod}$. Moreover we assume that $\zeta$ is $(\zeta_0,\tau_{mod})$-regular and that $\zeta$ and $\tau_{mod}$ are $\iota$-invariant, see Definition \ref{spanning} and Section \ref{visual boundary}. For fixed $p \in X$ and $\zeta$, the $\zeta$-angle provides a metric on $\Flagt$ by viewing it as embedded in the tangent space at $p$ and restricting the angle metric $\angle_p$ to the vectors of type $\zeta$. The $\zeta$-angle also makes sense for $\tau_{mod}$-regular directions by projecting to $\Flagt$. 

The organization of the section is as follows. In subsection \ref{section Bp} we relate the Riemannian metric on $\X$ to algebraic data on $\mfg$, e.g.\ the Killing form $B$ and the canonical inner product $B_p$. In subsection \ref{section perturbation regular} we use the vector-valued triangle inequality to control the regularity of bounded perturbations of long regular geodesic segments. In subsection \ref{section zeta projection}, we prove Lemma \ref{zeta projection}, which allows us to bound $\angle_p^\zeta(x,y)$ in terms of $\alpha_0,\zeta_0$ and $\angle_p(x,y)$. In subsection \ref{projecting curves in G} we prepare a technique for the subsequent subsections, where we bound the lengths of certain non-geodesic curves in $\X$ which are images of curves in $G$ under the orbit map. In subsection \ref{section cone rotation}, the curve lies in the subgroup stabilizing a point, and we bound the distance the midpoint of a segment can move when we move one endpoint a bounded amount, assuming the segment is long enough. Subsection \ref{section strongly asymptotic rays} is roughly similar; there we bound the distance between points far along on strongly asymptotic geodesic rays (so the curve in $G$ lies in a unipotent horocyclic subgroup). These combine to yield a crucial estimate in Corollary \ref{midpoint projection}, which implies that if a pair of points are in the $D$-neighborhood of a diamond, then their midpoint is close to the diamond; moreover the distance from the midpoint to the diamond becomes arbitrarily small as the points move farther apart. In the remaining subsections, we show that distance to a corresponding parallel set controls the corresponding $\zeta$-angles (Corollary \ref{distance to angles near pi}) and vice-versa (Lemma \ref{angle to distance}). Along the way we provide some control for the Lie derivatives of gradients of Busemann functions with respect to Killing vector fields, see the proofs of Lemma \ref{simplex displacement} and Lemma \ref{angle to distance}.

\subsection{Useful properties of the inner product $B_p$ on $\mfg$}\label{section Bp}

We remind the reader that our convention is that the Riemannian metric on $\X$ is the one induced by the Killing form, see Equation \ref{metric is Killing form}. Recall that each point $p \in \X$ induces an inner product $B_p$ on $\mfg$ and the evaluation map $\evp \colon \mfg \to T_p \X$, see Section \ref{cartan decomposition}. We first relate the inner product $B_p$, the Killing form $B$ on $\mfg$, and the Riemannian metric $\langle\cdot ,\cdot \rangle$ at $p$. 

\begin{lem}\label{metric and Bp}
	For any $X,Y \in \mfg$ and $p \in \X$, 
	$$ 2 \langle \evp X, \evp Y \rangle = B(X,Y) + B_p(X,Y) .$$
	In particular, any $U$ in $\mfn_\tau$ or $\mfg_\alpha$ is $\ad$-nilpotent, so $B(U,U)=0$ and $\pnorm{U} = \sqrt{2} \abs{ \evp U}$, see Section \ref{iwasawa decomposition}.
\end{lem}
Recall that $\vartheta_p$ is a Lie algebra automorphism so $\vartheta_p [X,Y] = [\vartheta_p X,\vartheta_pY]$ and $B(\vartheta_pX,\vartheta_pY)=B(X,Y)$.
\begin{proof}
	The kernel of $\evp$ is the $+1$-eigenspace for $\vartheta_p$, so for any $X \in \mfg$, $2\evp X = \evp (X - \vartheta_p X )$ and
	\begin{align*}
	& 4 \langle \evp X,\evp Y \rangle_p = \langle \evp (X-\vartheta_pX), \evp (Y-\vartheta_pY) \rangle_p = B(X-\vartheta_p X, Y-\vartheta_p Y) \\
	& = B(X,Y) +B(\vartheta_p X,\vartheta_pY) - B(\vartheta_pX,Y) -B(X,\vartheta_pY) = 2B(X,Y) + 2B_p(X,Y). \qedhere
	\end{align*} 
\end{proof}

Next we show that the transpose on $\End \mfg$ with respect to $B_p$ restricts to $-\vartheta_p$ on the image of the adjoint representation.
\begin{lem}\label{transpose}
	For $X,Y,Z \in \mfg$, $B_p(\ad X (Y),Z) = B_p(Y, \ad (-\vartheta_p X)(Z))$.
\end{lem}
\begin{proof}
	We have 
	\begin{align*}
	& B_p(\ad X (Y),Z) =-B(\vartheta_p \ad X (Y),Z) = -B(\ad(\vartheta_pX)(\vartheta_pY),Z) \\
	& = -B(\vartheta_pY,\ad(-\vartheta_pX) (Z)) = B_p(Y, \ad(-\vartheta_p X) (Z)) 
	\end{align*}
	where we have used that $\ad \vartheta_p X$ is skew-symmetric relative to $B$. 
\end{proof}

Third, we bound $B(\ad X (Y),Z)$ by the product of the $B_p$-norms of $X,Y$ and $Z$ and bound the operator norm of $\ad X$ by $\pnorm{X}$ along the way.

\begin{lem}\label{ad bound}
	Let $X,Y,Z \in \mfg$ and let $p \in \X$ induce the inner product $B_p$ on $\mfg$. Consider the operator norm $\abs{\cdot}_{op}$ and Frobenius norm $\abs{\cdot}_{Fr}$ on $\End \mfg$ induced by $B_p$. Then	
	\begin{enumerate}
		\item $\abs{ \ad Y}_{op} \le \abs{\ad Y}_{Fr}=\pnorm{Y}$,
		\item $B(X,\ad Y (Z)) \le \pnorm{X}\pnorm{Y}\pnorm{Z} $, and
		\item For $Y \in \mfp$, $\pnorm{[Y,X]} \le \kappa_0 \pnorm{Y} \pnorm{X}$.
	\end{enumerate} 
\end{lem}
\begin{proof}
	Recall that the operator norm of a linear transformation is the largest singular value, while the Frobenius norm is the square root of the sum of the singular values squared. Therefore
	$$	\abs{\ad X}^2_{op} \le \abs{\ad X}^2_{Fr} = \trace_{\mfg} (\ad (-\vartheta_p X) \circ \ad X) = B_p(X,X) $$
	by Lemma \ref{transpose}, proving the first claim. Using this, we have
	$$ B(X,\ad Y (Z)) = -B_p(\vartheta_p X, \ad Y (Z)) \le \pnorm{\vartheta_pX} \pnorm{\ad Y (Z)} \le \pnorm{X} \abs{\ad Y}_{op} \pnorm{Z} \le \pnorm{X} \pnorm{Y} \pnorm{Z} .$$ 
	If $Y \in \mfp$, we may choose a maximal abelian subspace $\mfa$ of $\mfp$ containing $Y$ and decompose $X=\sum_{\alpha \in \Lambda \cup \{0\}} X_\alpha$ according to the associated restricted root space decomposition, which is $B_p$-orthogonal. Therefore
	$$ \pnorm{[Y,X]}^2 = \pnorm{ \sum_{\alpha \in \Lambda} \alpha(Y) X_\alpha }^2 = \sum_{\alpha \in \Lambda} \alpha(Y)^2 \pnorm{X_\alpha }^2 \le \kappa_0^2 \pnorm{Y}^2 \pnorm{X}^2 $$
	where $\kappa_0$ is the maximum of $\{\alpha(A)\mid \alpha \in \Lambda, A \in \mfa, \abs{A}=1 \}$ see Section \ref{curvature}.
\end{proof}

Fourth, we need to compare the norms induced by $p,q \in \X$ in terms of $d(p,q)$. 
\begin{lem}\label{Bp comparison}
	Let $p,q \in \X$, $g \in G$ and $X \in \mfg$. Then
	\begin{enumerate}
		\item $ \vartheta_{gp} \circ \Ad(g) = \Ad(g) \circ \vartheta_p ,$
		\item $ \pnorm{X} = \abs{\Ad(g) X}_{B_{gp}} ,$
		\item $ \pnorm{X} \le e^{\kappa_0 d(p,q)} \qnorm{X} .$
	\end{enumerate}
\end{lem}
\begin{proof}
	The point stabilizer $G_{gp}$ is $g G_p g^{-1}$ and it follows that $\Ad(g)$ takes $\vartheta_p$ to $\vartheta_{gp}$. This, together with the $\Ad$ invariance of the Killing form implies $2$. For the last point, choose a maximal flat $F$ containing $p$ and $q$ and let $\mfg = \mfg_0 \oplus \bigoplus \mfg_\alpha$ be the restricted root space decomposition corresponding to $p$ and $F$. Then we may choose $A \in \mfa$, the maximal abelian subspace $\mfa$ of $\mfp$ corresponding to $F$, such that $e^Ap=q$, and then 
	$$ \pnorm{X} = \qnorm{ e^{\ad A} X} = \qnorm{\sum_{\alpha \in \Lambda \cup \{0\}} e^{\alpha(A)}X_\alpha} \le e^{\kappa_0 d(p,q)} \qnorm{X} ,$$
	using the restricted root space decomposition of $X$ and the fact that the restricted root space decomposition is $B_q$-orthogonal. 
\end{proof}

\subsection{Perturbations of long, regular segments}\label{section perturbation regular}

We will need to control the regularity of bounded perturbations of long regular geodesic segments. The following Lemma is an explicit version of Lemma 3.6 in \cite{KLP18b}. This assertion also appears in the proof of Lemma 7.10 in \cite{KLP14}.

\begin{lem}\label{regular projections}
	Suppose $xy$ is an $(\alpha_0,\tau_{mod})$-regular geodesic segment with $d(x,y)\ge l$ and let $x',y'$ be points in  $\X$ satisfying $d(x,x')\le \delta_x$ and $d(y,y') \le \delta_y$. If 
	$$ \alpha_0 - \frac{(\delta_x+\delta_y)(\alpha_0+\kappa_0)}{l-\delta_x-\delta_y} \ge \alpha_0' $$
	then $x'y'$ is $(\alpha_0',\tau_{mod})$-regular.
\end{lem}

We will often apply this lemma in the case $\delta_x=\delta_y=D$.

\begin{proof}
	We apply Corollary \ref{vector valued triangle inequality}, the triangle inequality for $\vec{d}$-distances: 
	$$ \abs{\vec{d}(x,y)-\vec{d}(x',y')} \le d(x,x') +d(y,y') \le \delta_x+\delta_y .$$
	Similarly, $\abs{d(x,y)-d(x',y')} \le d(x,x') +d(y,y') \le \delta_x+\delta_y ,$ so $d(x',y') \ge l-\delta_x+\delta_y$ and 
	$$ \frac{d(x,y)}{d(x',y')} \ge 1- \frac{\delta_x+\delta_y}{d(x',y')} \ge   1- \frac{\delta_x+\delta_y}{l-\delta_x-\delta_y}  . $$
	For any $\alpha \in \Delta_{\tau_{mod}}^+$,
	$$ \frac{\alpha(\vec{d}(x',y'))}{d(x',y')} \ge \frac{\alpha_0 d(x,y) - \delta_x \kappa_0-\delta_y \kappa_0}{d(x',y')} \ge \alpha_0 \left( 1- \frac{\delta_x+\delta_y}{l-\delta_x-\delta_y} \right) - \frac{ (\delta_x+\delta_y) \kappa_0}{l-\delta_x-\delta_y} = \alpha_0 - \frac{(\delta_x+\delta_y)(\alpha_0+\kappa_0)}{l-\delta_x-\delta_y} \ge \alpha_0' .$$
\end{proof}

\subsection{Angle comparison to Euclidean space}

When $p,q,r$ are points in $\X$ such that $d(p,q)$ is much larger then $d(q,r)$, we provide an upper bound for the Riemannian angle $\angle_p(q,r)$ by comparing to Euclidean space. The following estimate is surely not new, but we could not find a direct reference so we give a proof. 
%Comment: it looks like this estimate is the best of its kind. It's not easy to find a reference for some reason. What's nice here is that the third distance d(p,r) doesn't need to be mentioned. 

\begin{lem}\label{euclidean angle}
	Let $p,q,r$ be non-collinear points in $\X$. Then 
	$$ \sin \angle_p(q,r) \le \frac{d(q,r)}{d(p,q)}. $$
\end{lem}
The convenience of this estimate is that the third possible distance $d(p,r)$ does not appear.
\begin{proof}
	Let $X,Y \in \mfp$ such that $e^Xp=q$ and $e^Yp=r$. Then $\abs{X}=d(p,q)$ and $d(X,Y) \le d(q,r)$ and we may assume that $d(p,q) > d(q,r)$. In Euclidean space, the comparison holds: among vectors $Y'$ with $d(X,Y') \le d(X,Y)$, the largest angle occurs for a vector $Y'$ forming a right triangle with $X$ as hypotenuse.  Then 
	$$\sin \angle(X,Y) \le \sin \angle(X,Y') = \frac{d(X,Y')}{\abs{X}} \le \frac{d(q,r)}{d(p,q)} .$$
\end{proof}

% I had the footnote: \footnote{The function $v \mapsto \left\langle \frac{X}{\abs{X}},\frac{v}{\abs{v}} \right\rangle$ has gradient $(\abs{v}^2 X- \left\langle v,X \right\rangle v)/\abs{v}^3$ and is minimized on the sphere centered at $X$ with radius $d(X,Y)$ at a point $w$ with $\langle w,X \rangle =\abs{w}^2$.}

\subsection{Projecting regular vectors to flag manifolds}\label{section zeta projection}

Recall that we have a fixed type $\zeta=\zeta_{mod}$ which is $(\zeta_0,\tau_{mod})$-spanning. For a $\tau_{mod}$-regular $X \in \mfp$, define $\zeta(X)$ to be the unique vector in a common closed Weyl chamber as $X$ of type $\zeta$. Note that $\zeta(X)$ is the unique maximizer for $B(X,\cdot) \colon \Ad(K)Z \to \R$ where $Z \in \mfp$ is any vector of type $\zeta$ by Corollary \ref{morse function estimate}. This map $\zeta$ from $\tau_{mod}$-regular elements of $\mfp$ to $\Ad(K)Z$ is a smooth fiber bundle. In the next lemma we show that nearby $\tau_{mod}$-regular points project to nearby points on $\Ad(K)Z$ in the metric induced by viewing $\Ad(K)Z$ as a Riemannian submanifold of $\mfp$. Note that one expects a local Lipschitz constant proportional to $\frac{1}{\alpha_0}$ by considering vectors near the walls $\ker \alpha$ for $\alpha \in \Delta_\tau^+$. 

% Why is \zeta smooth? I can't find a proof in KLP. They prove that the analogous map on the visual boundary is continuous and state without proof that it is a fiber bundle. Should I prove this in the appendix?
% The smoothness should follow from a standard result about "G-slices"

\begin{lem}\label{zeta projection}
	Let $X,X'$ be $(\alpha_0,\tau)$-regular unit vectors in $\mathfrak{p}$ with $d_\mathfrak{p}(X,X') \le \alpha_0$. Write $Z=\zeta(X)$ and $Z'=\zeta(X')$. Then the Riemannian distance on $\Ad(K)Z$ from $Z$ to $Z'$ is bounded by the distance in $\mfp$ from $X$ to $X'$:
	$$ d_{\Ad(K)Z}(Z,Z') \le \frac{1}{\alpha_0 \zeta_0}d_{\mfp}(X,X') .$$
\end{lem}

\begin{proof}
	Let $t \mapsto X_t$ be a unit-speed line segment from $X$ to $X'$ in $\mfp$. Let $\{X^i \}_{i=1}^{\dim \mfp}$ be linear coordinates on $\mfp$, and we may assume that the derivative of $t \mapsto X_t$ is $\pdv{X^1}$. Since $d_{\mfp}(X,X') \le \alpha_0$ each $X_t$ is $(\frac{\alpha_0}{2},\tau_{mod})$-regular. Write $Z_t=\zeta(X_t)$ and note that $t \mapsto Z_t$ is a smooth curve on $\Ad(K)Z$. To prove the claim we will show that $\abs{\dv{Z_t}{t}} \le \frac{1}{\alpha_0\zeta_0}$, where we restrict the inner product on $\mfp$ to a Riemannian metric on $\Ad(K)Z$. 
	
	Restricting the domain of $B$, we write $B \colon \mfp \times \Ad(K)Z \to \mathbb{R}$. Near $(X_0,Z_0)=(X_{t_0},Z_{t_0})$, we have coordinates $\{Z^j\}_{j=1}^{\dim \Ad(K)Z}$ on $\Ad(K)Z$. We may assume that $Z_t$ is an immersion at $Z_0$ because the set $\{t \mid \abs{\dv{Z_t}{t}}=0\}$ does not contribute to the arclength of $Z_t$ and furthermore up to a change of coordinates we may assume that $\dv{Z_t}{t}=\pdv{Z^1}$. On this coordinate patch $U$, we obtain the function $B_j \colon \mfp \times U \to \mathbb{R}$ defined by $B_j(X'',Z'')\coloneqq \d{B}_{(X'',Z'')}(\pdv{Z^j})$. Along the curve $t \mapsto (X_t,Z_t)$, the function $B_j$ is identically $0$ (where defined) since $Z_t$ maximizes $B(X_t,\cdot)$ on $\Ad(K)Z$. Differentiating $B_j(X_t,Z_t)=0$ in $t$, we obtain
	$$ 0 = \d{B_j}_{(X_t,Z_t)} \left( \pdv{X^1},\pdv{Z^1} \right) = \pdv{B_j}{X^1}+\pdv{B_j}{Z^1}  .$$
	Observe that
	$$ \pdv{B_j}{Z^1}_{(X_t,Z_t)} = \left( \pdv{Z^1}\pdv{Z^j} B \right)_{(X_t,Z_t)} = \Hess(B) \left(\pdv{Z^1},\pdv{Z^j} \right)_{(X_t,Z_t)} = \Hess(B(X_t,\cdot))\left(\pdv{Z^1},\pdv{Z^j}\right)_{Z_t} $$
	so by Corollary \ref{morse function estimate} we have 
	$$ \abs{ \pdv{B_j}{Z^1} } \ge \alpha_0 \zeta_0 \abs{ \left\langle \pdv{Z^1},\pdv{Z^j} \right\rangle } .$$
	In particular, along $(X_t,Z_t)$ and setting $j=1$, we have 
	$$ \alpha_0\zeta_0 \abs{ \pdv{Z^1} }^2 \le \abs{ \pdv{B_1}{X^1}_{(X_t,Z_t)}} = \abs{B_1 \left( \pdv{X^1},Z_t \right)} = \abs{B \left( \pdv{X^1},\pdv{Z^1} \right) } \le \abs{\pdv{Z^1}} $$
	since $\pdv{X^1}$ is a unit vector. We obtain for all $t$
	$$ \abs{\pdv{Z^1}}\le \frac{1}{\alpha_0\zeta_0} $$
	and the the claim is proven. 
\end{proof}

\subsection{Projecting curves in $G$ to $\X$}\label{projecting curves in G}

In this subsection we prepare to estimate the length of curves in $\X$ which are images of curves in $G$ under the orbit map. We begin by comparing the speeds of two such curves related by right-translation. We apply this result in the next section to Lemma \ref{cone rotation} for a curve in $K$, and in the following section to Lemma \ref{strongly asymptotic geodesics} for a curve in the subgroup $N_\tau$. 

For an element $g \in G$, we let $l_G \colon G \to G, l_g(h)=gh$ denote left translation and $r_g \colon G \to G, r_g(h)=hg$ denote right translation. We denote by $\conj_g \colon G \to G$ the conjugation map $\conj_g(h) = ghg^{-1}$.

\begin{lem}\label{right translated curves}
	Let $g \colon \mathbb{R} \to G$ be a curve in $G$, let $h \in G$ and let $p \in \X$. Write $q_h(s)=g(s)hp$. If  $\dot{g}(s) = (\dd{l_{g(s)}})_1 X_s$ then
	$$ \abs{\dot{q_h}(s)} = \abs{\evp \Ad(h^{-1}) X_s} .$$
\end{lem}

\begin{proof}
	The curve $q_h(s)=g(s)hp$ has the same speed as $c_h(s)=h^{-1}g(s)hp$ since $h^{-1}$ is an isometry. Writing 
	$$ c_h(s) = p \circ \conj_{h^{-1}} \circ g (s) $$
	and differentiating with respect to $s$ we have 
	$$ \dot{c_h}(s) = (\dd{\orb_p})_{h^{-1}gh} \circ (\dd{ \conj_{h^{-1}}})_{g(s)} \circ \dot{g}(s) .$$
	For any $a,b \in G$ and $X \in T_1G$ we have 
	$$ (\dd{\conj_a})_b (\dd{l_b})_1 X = (\dd{l_a})_{ba^{-1}} (\dd{r_a^{-1}})_b (\dd{l_b})_1 X = (\dd{l_a})_{b a^{-1}}(\dd{l_b})_{a^{-1}} (\dd{l_a^{-1}})_1 (\dd{l_a})_{a^{-1}} (\dd{r_a^{-1}})_1 X   = \dd{l_{aba^{-1}}} \Ad(a) X $$
	We also have $(\dd{\orb_p})_a (\dd{l_a})_1 =\dd{a}_p (\dd{\orb_p})_1$, so if $\dot{g}(s) = \dd{l_{g(s)}} X_s$, then 
	$$ \dot{c_t}(s) = (\dd{\orb_p})_{h^{-1}gh} \circ (\dd{l_{h^{-1}gh}})_1 \Ad(h^{-1}) X_s  = (\dd{h^{-1}gh})_p (\dd{\orb_p})_1 \Ad(h^{-1}) X_s .$$
	This implies 
	$$ \abs{\dot{q_h}(s)} = \abs{ \dot{c_h}(s) } = \abs{(\dd{\orb_p})_1 \Ad(h^{-1}) X_s} = \abs{\evp \Ad(h^{-1}) X_s} $$
	and completes the proof.
\end{proof}

\subsection{Weyl cones forming small angles}\label{section cone rotation}

In this subsection, we show that if $q \in V(p,\st(\tau),\alpha_0)$ and $ r \in V(p,\st(\tau'),\alpha_0)$ with $d(p,q)$ much larger than  $d(q,r)$, the midpoint of $pq$ is close to $V(p,\st(\tau'),\alpha_0)$. 

\begin{figure}[h]
	\centering
	\begin{tikzpicture}
	
	\node[circle,scale=0.5,fill=black,label=left:{$p$}] (p) at (2,3) {};
	\coordinate[label=right:{$\tau_D$}] (tauD) at (14,6);
	\coordinate[label=right:{$\tau_t$}] (taut) at (14.5,3.5);
	\coordinate[label=right:{$\tau$}] (tau) at (14,0.5) {};
	
	\node[circle,fill,scale=0.5,label=right:{$r$}] (r) at (10,5) {};
	
	\draw[thin,dashed,-{latex}] (p) -- (tau)
	node[circle,fill,scale=0.5,pos=0.4,label=below:{$m$}] (m) {}  node[circle,fill,scale=0.5,pos=0.8,label=below:{$q$}] (q) {};
	
	\draw[draw=none] (p) -- (tauD)
	node[circle,fill,scale=0.5,pos=0.4,label=above:{$m(D)$}] (mD) {}  node[circle,fill,scale=0.5,pos=0.8,label=left:{$q(D)$}] (qD) {};
	
	\draw[thin,dashed,-{latex}] (q) .. controls (12.3,3.3) .. (qD) node[circle,fill,scale=0.5,pos=0.52,label=right:{$q(t)$}] (qt) {};
	\draw[thin,dashed,-{latex}] (m) .. controls (7.3,3.3) .. (mD)
	node[circle,fill,scale=0.5,pos=0.47,label=right:{$m(t)$}] (mt) {};
	\draw[thin,dashed,-{latex}] (q) -- (r)
	node[circle,fill,scale=0.5,pos=0.57,label=left:{$c(t)$}] (ct) {};
	
	\draw (14.5,1) -- (p) -- (13.5,0); % draw a cone at tau
	\draw (14.5,6.5) -- (p) -- (13.5,5.5); % draw a cone at tauD
	\draw (15,4) -- (p) -- (14,3); % draw a cone at tauD
	\end{tikzpicture}
	\caption{Weyl cones forming a small angle}
\end{figure}

\begin{lem}\label{cone rotation}
	Let $p,q,r \in \X$ with $pq$ an $(\alpha_0,\tau)$-regular geodesic ray with $d(p,q)\ge 2l$ and $d(q,r) \le D$. If $m=\midp(p,q)$, $K=\Stab_G(p)$ and
	$$ \alpha_0 - \frac{D(\kappa_0+\alpha_0)}{2l-D} \ge \alpha_0' >0 $$
	and 
	$$ \frac{1}{2}(e^{2\kappa_0D}-1)[\sinh(\alpha_0'(2l-D))]^{-2} \le 3e^{2\kappa_0D} $$
	then there exists $k \in K$ such that $km \in V(p,\st(\tau(pr)),\alpha_0)$  and $d(m,km)$ is at most $2De^{\kappa_0D-\alpha_0l}$.
\end{lem}

The first inequality guarantees that $pr$ is $\tau_{mod}$-regular so that $\tau(pr)$ is well-defined. The second requirement looks strange and involves an arbitrary choice, but is extremely mild and serves our purposes well. (When we apply this Lemma, we will have a bounded $D$ and a large $l$.) Compared to other variations of Lemma \ref{cone rotation} we could present here, the given version has a less cumbersome upper bound in the conclusion of the Lemma.
%Comment: I did actually check that in my example the stronger inequality doesn't improve the final estimate.

\begin{proof}
	We may assume that $d(p,q)=2l$ and $d(q,r)=D$. Let $c \colon [0,D]\to \X$ be the unit-speed geodesic from $q$ to $r$. We have $l$ large enough that Lemma \ref{regular projections} implies that each ray $pc(t)$ is $(\alpha_0',\tau_{mod})$-regular and defines a simplex $\tau_t\coloneqq \tau(pc(t))$. We may decompose 
	$$ \dot{c}(t) = N_{c(t)}+T_{c(t)} $$
	so that $T_{c(t)}$ is tangent to $V_t\coloneqq V(p,\ost(\tau_t))$ and $N_{c(t)}$ is normal to $V_t$. There is a unique $X_t \in \mfk^{\tau_t} \subset T_1K$ such that $\eva_{c(t)}X_t=N_{c(t)}$, and we extend each $X_t$ to a \textit{right}-invariant vector field on $K$. We may view this time-dependent vector field as vector field supported on a compact neighborhood of $[0,D] \times K$, so it defines a flow and in particular a curve $k\colon [0,D] \to K$ with $k(0)=1$ and $\dot{k}(t)=(X_t)_{k(t)} = (\dd{r_{k(t)}})_1 X_t$. 
	
	Viewing $\mathfrak{k}$ as $T_1K$, it is convenient to set $X_t = \Ad(k(t))Y_t$ and work with the time-dependent tangent vector $Y_t \in \mfk^\tau$. We have $\dot{k}(t)=(\dd{l_{k(t)}})_1 Y_t$, so we may extend $Y_t$ to the unique \textit{left}-invariant vector field agreeing with $X_t$ along $k(t)$. 
	
	We may now write $c(t) = k(t) v(t)$ where $v(t) \in V(p,\st(\tau),\alpha_0')$. Since $T_{c(t)}=\dd{k(t)} \dot{v}(t)$ we have $\abs{\dot{v}} \le \abs{\dot{c}}$, so 
	$$ d(k(t)v(0),k(t)v(t)) = d(v(0),v(t)) \le t \le D .$$
	
	Setting $q(t)=k(t)q$, we have $\abs{\dot{q}(t)}=\abs{ \eva_q Y_t}$ by Lemma \ref{right translated curves} and by Lemma \ref{Bp comparison}.3 we have
	\begin{equation}\label{equation cone rotation}
	2 \abs{\eva_q Y_t}^2 - \abs{Y_t}^2_B = \abs{Y_t}^2_{B_q} \le e^{2\kappa_0t} \abs{Y_t}^2_{B_{v(t)}} = e^{2 \kappa_0t} \left( 2 \abs{ \eva_{v(t)} Y_t}^2 - \abs{Y_t}^2_B \right)
	\end{equation}
	where $\abs{Y_t}^2_B=B(Y_t,Y_t)$ is nonpositive. 
	
	For large $l$, the evaluation of $Y_t$ at $v$ bounds the Killing form norm of $Y_t$: We choose a maximal flat containing $p$ and $v=e^Ap$ and, suppressing $t$, write $Y_t=\sum_{\alpha \in \Lambda_\tau^+} Y_\alpha+ Y_{-\alpha}$ with $Y_\alpha \in \mfg_\alpha$ and compute  
	\begin{align*}
	\abs{ \eva_v Y_t }^2 & = \abs{ \evp \Ad(e^{-A}) Y_t }^2 && \text{ by Lemma \ref{Bp comparison}.2 } \\
	& = \sum_{\alpha \in \Lambda_\tau^+} \abs{  \left(e^{\alpha(A)}-e^{-\alpha(A)} \right) \evp Y_\alpha }^2 && \text{ since } Y_\alpha+Y_{-\alpha} \in \ker \evp \\
	& = \frac{1}{2} \sum_{\alpha \in \Lambda_\tau^+} \left(e^{\alpha(A)}-e^{-\alpha(A)} \right)^2 \pnorm{Y_\alpha}^2  \\
	& \ge \frac{1}{2} \sum_{\alpha \in \Lambda_\tau^+} \left[ 2 \sinh( \alpha_0' (2l-t)) \right]^2 \pnorm{Y_\alpha}^2 && \text{ since } \alpha(A) \ge \alpha_0'(2l-t) \text{ by regularity}\\
	& = \frac{1}{2} \left[ 2 \sinh( \alpha_0' (2l-t)) \right]^2 \sum_{\alpha \in \Lambda_\tau^+} \pnorm{Y_\alpha}^2 \\
	& = \frac{1}{4} \left[ 2 \sinh( \alpha_0' (2l-t)) \right]^2 (-\abs{Y}^2_B) .
	\end{align*}	
	This bound $-[\sinh( \alpha_0' (2l-t) )]^2 \abs{Y_t}^2_B \le \abs{ \eva_{v(t)} Y_t}^2 $ together with (\ref{equation cone rotation}) implies
	$$ 2 \abs{\eva_q Y_t}^2 \le e^{2 \kappa_0t} 2 \abs{ \eva_{v(t)} Y_t}^2  -(e^{2 \kappa_0Dt}-1) \abs{Y_t}^2_B  \le 2 \abs{ \eva_{v(t)} Y_t}^2 \left[ e^{2\kappa_0t} + \frac{1}{2}(e^{2\kappa_0t}-1)[ \sinh(\alpha_0'(2l-t))]^{-2} \right] .$$	
	
	We now write $m(t)=k(t)m$ where $m=\midp(p,q) = e^{lW}p$ for $W \in \mfp$. For $t \ge 0$, using $\alpha(W) \ge \alpha_0 >0$ for all $\alpha \in \Lambda_\tau^+$ and Lemma \ref{right translated curves}, we have 
	\begin{align*}
	\abs{ \dot{m}(t) }^2 & = \abs{ \evp \Ad(e^{-lW}) Y_t }^2 \\
	& = \frac{1}{2} \sum_{\alpha \in \Lambda_\tau^+} \left(e^{l\alpha(W)}-e^{-l\alpha(W)} \right)^2 \pnorm{Y_\alpha}^2 \\
	& \le \frac{1}{2} \sum_{\alpha \in \Lambda_\tau^+} \left[ \left( e^{2l\alpha(W)} - e^{-2l\alpha(W)} \right) e^{-l\alpha(W)} \right]^2 \pnorm{Y_\alpha}^2 \\
	& \le \frac{1}{2} \sum_{\alpha \in \Lambda_\tau^+} \left( e^{2l\alpha(W)} - e^{-2l\alpha(W)} \right)^2 e^{-2\alpha_0 l} \pnorm{Y_\alpha}^2 \\
	& = e^{-2\alpha_0 l} \abs{\dot{q}(t)}^2 \\
	& \le e^{-2\alpha_0 l} \left[ e^{2\kappa_0t} + \frac{1}{2}(e^{2\kappa_0t}-1)[ \sinh(\alpha_0'(2l-t))]^{-2} \right]
	\end{align*}
	
	The length of $m$ is then 
	\begin{align*}
	\int_0^D \abs{\dot{m}(t)} \dd{t} & \le \int_0^D e^{-\alpha_0l} \sqrt{ e^{2\kappa_0t} + \frac{1}{2}( e^{2\kappa_0t}-1)[\sinh( \alpha_0'(2l-t))]^{-2} } \dd{t} \\
	& \le \int_0^D e^{-\alpha_0l} \sqrt{ e^{2\kappa_0D} + 3e^{2\kappa_0D} } \dd{t} \le 2D e^{\kappa_0D - \alpha_0l} 
	\end{align*}
	and $k(D)$ is the desired isometry. 
\end{proof}

It is possible to give a slightly stronger upper bound in Lemma \ref{cone rotation}, but the improvement would be inconsequential when we apply this Lemma in Section \ref{local to global section} while making the already cumbersome statements even harder to read.

\subsection{Strongly asymptotic geodesics and Weyl cones}\label{section strongly asymptotic rays}

The next estimate says that a point far along an $(\alpha_0,\tau)$-regular geodesic ray gets arbitrarily close to any given parallel set $P(\hat{\tau},\tau)$. The following Lemma is a quantified version of Lemma 2.39 in \cite{KLP14}. 

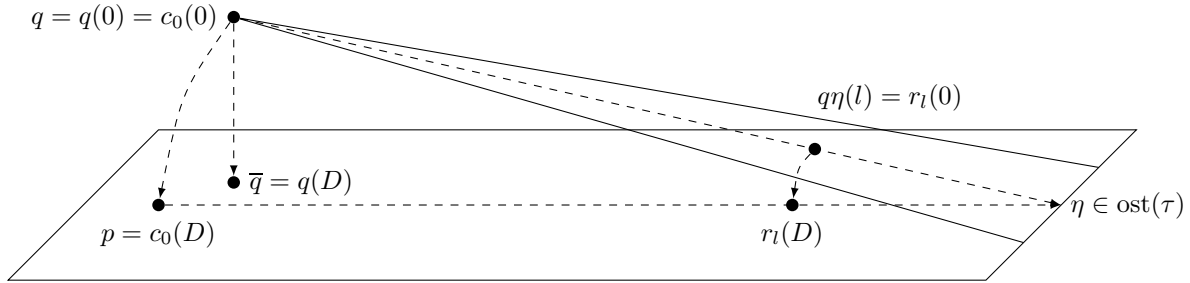
\begin{figure}[h]
	\centering
	\begin{tikzpicture}
	
	\node[circle,scale=0.5,fill=black,label=below:{$p=c_0(D)$}] (p) at (2,0.5) {};
	\coordinate[label=right:{$\eta \in \ost(\tau)$}] (tau) at (14,0.5) {};
	\node[circle,fill,scale=0.5,label=left:{$q=q(0)=c_0(0)$}] (q) at (3,3) {};
	\node[circle,fill,scale=0.5,label=right:{$\overline{q}=q(D)$}] (qb) at (3,0.8) {};
	
	\draw[thin,dashed,-{latex}] (q) -- (qb) {};
	\draw[thin,dashed,-{latex}] (q) -- (tau)
	node[circle,fill,scale=0.5,pos=0.7,label={[xshift=1.0cm, yshift=0.3cm]:$q\eta(l)=r_l(0)$}] (r) {};
	\draw[thin,dashed,-{latex}] (q) .. controls (2.3,2) .. (p) {};
	\draw[thin,dashed,-{latex}] (p) -- (tau)
	node[circle,fill,scale=0.5,pos=0.7,label=below:{$r_l(D)$}] (rlD) {};
	\draw[thin,dashed,-{latex}] (r) .. controls (10.5,1) .. (rlD) {};
	
	\draw (15,1.5) -- (2,1.5) -- (0,-0.5) -- (13,-0.5) -- cycle ; % draw a parallel set containing p and tau
	\draw (14.5,1) -- (q) -- (13.5,0); % draw a cone with tip q
	\end{tikzpicture}
	\caption{Strongly asymptotic geodesics get close at an exponential rate}
\end{figure}

\begin{lem}\label{strongly asymptotic geodesics}
	Let $q \in \X$ and let $\eta \in \partial X$ be $(\alpha_0,\tau)$-regular. Let $P=P(\hat{\tau},\tau)$ be a parallel set with $d(q,P) \le D$, and let $p \in P$ be the unique point on the horocyle $H(q,\tau)$. Then for all $l \ge 0$ the geodesic rays $p\eta$ and $q\eta$ satisfy 
	$$ d(p\eta(l),q\eta(l)) \le De^{\kappa_0D-\alpha_0l} .$$
\end{lem}

It is possible to prove (a slightly weaker variation of) Lemma \ref{strongly asymptotic geodesics} as a limiting case of Lemma \ref{cone rotation}, or to construct a curve in $N_\tau$ in a similar way as we constructed a curve in $K$ in Lemma \ref{cone rotation}. However, we give a direct proof here using the generalized Iwasawa decomposition, see Section \ref{iwasawa decomposition}.

\begin{proof}
	We may assume that $d(q,P)=D$. By abuse of notation let $q \colon [0,D] \to \X$ be the unit speed geodesic segment from $q$ to its nearest point $\bar{q} \in P$. Let $G=N_\tau A_\tau K$ be the generalized Iwawsawa decomposition associated to $p$ and $\tau$, see Section \ref{iwasawa decomposition}. Since $N_\tau \times A_\tau \to M, (u,a) \mapsto uap$ is a diffeomorphism, we may write $q(s) = u(s)a(s)p$ for unique curves $u \colon [0,D] \to N_\tau$ and $a \colon [0,D] \to A_\tau$. Note that $u(D)=1=a(0)$, since horocycles at $\tau$ meet parallel sets $P(\hat{\tau},\tau)$ in exactly one point. 
	
	Writing $c_t(s) = C(s,t) = u(s)a(t)p$ we have $q(s) = C(s,s)=c_s(s)$ so
	$$ \dot{q}(s_0) = \eval{\pdv{C}{s}}_{s_0,s_0} + \eval{\pdv{C}{t}}_{s_0,s_0} = \dot{c}_{s_0}(s_0) + \eval{\pdv{C}{t}}_{s_0,s_0} $$
	and these vectors are orthogonal, so each has norm bounded by $1$. The curve $t \mapsto a(t)p$ has speed bounded by $1$ since 
	$$ \eval{\pdv{C}{t}}_{s_0,t_0} = \dd{u(s_0)} \eval{ \dv{t} a(t)p }_{t=t_0} ,$$
	so $d(p, a(t)p) \le t \le D$. We write $\dot{u}(s) = \dd{l_{u(s)}} U_s$ and use Lemmas \ref{metric and Bp}, \ref{Bp comparison} and \ref{right translated curves} to obtain 
	$$ \abs{ \dot{c}_0(s)} = \abs{ \evp U_s} = \frac{1}{\sqrt{2}} \pnorm{U_s} \le \frac{1}{\sqrt{2}} e^{\kappa_0 d(p,a(s)p)} \abs{ U_s}_{B_{a(s)p}} = e^{\kappa_0 d(p,a(s)p)} \abs{ \dot{c}_s(s)} \le e^{\kappa_0s} .$$
	
	We next need to push this horocyclic curve towards $\tau$ and check that the length shrinks by at least $e^{-\alpha_0 l}$. Let $X \in \mfp$ be the unit vector so that $q\eta(t) = u(0)e^{tX}p$. By abuse of notation define the curve $r_t(s)=u(s)e^{tX}p$ from $q \eta(t)$ to $p \eta(t)$ and note that $r_l(0)=u(0)e^{lX}p=q\eta(l)$. We've shown that the speed of $r_0 = c_0$ is at most $e^{\kappa_0s}$, and we may conclude after we show that 
	$$ \abs{ \dot{r_t}(s) } \le e^{-\alpha_0 t} \abs{ \dot{r_0}(s) } $$
	in the next paragraph.
	
	Define curves $U_\alpha(s) \in \mfg_\alpha$ by  $\dot{u}(s) = (\dd{l_{u(s)}})_1 \sum_{\alpha \in \Lambda_\tau^+}U_\alpha(s)$ and using Lemma \ref{right translated curves} write 
	\begin{align*}
	\abs{\dot{r_t}(s)}_{T_{r_t(s)}\X} & = \abs{ \evp \Ad(e^{-tX}) \sum_{\alpha \in \Lambda_\tau^+} U_\alpha (s) }_{T_p\X} \\
	& = \abs{ \evp \sum_{\alpha \in \Lambda_\tau^+} e^{-t\alpha(X)} U_\alpha (s) }_{T_p\X} \\
	& = \frac{1}{\sqrt{2}}\pnorm{ \sum_{\alpha \in \Lambda_\tau^+} e^{-t\alpha(X)} U_\alpha (s) } \\
	& \le \frac{1}{\sqrt{2}} e^{-\alpha_0 t}\pnorm{ \sum_{\alpha \in \Lambda_\tau^+} U_\alpha (s) } \\
	& = e^{-\alpha_0 t} \abs{ \dot{r_0}(s) }_{T_{c(s)}\X}  
	\end{align*}
	Integrating this inequality bounds the length of $r_l$ by $D e^{\kappa_0D - \alpha_0 l}$ and completes the proof.
\end{proof}

It is possible to give a slightly stronger upper bound in Lemma \ref{strongly asymptotic geodesics}, but the improvement would be inconsequential when we apply this Lemma in Section \ref{local to global section} while making the already cumbersome statements even harder to read.
%e.g. the upper bound  $\kappa_0^{-1}(e^{\kappa_0D}-1)e^{-\alpha_0 l}$ is easy

% I added a condition. It might be possible to avoid this condition by some continuity argument and that \tau is isolated in the tau_mod boundary of P(tau',tau). If I can't, then I need to rearrange so that the lemmas come in the right order.
% I rearranged the lemmas. If I want to look for a continuity argument, I could take a look at the KLP Morse Lemma paper. They have an argument about longitudinal directions being open and closed that seems potentially helpful

The following Lemma is a quantified version of Lemma 2.40 in \cite{KLP14}. 
\begin{lem}\label{points far along asymptotic Weyl cones}
	Let $p,q,x \in \X$ with $pq$ an $(\alpha_0,\tau)$-regular geodesic segment and $d(p,q) \ge l$ and $d(p,x) \le D$. If 
	$$ \alpha_0 - \frac{D(\alpha_0+\kappa_0)}{l-D} \ge \alpha_0' \quad \text{ and } \quad \frac{1}{\alpha_0'\zeta_0}\frac{D}{l} \le \frac{\zeta_0^2}{\kappa_0^2} $$
	then 
	$$ d(q, V(x,\st(\tau),\alpha_0')) \le De^{\kappa_0D-\alpha_0l} .$$
\end{lem}

\begin{proof}
	Let $\eta \in \ost(\tau)$ such that $pq(+\infty) = \eta$. Let $y$ be the unique point in the intersection $P(S_x\tau,\tau) \cap H(p,\tau)$. The point $q'$ on the image of $y\eta$ such that $\vec{d}(y,q')=\vec{d}(p,q)$ satisfies $d(q,q') \le D e^{\kappa_0D - \alpha_0l}$ by Lemma \ref{strongly asymptotic geodesics}. We will prove the Lemma by showing that $xq'$ is $(\alpha_0',\tau)$-regular.
	
	Choose chambers $\sigma,\sigma'$ so that $yq' \in V(y,\sigma)$ and $xq' \in V(x,\sigma')$. Then there is a unique (restricted) isometry $g\colon V(y,\sigma) \to V(x,\sigma')$ by Theorem \ref{Vmod welldefined} and
	$$ d(gq',q')= \abs{\vec{d}(x,gq') - \vec{d}(x,q')} = \abs{\vec{d}(y,q')-\vec{d}(x,q')} \le d(x,y) \le D .$$
	Now both $q'$ and $gq'$ lie in the same Euclidean Weyl cone $V(x,\sigma')$ with $d(q',gq') \le D$ and the geodesic segment from $x$ to $gq'$ is length at least $l$ and $(\alpha_0,\tau_{mod})$-regular, so Lemma \ref{regular projections} implies that $xq'$ is $(\alpha_0',\tau_{mod})$-regular. 
	
	We conclude by showing that $xq'$ is $\tau$-regular. By Lemma \ref{zeta projection} and Lemma \ref{euclidean angle} we have that $\angle^\zeta_{q'}(x,y) \le \frac{1}{\alpha_0'\zeta_0}\frac{D}{l} \le \frac{\zeta_0^2}{\kappa_0^2}$, so $\angle^\zeta_{q'}(x,\tau) \ge \pi - \varepsilon(\zeta_{mod})$ by Lemma \ref{angle implies antipodal}. Since $S_x\tau=S_{q'}\tau$ is the unique antipode of $\tau$ in the boundary of $P(S_x\tau,\tau)$, it follows that $xq'$ is $\tau$-regular. 
\end{proof}

\subsection{Projecting midpoints to Weyl cones}

We combine the previous Lemmas \ref{cone rotation}, \ref{strongly asymptotic geodesics} and \ref{points far along asymptotic Weyl cones} to show that a long regular geodesic segment in a bounded neighborhood of a Weyl cone has its midpoint arbitrarily close to the Weyl cone. 

\begin{cor}\label{midpoint projection}
	Let $p,q,x \in \X$ with $pq$ an $(\alpha_0,\tau_{mod})$-regular geodesic segment with midpoint $m$, let $\tau \in \Flagt$ and let $V=V(x,\st(\tau))$. Assume that $d(p,x) \le D, d(q,V) \le D$ and $d(p,q)\ge 2l$. Suppose that
	\begin{enumerate}
		\item\label{midpoint 1} $$ \alpha_0 - \frac{2D(\alpha_0+\kappa_0)}{l-2D} \ge \alpha_0' >0,  $$
		\item\label{midpoint 2} $$ \frac{1}{2}(e^{4\kappa_0D}-1)[\sinh(\alpha_0'(2l-2D))]^{-2} \le 3e^{4\kappa_0D}, \text{ and }$$
		\item\label{midpoint 3} $$ \frac{2}{\alpha_0'\zeta_0}\frac{D}{l} \le \frac{\zeta_0^2}{\kappa_0^2} $$
	\end{enumerate}
	then 
	$$ d(m,V(x,\st(\tau),\alpha_0')) \le 5De^{2\kappa_0D-\alpha_0l} .$$
\end{cor}

% This Lemma has been corrected. 

\begin{proof}
	Since $d(q,V) \le D$ and the Hausdorff distance from $V$ to $V(p,\st(\tau))$ is at most $D$, we have $d(q,V(p,\st(\tau)))\le2D$. We may now apply Lemma \ref{cone rotation} together with assumptions \ref{midpoint 1} and \ref{midpoint 2} to see that there exists $m' \in V(p,\st(\tau),\alpha_0)$ with $d(m,m') \le 4De^{2\kappa_0D-\alpha_0l}$ and $d(p,m')=d(p,m)\ge l$. 
	
	Assumptions \ref{midpoint 1} and \ref{midpoint 3} allow us to apply Lemma \ref{points far along asymptotic Weyl cones} to see that $d(m',V(x,\st(\tau),\alpha_0')) \le D e^{\kappa_0D-\alpha_0l}$.	By the triangle inequality, 
	$$ d(m,V(x,\st(\tau),\alpha_0')) \le d(m,m')+d(m',V(x,\st(\tau),\alpha_0')) \le 4De^{2\kappa_0D-\alpha_0l} + D e^{\kappa_0D-\alpha_0l} \le 5D e^{2\kappa_0D-\alpha_0l} .$$	
\end{proof}

% comment on the fact that a better inequality is possible? It is obvious...

\subsection{Simplex displacement after a short flow}

Recall that we have fixed a model type $\zeta=\zeta_{mod}$ spanning $\tau_{mod}$, see Definition \ref{spanning} and Section \ref{section regularity ideal}. 

\begin{lem}\label{simplex displacement}
	For any point $p \in \X$, simplex $\tau \in \Flag(\tau_{mod})$ and transvection vector $X \in \mfp$, it holds that
	$$ \sin \frac{1}{2} \angle_p^\zeta(\tau,e^X\tau)\le \frac{\kappa_0}{2} \pnorm{X} .$$
\end{lem}

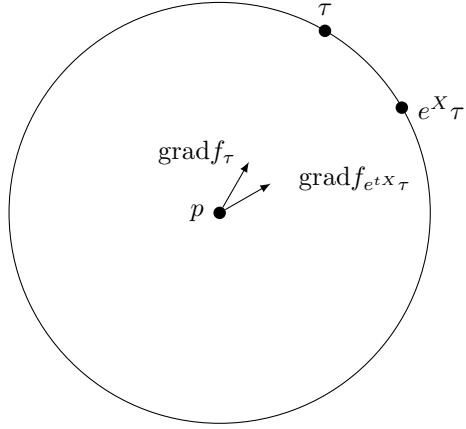
\begin{figure}[h]
	\centering
	\begin{tikzpicture}[scale=0.7]
	\node[circle,scale=0.5,fill=black,label=left:$p$] (p) at (0,0) {};
	\draw (p) circle (4);
	\node[circle,scale=0.5,fill=black,label=above:{$\tau$}] (t) at (2,3.46) {};
	\node[circle,scale=0.5,fill=black,label=right:{$e^X\tau$}] (ext) at (3.46,2) {};
	\node[label=left:{$\grad f_{\tau}$}] (tv) at (0.66,1.15) {};
	\node[label=right:{$\grad f_{e^{tX}\tau}$}] (extv) at (1.15,0.66) {};
	\draw[-{latex}] (p) -- (tv);
	\draw[-{latex}] (p) -- (extv);		
	\end{tikzpicture}
	\caption{Simplex displacement}
\end{figure}

\begin{proof}
	Denote by $f_\tau$ the Busemann function associated to the ray from $p$ to $\zeta(\tau)$ and write $\grad f_\tau$ for its gradient. Then 
	$$ \angle_p^\zeta(\tau,e^X\tau) = \angle_p(\grad f_\tau, \grad f_{e^X \tau}) $$
	and 
	$$ \sin \frac{1}{2} \angle_p(\grad f_\tau, \grad f_{e^X \tau} ) = \frac{1}{2} d_{T_p\X} (\grad f_\tau, \grad f_{e^X\tau} ) .$$
	
	Let $Z \in \mfp$ be the unit vector so that $\evp Z = (\grad {f_\tau})_p .$ Decompose $X=K+Y$ according to the generalized Iwasawa decomposition $\mfg = \mfk + \mfa_\tau + \mfn_\tau$ so that flowing by $Y$ fixes $\tau$ and therefore commutes with $\grad f_\tau$, and flowing by $K$ fixes $p$, see Section \ref{iwasawa decomposition}. We may write $X = A + \sum_{\alpha \in \Lambda^+} (-X_\alpha + \vartheta_p X_\alpha)$ and $K = \sum_{\alpha \in \Lambda^+} (X_\alpha + \vartheta_p X_\alpha)$ so $\pnorm{K} \le \pnorm{X}$. At $p$ we have
	\begin{align*}
	\eval{\dv{t} \left(\grad{f_{e^{tX}\tau}} \right)_p}_{t=0}  & =  \eval{\dv{t} \left( \left(e^{tX}\right)_\ast \grad{f_\tau} \right)_p }_{t=0} \\
	& = (\lied_{-X^\ast} \grad{f_\tau})_p \\
	& = [-X^\ast, \grad{f_\tau}]_p \\
	& = [(-X + Y)^\ast, \grad{f_\tau}]_p \\
	& = [-K^\ast,\grad{f_\tau} ]_p \\
	& = (\lied_{-K^\ast} \grad f_\tau)_p \\
	& = \lim_{t \to 0} \frac{ (\d{e^{tK}}) (\grad f_\tau)_{e^{-tK}p} - (\grad f_\tau)_p}{t} \\
	& = \lim_{t \to 0} \frac{ (\d{e^{tK}}) (\grad f_\tau)_{p} - (\grad f_\tau)_p}{t} \\
	& = \lim_{t \to 0} \frac{ (\d{e^{tK}}) \evp Z - \evp Z }{t} \\
	& = \lim_{t \to 0} \frac{  \evp \Ad(e^{tK}) Z - \evp Z }{t} \\
	& = \evp [K,Z]_{\mfg}
	\end{align*}
	Since we assumed nothing about the relationship of $X$ and $\tau$ we see that for all $t' \in [0,1]$, 
	$$ \Bigg\lvert \eval{\dv{t} \left(\grad{f_{e^{tX}\tau}} \right)_p}_{t=t'} \Bigg\rvert = \Bigg\lvert \eval{\dv{t} \left( (e^{tX})_\ast \grad{f_{e^{t'X}\tau}} \right)_p }_{t=0} \Bigg\rvert \le \pnorm{[K,Z]} \le \kappa_0 \pnorm{K} \le \kappa_0 \pnorm{X} $$
	where we used Lemma \ref{ad bound} in the second inequality.
	Finally we obtain
	$$ \abs{\grad{f_\tau}-\grad{f_{e^X \tau}}}_{T_p \X} \le \int_0^1 \abs{ \dv{t} \grad{f_{e^{tX}\tau}} }_{T_p \X} \dd{t} \le  \kappa_0 \pnorm{X} $$
	which completes the proof. 
\end{proof}

\subsection{The distance to a parallel set bounds the $\zeta$-angle}

\begin{cor}\label{distance to angle}
	Let $p,q$ be points $\X$ and $\tau,\tau' \in \Flag(\tau_{mod})$. If $d(p,q) \le \frac{2}{\kappa_0}$ then
	$$ \abs{\angle_p^\zeta(\tau,\tau') - \angle_q^\zeta(\tau,\tau')} \le 4 \sin^{-1} \left(\frac{\kappa_0}{2} d(p,q) \right) .$$
\end{cor}

\begin{proof}
	Write $q=e^{-X}p$ for $X \in \mfp$. We use that $\zeta$-angles are $G$-invariant, the triangle inequality for quadruples in $(\Flag(\tau_{mod}),\angle_p^\zeta)$ and the simplex displacement estimate given by Lemma \ref{simplex displacement}.
	$$ \abs{\angle_p^\zeta(\tau,\tau') - \angle_q^\zeta(\tau,\tau')} = \abs{\angle_p^\zeta(\tau,\tau') - \angle_p^\zeta(e^X\tau,e^X\tau')} \le \angle_p^\zeta(\tau,e^X\tau)+\angle_p^\zeta(\tau',e^X\tau') \le 4 \sin^{-1} \left( \frac{\kappa_0}{2} \pnorm{X} \right) .$$
	Since $\pnorm{X}=d(p,q)$ we are done.
\end{proof}

We will often apply Corollary \ref{distance to angle} in the following form. This result is a quantified version of Lemma 2.43.(i) in \cite{KLP14}. 

\begin{cor}\label{distance to angles near pi}
	Let $\tau_+,\tau_-$ be antipodal simplices in $\Flag(\tau_{mod})$ and let $P=P(\tau_-,\tau_+)$ be the parallel set joining them. Let $p$ be any point in $\X$ such that $ d(p,P) \le \frac{2}{\kappa_0} $. Then
	$$ \angle_p^\zeta(\tau_-,\tau_+) \ge \pi - 4 \sin^{-1} \left(\frac{\kappa_0}{2} d(p,P) \right) .$$
\end{cor}

\begin{proof}
	Since $\angle_q^\zeta(\tau_-,\tau_+)=\pi$ for any $q \in P$, and in particular the projection of $p$ to $P$, the assertion follows immediately from Corollary \ref{distance to angle}.	
\end{proof}

\subsection{The $\zeta$-angle bounds the distance to the parallel set}

We continue to work with a fixed $(\zeta_0,\tau_{mod})$-spanning type $\zeta=\zeta_{mod}$ and from now on assume that $\zeta$ is $\iota$-invariant, see the discussion after Theorem \ref{Vmod welldefined}. The next lemma complements Corollary \ref{distance to angles near pi}: when the $\zeta$-angle at $q \in \X$ between simplices $\tau_\pm \in \Flagt$ is near $\pi$, the point $q$ is near the parallel set $P(\tau_-,\tau_+)$. In the proof we use the fact that a vector field $X$ is Killing (if and) only if for all vector fields $V,W$ on $\X$, we have 
$$ X\langle V,W \rangle = \langle [X,V],W \rangle + \langle V, [X,W] \rangle ,$$
see \cite[9.25]{ONe83}. The following result is a quantified version of Lemma 2.43.(ii) in \cite{KLP14}. 

% Could add a paragraph explaining why this approach "has to" work. E.g. Busemann functionss are asymptotically linear and also X being homogeneous.
% Decided not to for now. It's kind of a pain to make these reasons rigorous.

\begin{lem}\label{angle to distance}
	Let $\tau_\pm \in \Flag(\tau_{mod})$ and let $q \in \X$. If $\delta \le \frac{\zeta_0^2}{2\kappa_0^2}$ and $\angle_q^\zeta(\tau_-,\tau_+) \ge \pi - \delta$ then $\tau_\pm$ are antipodal and $d(q,P(\tau_-,\tau_+)) \le \delta/\zeta_0$.
\end{lem}

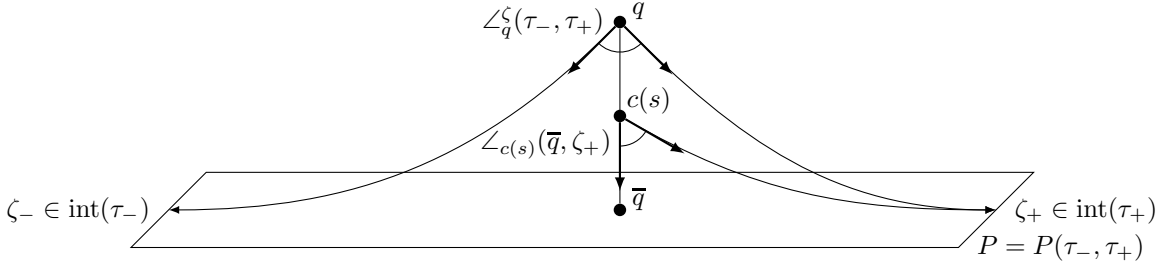
\begin{figure}[h]
	\centering
	\begin{tikzpicture}
		\node[circle,scale=0.5,fill=black,label=left:{$\angle_q^\zeta(\tau_-,\tau_+)$}] (q) at (0.5,3) {}; % draw q
		\node[circle,scale=0.5,fill=black,label={[xshift=7,yshift=-6]:{$\overline{q}$}}] (bq) at (0.5,0.5) {}; % draw bar{q}
		\node[label=right:{$P=P(\tau_-,\tau_+)$}] (P) at (5,0) {}; %label the parallel set
		\draw (-6,0) -- (P.center) -- (6,1) -- (-5,1) -- cycle; % draw the parallel set 
		\node[label=right:{$\zeta_+ \in \interior(\tau_+)$}] (z) at (5.5,0.5) {}; %label the ideal point
		\node[label=left:{$\zeta_- \in \interior(\tau_-)$}] (zm) at (-5.5,0.5) {}; %label the other ideal point	
		\node[label={[xshift=7,yshift=-6]:{$q$}}] at (q) {};
		\draw (bq) -- (q) node[circle,scale=0.5,fill=black,pos=0.5,label={[xshift=-28,yshift=-22]:{$\angle_{c(s)}(\overline{q},\zeta_+)$}}] (cs) {}; % draw the geodesic segment
		\node[label={[xshift=11,yshift=-6]:{$c(s)$}}] at (cs) {};
		\draw[-{latex}] (q) to [out=-45,in=180] (z.center); % draw the ray from q to z
		\draw[-{latex}] (q) to [out=-135,in=0] (zm.center); % draw the ray from q to zm
		\draw[-{latex}] (cs) to [out=-30,in=180] (z.center);% draw the ray from c(s) to z
		\draw[-{latex},thick] (q) --+ (0.7,-0.7) {}; %node[label=right:{$(\grad f_+)_q$}] {}; % draw the tangent vector at q
		\draw[-{latex},thick] (q) --+ (-0.7,-0.7) {}; %node[label=left:{$(\grad f_-)_q$}] {}; % draw the tangent vector at q
		\draw[-{latex},thick] (cs) --+ (0.86,-0.5) {}; %node[label=right:{$(\grad f_+)_{c(s)}$}] {}; % draw the tangent vector at c(s)
%		\draw[-{latex},thick] (q) --+ (0,-1) node[label=left:{$-X^\ast_q$}] {}; % draw the tangent vector at q
		\draw[-{latex},thick] (cs) --+ (0,-1) {}; %node[label=left:{$-X^\ast_{c(s)}$}] {}; % draw the tangent vector at c(s)
		\draw (q) +(225:0.4) arc (225:315:0.4); % draw the angle
		\draw (cs) +(-90:0.4) arc (-90:-30:0.4); % draw the angle
	\end{tikzpicture}
	\caption{The $\zeta$-angle at $q$ bounds the distance to $P$}
\end{figure}

% Had to reprove this, got the bound 2\kappa_0^2. Should probably be \kappa_0^2 or less for sl(n,R). Comment somewhere that this is probably not optimal? 
% I decided no, because I didn't actually prove that bound for sl(n,R).

\begin{proof}
	Since $\angle_q^\zeta(\tau_-,\tau_+) \ge \pi - \frac{\zeta_0^2}{2\kappa_0^2} > \pi - \frac{\zeta_0^2}{\kappa_0^2}$, Lemma \ref{angle implies antipodal} implies that the simplices $\tau_-,\tau_+$ are antipodal. 
	
	Write $\zeta_\pm$ for the unique ideal points $\tau_\pm$ of type $\zeta$, and choose Busemann functions $f_\pm$ at $\zeta_\pm$. For all $p \in \X$ we have $\cos \angle_p^\zeta(\tau_-,\tau_+)=\cos \angle_p(\zeta_-,\zeta_+) = \langle \grad f_-, \grad f_+ \rangle_p$. Let $\bar{q} \in P=P(\tau_-,\tau_+)$ be the nearest point on $P$ to $q$, and let $X \in \mfp_{\bar{q}}$ such that $c(t)=e^{tX}\bar{q}$ is the unit-speed geodesic from $\bar{q}$ to $q$. Either $\angle_q(\zeta_-,\bar{q}) \ge \frac{\pi}{2} - \frac{\delta}{2}$ or $\angle_q(\bar{q},\zeta_+) \ge \frac{\pi}{2} - \frac{\delta}{2}$, so without loss of generality we may assume the second inequality holds. Let $f \colon (-\infty,\infty) \to [-1,1]$ be defined by $f(s) = \langle -X^\ast,\grad f_+ \rangle_{c(s)}$ and note that $f(s)=\cos \angle_{c(s)}(\bar{q},\zeta_+)$ for all $s >0$. We first show that $f'(s) \ge 0$ for all $s$, so $f$ is (weakly) monotonic.
	
	At the point $c(s)$, we have $X \in \mathfrak{p}_{c(s)}$ since $X$ is a transvection along $c$. The point $c(s)$ together with a fixed choice of chamber containing $\tau_+$ allows us to decompose $X$ according to the restricted root space decomposition. Suppressing the dependence on $s$, we have $X = A + \sum_{\alpha \in \Lambda^+} -X_\alpha + \vartheta X_\alpha$. Then for $K = \sum_{\alpha \in \Lambda^+} X_\alpha + \vartheta X_\alpha$ and the unit vector $Z \in \mfp_{c(s)}$ pointing to $\zeta_+$ we see that
\begin{align*}
	f'(s) & = X^\ast \langle - X^\ast , \grad f_+ \rangle_{c(s)} \\ % X restricts to \dot c
	& = \langle -X^\ast, [X^\ast,\grad f_+] \rangle_{c(s)} \\ % X is Killing 
	& = \langle -X^\ast, [K,Z]^\ast \rangle_{c(s)} \\ %simplex displacement 
	& = B( -X, -[K,Z]_{\mfg} ) \\ % \mfg is killing fields 
	& = B(A + \sum_{\beta \in \Lambda^+} -X_\beta + \vartheta X_\beta , \sum _{\alpha \in \Lambda^+} \alpha(Z) (X_\alpha - \vartheta X_\alpha) ) \\ % dp: (\mfp,B) \to (T_pM,\langle,\rangle) is an isometry, and apply roots.
	& =  \sum_{\alpha \in \Lambda^+} \alpha(Z) B(X_\alpha-\vartheta X_\alpha,X_\alpha-\vartheta X_\alpha) \\ % \mfp_\alpha decomp is B_p orthogonal 
	& \ge \zeta_0 \sum_{\alpha \in \Lambda_\tau^+} \abs{-X_\alpha+ \vartheta X_\alpha}_B^2.
\end{align*}
	The third line follows from the reasoning in the proof of Lemma \ref{simplex displacement}. This calculation shows that $f'(s) \ge 0 $ for all $s$. Moreover, since $X^\ast$ is orthogonal to $P(\bar{q},\tau)$ at $s=0$, we have $1=\abs{X^\ast_{\bar{q}}}^2 = \sum_{\alpha \in \Lambda_\tau^+} \abs{-X_\alpha + \vartheta X_\alpha}_B^2$, so $f'(0) \ge \zeta_0$. 
	
	We next bound the norm of 
	\begin{align*}
		f''(s) & = X^\ast (X^\ast \langle -X^\ast, \grad f_+ \rangle )_{c(s)}  \\ 
		& = \langle -X^\ast , [X^\ast,[X^\ast, \grad f_+]] \rangle_{c(s)} \\
		& = \langle -X^\ast , [X^\ast,[K^\ast, \grad f_+]] \rangle_{c(s)} \\
		& = \langle -X^\ast , [K^\ast,[X^\ast, \grad f_+]] \rangle_{c(s)} - \langle X^\ast , [[K^\ast,X^\ast], \grad f_+] \rangle_{c(s)} \\
		& = \langle -X^\ast , [K^\ast,[K^\ast, \grad f_+]] \rangle_{c(s)} + \langle X^\ast , [K'^\ast, \grad f_+] \rangle_{c(s)} \\
		& = B_{c(s)}( -X , [K,[K, Z]] ) + B_{c(s)} (X , [K', Z] ) \\
		& =  B_{c(s)}( [K,X] , [K, Z] ) + B_{c(s)} ([X,Z] , K' )
	\end{align*}
	where $[K,X]=K'+A'+N'$ according to the KAN decomposition for $c(s)$ and $\tau_+$. We get the bound
	\begin{align*}
		|f''(s)| & = | B_{c(s)}( [K,X] , [K, Z] ) + B_{c(s)} ([X,Z] , K' )| \\
		& \le |B_{c(s)}( [K,X] , [K, Z] )| + |B_{c(s)} ([X,Z] , K' )| \\
		& \le \abs{[K,X]}_{B_{c(s)}} \abs{[K,Z]}_{B_{c(s)}} + \abs{[X,Z]}_{B_{c(s)}} \abs{K'}_{B_{c(s)}} \\
		& \le 2 \kappa_0^2 \\
	\end{align*}
	by applying Lemma \ref{ad bound}.

	Since $f'(0) \ge \zeta_0$ and $\abs{f''(s)} \le 2 \kappa_0^2$, we have $f(s) \ge s\zeta_0 - \kappa_0^2 s^2$. Since $f$ is monotonic, if $s \ge \frac{\zeta_0}{2\kappa_0^2}$ then $f(s) \ge f\left( \frac{\zeta_0}{2\kappa_0^2} \right) \ge \frac{ \zeta_0^2}{4\kappa_0^2}.$ On the other hand, if $s \le \frac{\zeta_0}{2\kappa_0^2}$ we have $f(s) \ge \zeta_0 s - \kappa_0^2 \left(\frac{\zeta_0}{2\kappa_0^2}\right)s \ge \frac{1}{2} \zeta_0 s.$ This implies
	$$ \frac{1}{2} \zeta_0 d(q,P) \le f(d(q,P)) = \cos \angle_q^\zeta(\bar{q},\tau_+)  \le \cos( \frac{\pi}{2} - \frac{\delta}{2} ) = \sin (\frac{\delta}{2}) \le \frac{\delta}{2} $$
	unless $d(q,P) > \frac{\zeta_0}{2\kappa_0^2}$, which yields $\frac{\zeta_0^2}{2\kappa_0^2} < \delta$ and contradicts our assumption. 
\end{proof}

\section{Quantified local-to-global principle}\label{local to global section}

In this section we augment the theorems of \cite[Section 7]{KLP14} with quantitative estimates. We obtain a precise version of the local-to-global principle which allows us to perturb known Anosov representations by a definite amount, producing new Anosov representations in Section \ref{section applications}. 

In rank one, local quasigeodesics of sufficiently good quality are global quasigeodesics, as a consequence of the Morse lemma. The Morse Lemma fails in the Euclidean plane, hence in higher rank, so we must use \textit{Morse quasigeodesics} as defined in \cite{KLP14}. The strategy here, as in \cite{KLP14}, is to show that local Morse quasigeodesics of sufficiently good quality have straight and spaced midpoint sequences which are then globally Morse quasigeodesics. First we give an explicit local criteria for a sequence to be a Morse quasigeodesic.

\subsection{Sufficiently straight and spaced sequences are Morse quasigeodesics}\label{section straight and spaced sequences}

We recall some definitions from \cite{KLP14}. A sequence of points $(x_n)$ in $\X$ is $(\alpha_0,\tau_{mod},\epsilon)$\textit{-straight} if each geodesic segment $x_nx_{n+1}$ is $(\alpha_0,\tau_{mod})$-regular and if 
$$ \angle_{x_n}^\zeta(x_{n-1},x_{n+1}) \ge \pi- \epsilon $$
for all $n$. The sequence is $s$-\textit{spaced} if $d(x_n,x_{n+1})\ge s$ for all $n$. We say a sequence $(x_n)$ \textit{moves $\epsilon$-away from a simplex }$\tau$ if for all $n$
$$ \angle_{x_n}^\zeta(\tau,x_{n+1}) \ge \pi - \epsilon. $$

In this paper we are only interested in discrete sequences of points in $\X$. For us, a $(c_1,c_2,c_3,c_4)$-\textit{quasigeodesic} is a sequence $(x_n)$ (possibly finite, infinite, or biinfinite) such that 
$$ \frac{1}{c_1} \abs{N}-c_2 \le d(x_n,x_{n+N}) \le \abs{N}c_3+c_4 .$$
A sequence $(x_n)$ is $(c_1,c_2)$-\textit{coarsely spaced} (or \textit{lower-quasigeodesic}) if
$$ \frac{1}{c_1} \abs{N}-c_2 \le d(x_n,x_{n+N}).$$
Likewise $(x_n)$ is $(c_3,c_4)$-\textit{coarsely Lipschitz} (or \textit{upper-quasigeodesic}) if
$$ d(x_n,x_{n+N}) \le \abs{N}c_3+c_4 .$$

For an $(\alpha_0,\tau_{mod})$-regular segment $pq$, the $(\alpha_0,\tau_{mod})$-\textit{diamond} is the intersection
$$\diamondsuit_{\alpha_0}(p,q) \coloneqq V(p,\st(\tau(pq)),\alpha_0) \cap V(q,\st(\tau(qp)),\alpha_0) .$$

\begin{figure}[h]
	\centering
	\begin{tikzpicture}
		\node[circle,scale=0.5,fill=black,label=left:$p$] (p) at (-5,0) {};
		\node[circle,scale=0.5,fill=black,label=right:$q$] (q) at (5,0) {};
		\draw (5.5,2) -- (p) -- (4.5,-2);
		\draw (-4.5,2) -- (q) -- (-5.5,-2);
		\node (Vq) at (-3.5,2.2) {{$V(q,\st(\tau(qp)),\alpha_0)$}};
		\node (Vp) at (3.5,2.2) {{$V(p,\st(\tau(pq)),\alpha_0)$}};
		\node (d) at (0,0) {$\diamondsuit_{\alpha_0}(p,q)$} ;
		\draw (-7,-3) -- (6,-3) -- (7,3) -- (-6,3) -- cycle;
	\end{tikzpicture}
	\caption{The $(\alpha_0,\tau_{mod})$-diamond with endpoints $p$ and $q$}
\end{figure}
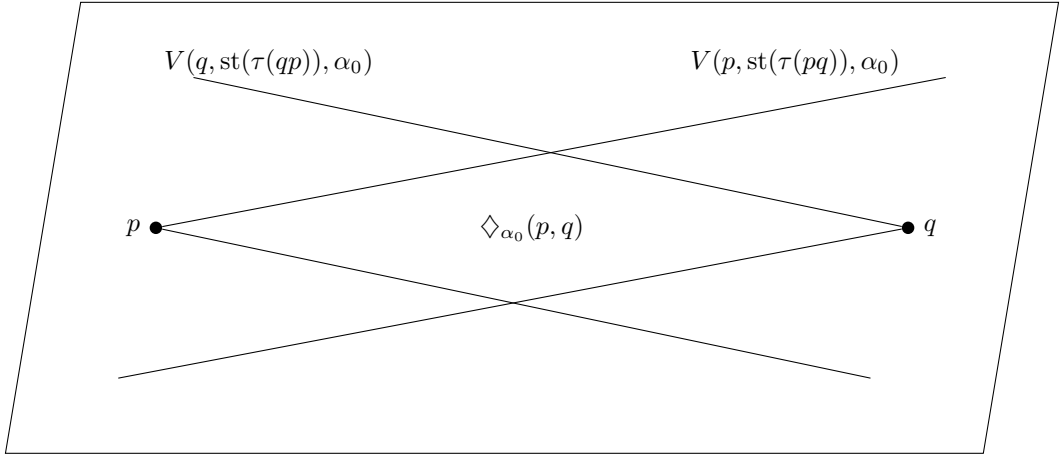

A quasigeodesic is $(\alpha_0,\tau_{mod},D)$-\textit{Morse} if for all $x_n,x_m$ there exists a diamond $\diamondsuit_{\alpha_0}(p,q)$ such that $d(p,x_n),d(q,x_m) \le D$ and for all $n \le i \le m$, $d(x_i,\diamondsuit) \le D$. In hyperbolic space, quasi-geodesics are automatically Morse by the Morse lemma. In higher rank symmetric spaces of noncompact type, the following theorem allows us to construct Morse quasigeodesics from sufficiently straight and spaced sequences. 

There are a few variations of the precise definition of Morse quasi-geodesic in the literature. The definition of Morse quasi-geodesic here is the same as that given in \cite[Definition 5.50]{KLP17}, except that we keep track of more constants in the definition of quasigeodesic. This is the same as \cite[Definition 7.14]{KLP14} except that we work with sequences rather than paths. Likewise \cite[Definition 6.13]{KL18b} defines paths to be Morse quasigeodesics when they satisfy a similar and equivalent, but not identical, property as the one we have given here (the constants will be different). %We warn the reader that the notion of Morse quasigeodesic here is not equivalent to the notion of $1$-dimensional Morse quasiflat given in \cite{HKS19,HKS20}. Rather, a maximal flat in a symmetric space is an example of a Morse quasiflat.
%%%: also Russel et al. also google "Morse quasigeodesic." Maybe even more uses.

Define the constant 
$$ c_0 \coloneqq \sum_{\alpha \in \Lambda_{\tau_{mod}}^+} \dim \mfg_\alpha,$$
equal to the codimension of any parallel set of type $\tau_{mod}$. The inequality $c_0 \ge 1$ always holds. Theorem \ref{straight and spaced sequences are morse} is a quantified version of Theorem 7.2 in \cite{KLP14}. 

\begin{thm}\label{straight and spaced sequences are morse}
	Fix $\alpha_{new}<\alpha_0,\delta$ and assume $\epsilon$ is small and $s$ is large. Precisely, we assume that:
	%3.1.2021: corrected
	\begin{enumerate}
		\item\label{1} $5 \epsilon \le \frac{\zeta_0^2}{2\kappa_0^2}$, so that we may apply the angle-to-distance estimate in Lemma \ref{angle to distance}; 
		\item\label{2} 
		$$\frac{\epsilon\kappa_0}{\zeta_0} e^{2 \kappa_0 \epsilon/\zeta_0-\alpha_0s} \le \sin \left( \frac{\epsilon}{4} \right) $$
		so that we may apply the distance-to-angle estimate Lemma \ref{distance to angles near pi};
		\item\label{3}
		$$ \frac{5\epsilon}{\zeta_0} \le \delta $$
		to control the distance from the sequence to the parallel set;
		\item\label{4}
		$$ \alpha_0 - \frac{2\delta(\alpha_0+\kappa_0)}{s-2\delta} \ge \alpha_{new} $$
		so that certain projections are $(\alpha_{new},\tau_{mod})$-regular by Lemma \ref{regular projections};
		\item\label{5} 
		$$ 2\epsilon + \sin^{-1} \left( \frac{2 \delta}{ \alpha_0 \zeta_0 s} \right) < \varepsilon(\zeta) $$
		so that certain simplices are antipodal, see Section \ref{section angles}. 
	\end{enumerate}

	Then every $(\alpha_0,\tau_{mod},\epsilon)$-straight $s$-spaced sequence $(x_n)$ in $\X$ is $\delta$-close to a parallel set $P(\tau_-,\tau_+)$ such that 
	$$ \overline{x}_{n \pm m} \in V(\overline{x}_n,\st(\tau_\pm),\alpha_{new}) $$
	for all $n$ and $m \ge 1$. It follows that the sequence is coarsely spaced: 
	$$ d(x_n,x_{n\pm m})\ge 2 \alpha_{new} \zeta_0 c_0 (s-2\delta)m-2\delta, $$
	and if $(x_n)$ is coarsely Lipschitz it is then a $(\alpha_{new},\tau_{mod},\delta)$-Morse quasigeodesic. 
\end{thm}

Our proof closely follows \cite[Section 7]{KLP14}, who prove the same theorem without the explicit assumptions \ref{1} through \ref{5} and without the explicit estimates we obtained in Section \ref{section estimates}. Note that the resulting sequence will always be $\frac{\zeta_0}{2\kappa_0^2}$-close to the parallel set, even if $\delta$ is chosen larger than that quantity. %corrected

\begin{proof}
	Step 1: Propagation cf.\ \cite[Lemma 7.6]{KLP14}. We show that for sufficiently straight and spaced sequences, the property of moving away from a simplex propagates along the sequence. 
	
	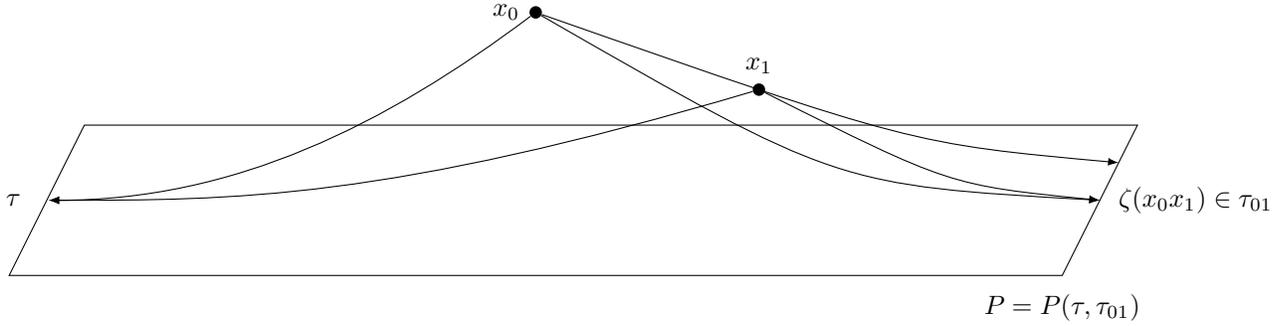
\begin{figure}[h]
		\centering
		\begin{tikzpicture}
		\node[label=below:{$P=P(\tau,\tau_{01})$}] (f) at (14,0) {}; % label the parallel set
		\node[circle,scale=0.5,fill=black,label=left:{$x_0$}] (x0) at (7,3.5) {}; % define x0
		\node[label=right:{$\zeta (x_0x_1) \in \tau_{01}$}] (t01) at (14.5,1) {}; % define tau01
		\node (x01) at (14.75,1.5) {}; %define endpoint of ray x0x1
		\node[label=left:{$\tau$}] (t) at (0.4,1) {}; %define tau
		\draw (0,0) -- (14,0) -- (15,2) -- (1,2) -- cycle; %draw parallel set
		\draw[-{latex[scale=2.5,length=2,width=3]}]
		 (x0) .. controls  (12,1.75) .. (x01.center) 
		 node[circle,fill,scale=0.5,pos=0.25,label=above:{$x_1$}] (x1) {};  % draw geodesic from x0 thru x1
		 \draw[-{latex[scale=2.5,length=2,width=3]}]
		 (x0) .. controls  (11,1.2) .. (t01.center); 
		\draw[-{latex[scale=2.5,length=2,width=3]}]
		 (x0) .. controls (5,2) and (3,1) .. (t);
		\draw[-{latex[scale=2.5,length=2,width=3]}]
		 (x1) .. controls (5,1) and (3,1) .. (t);
		 \draw[-{latex[scale=2.5,length=2,width=3]}]
		 (x1) .. controls (12.5,1.2) .. (t01.center); %draw ray x1 zeta(tau01)
		\end{tikzpicture}
		\caption{Propagation}
	\end{figure}
	
	Assume that for some simplex $\tau$ in $\Flagt$ we have $\angle_{x_0}^\zeta(\tau,x_1) \ge \pi- 2\epsilon.$ Since $2 \epsilon < \frac{\zeta_0^2}{2\kappa_0^2}$ by assumption \ref{1}, Lemma \ref{angle implies antipodal} implies that the simplex $\tau_{01}$ containing $x_0x_1(+\infty)$ is antipodal to $\tau$ and together they define a parallel set $P=P(\tau,\tau_{01})$. Moreover, assumption \ref{1} and our angle-to-distance estimate Lemma \ref{angle to distance} imply that $d(x_0,P) \le \frac{2\epsilon}{\zeta_0}$. By Lemma \ref{strongly asymptotic geodesics}, the geodesic ray from $x_0$ through $x_1$ gets arbitrarily close to $P$ and in particular 
	$$ d(x_1, P) \le \frac{2\epsilon}{\zeta_0} e^{2 \kappa_0 \epsilon/\zeta_0-\alpha_0 s} $$ 
	and by assumption \ref{2} and the distance-to-angle estimate Corollary \ref{distance to angles near pi} we have 
	$$ \angle_{x_1}^\zeta(\tau,\tau_{01}) \ge \pi- 4 \sin^{-1} \left( \frac{\epsilon \kappa_0}{\zeta_0} e^{2 \kappa_0 \epsilon/\zeta_0-\alpha_0 s} \right) \ge \pi - \epsilon $$
	which then implies that $\angle^\zeta_{x_1}(\tau,x_0)=\pi-\angle^\zeta_{x_1}(\tau,\tau_{01}) \le \epsilon$. Straightness and an application of the triangle inequality for $(S(T_{x_1}\X),\angle_{x_1})$ implies $\angle^\zeta_{x_1}(\tau,x_2) \ge \pi-2\epsilon$. By induction we have that $\angle_{x_n}^\zeta(\tau,x_{n+1})\ge \pi-2\epsilon$ for all $n \ge 1$. 
	
	Step 2: Extraction cf.\ \cite[Lemma 7.7]{KLP14}. We extract antipodal simplices that the sequence moves away/towards. It follows that the sequence stays near the corresponding parallel set.\footnote{The simplices are unique when the sequence is biinfinite, see \cite[7.19,5.15]{KLP14}, but this theorem also applies when the sequence is finite or a Morse quasiray.}
	
	For each $n$ define the compact subsets $C_n^\pm \subset \Flagt$
	$$ C_n^\pm \coloneqq \{ \tau_\pm \mid \angle_{x_n}^\zeta(\tau_\pm,x_{n\mp 1}) \ge \pi- 2\epsilon \}. $$
	Each of these is nonempty since $\angle_{x_n}^\zeta(x_{n\mp 1}x_n,x_{n\mp 1})=\pi$ implies $\tau(x_{n \mp 1} x_n)\in C_n^\pm$. By step 1, $C_n^- \subset C_{n+1}^-$ so there exists $\tau_- \in \cap_n C_n^-$. Similarly, there exists some $\tau_+ \in \cap_n C_n^+$. Straightness and the triangle inequality imply $\angle^\zeta_{x_n}(\tau_-,\tau_+)\ge \pi - 5 \epsilon$, and by assumption \ref{1} we have $5\epsilon \le \frac{\zeta_0^2}{2\kappa_0^2}$. Therefore the angle-to-distance estimate Lemma \ref{angle to distance} implies that $\tau_\pm$ are antipodal and define the parallel set $P=P(\tau_-,\tau_+)$ and moreover 
	$$ d(x_n,P) \le \frac{5\epsilon}{\zeta_0} \le \delta $$
	with the last inequality from assumption \ref{3}. 
		
	Step 3: Morseness cf.\ \cite[Lemma 7.9, Lemma 7.10, Corollary 7.13]{KLP14}. We verify that the sequence is a Morse quasi-geodesic. We have already shown the angles are straight enough to guarantee that the projection to $P$ is bounded. We show that projected rays land in nested cones; it follows that projecting further to the $\zeta$-ray yields a monotonic sequence which makes progress bounded away from zero. 

\begin{figure}[h]
	\centering
	\begin{tikzpicture}
		\node[label=below:{$P=P(\tau_-,\tau_+)$}] (P) at (14,0) {}; %label the parallel set
		\node[circle,scale=0.5,fill=black,label=above:{$x_n$}] (xn) at (2,4.4) {}; %define x_n
		\node[label=right:{$\xi$}] (xi) at (14.67,2) {}; % define xi
		\node[label=right:{$\zeta(\overline{x}_n \overline{x}_{n+1})$}] (taup) at (14.5,1.5) {};
		\node[circle,scale=0.5,fill=black,label=left:{$\overline{x}_n$}] (bxn) at (2,1.5) {}; %projection of x_n
		\draw[dashed,-{latex}] (bxn) -- (xi.center) 	node[circle,scale=0.5,pos=0.5,fill=black,label=below:{$\overline{x}_{n+1}$}] (bxn1) at (7,3.7) {};  %draw the ray bxn bxn+1 and define bxn+1
		\node[circle,scale=0.5,fill=black,label={[xshift=15,yshift=-1]:{$x_{n+1}$}}] (xn1) at ([yshift=2.68cm]bxn1) {}; %define x_n+1
		\draw (15,5.5) -- (xn) -- (14,4); %draw a cone with tip xn containing xn+1
		\draw[dashed,-{latex}] (xn) .. controls (8,2.5) .. (xi.center); % draw the ray from xn to \xi
		\node[label=right:{$\zeta(x_n x_{n+1})$}] (znn1) at (14.5,4.75) {};
		\draw[thick,dotted,-{latex}] (xn) -- (znn1.center);
		\draw[thick,dotted,-{latex}] (xn) .. controls (8,2) .. (taup.center);  %draw the ray x_n \xi
		\draw (0,0) -- (14,0) -- (15,3) -- (1,3) -- cycle; %draw a parallelogram to represent the parallel set
		\draw (14.8,2.8) -- (bxn) -- (13.8,0.2); %draw a cone
	\end{tikzpicture}
	\caption{The projection $\overline{x}_{n+1}$ lands in the Weyl cone $V(\overline{x}_n,\st(\tau_+),\alpha_{new})$}
\end{figure}

	By assumption \ref{4}, and Lemma \ref{regular projections}, we have that the projections $(\overline{x}_n)$ to $P$ are $(\alpha_{new},\tau_{mod})$-regular. Let $\xi$ be the ideal point corresponding to the ray $\overline{x}_n \overline{x}_{n+1}$. Since the rays $x_n \xi$ and $\overline{x}_n\xi$ are asymptotic, their Hausdorff distance is at most $d(x_n,\overline{x}_n) \le \delta$, so $x_{n+1}$ is at most $2\delta$ from $x_n \xi$. Then
	$$ \angle_{Tits}^\zeta(\tau_-,\xi) \ge \angle_{x_n}^\zeta(\tau_-,\xi) \ge \angle_{x_n}^\zeta(\tau_-,x_{n+1})-\angle_{x_n}^\zeta(x_{n+1},\xi) \ge \pi - 2  \epsilon - \angle_{x_n}^\zeta (x_{n+1},\xi). $$
	By Lemma \ref{zeta projection} and Lemma \ref{euclidean angle} we may guarantee that 
	$$ \sin \angle_{x_n}^\zeta (x_{n+1},\xi) \le \frac{1}{\alpha_0\zeta_0} \frac{2\delta}{s} $$
	so by assumption \ref{5} this Tits angle is within $\varepsilon(\zeta)$ of $\pi$, so $\zeta(\tau_-)$ is antipodal to $\zeta(\xi)$, but the only simplex in $\partial P$ antipodal to $\tau_-$ is $\tau_+$, so $\tau(\xi)=\tau_+$ and
	$$ \angle_{\overline{x}_n}^\zeta(\tau_-,\overline{x}_{n+1}) = \angle_{\overline{x}_n}^\zeta(\tau_-,\xi) = \pi .$$
	We know that $\overline{x}_n \overline{x}_{n+1}$ is $(\alpha_{new},\tau_{mod})$-regular and $\angle_{\overline{x}_n}^\zeta(\tau_-,\xi) = \pi$ and these two properties are equivalent to $\overline{x}_{n+1} \in V(\overline{x}_n,\st(\tau_+),\alpha_{new})$. Using the convexity of Weyl cones and induction, we get that for all $n$ and all $m \ge 1$
	$$ \overline{x}_{n \pm m} \in V(\overline{x}_n,\st(\tau_\pm),\alpha_{new}). $$
	
	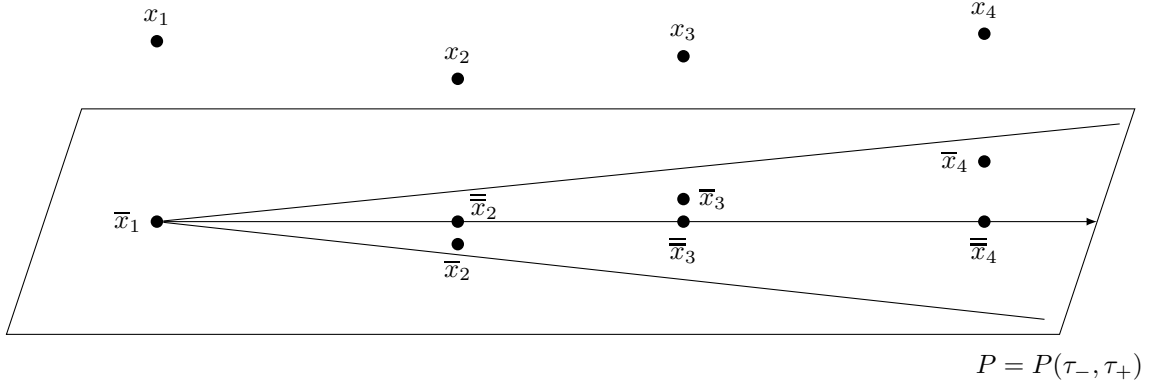
\begin{figure}[h]
		\centering
		\begin{tikzpicture}
		\node[label=below:{$P=P(\tau_-,\tau_+)$}] (f) at (14,0) {};
		\node[circle,scale=0.5,fill=black,label=above:{$x_1$}] (x1) at (2,3.9) {};
		\node[circle,scale=0.5,fill=black,label=above:{$x_2$}] (x2) at (6,3.4) {};
		\node[circle,scale=0.5,fill=black,label=above:{$x_3$}] (x3) at (9,3.7) {};
		\node[circle,scale=0.5,fill=black,label=above:{$x_4$}] (x4) at (13,4.0) {};
		\node[circle,scale=0.5,fill=black,label=left:{$\overline{x}_1$}] (bx1) at (2,1.5) {};
		\node[circle,scale=0.5,fill=black,label=below:{$\overline{x}_2$}] (bx2) at (6,1.2) {};
		\node[circle,scale=0.5,fill=black,label=right:{$\overline{x}_3$}] (bx3) at (9,1.8) {};
		\node[circle,scale=0.5,fill=black,label=left:{$\overline{x}_4$}] (bx4) at (13,2.3) {};
		\node[circle,scale=0.5,fill=black,label={[xshift=10,yshift=-5]:$\overline{\overline{x}}_2$}] (bbx2) at (6,1.5) {};
		\node[circle,scale=0.5,fill=black,label=below:{$\overline{\overline{x}}_3$}] (bbx3) at (9,1.5) {};
		\node[circle,scale=0.5,fill=black,label=below:{$\overline{\overline{x}}_4$}] (bbx4) at (13,1.5) {};
		\draw (0,0) -- (14,0) -- (15,3) -- (1,3) -- cycle; %draw a parallelogram to represent the parallel set
		\draw (14.8,2.8) -- (bx1) -- (13.8,0.2); %draw a cone
		\draw[-{latex}] (bx1) -- (14.5,1.5);  %draw a ray in the Pset
		\end{tikzpicture}
		\caption{Sufficiently straight and spaced sequences have monotonic projections to a geodesic ray}
	\end{figure}
	
	Finally, we want to show the sequence is coarsely spaced. The bound 
	$$ d(x_n,x_{n+m}) \ge 2 \alpha_{new} \zeta_0 c_0 (s-2\delta)m-2\delta $$
	will follow from 
	$$ d(\overline{x}_n,\overline{x}_{n+m}) \ge 2 \alpha_{new} \zeta_0 c_0 (s-2\delta)m .$$
	Indeed, the sequence $(\overline{x}_n)$ in $P$ is $(s-2\delta)$-spaced and has a monotonic projection $(\overline{\overline{x}}_n)$ to the geodesic line $\overline{x}_n \zeta(\tau_+)$ for any $n$ by the nestedness of Weyl cones. By \cite[2.14.5]{E96}, 
	$$ B(\zeta,\vec{d}(\overline{x}_n,\overline{x}_{n+1})) = \sum_{\alpha \in \Lambda} \alpha(\zeta)\alpha(\vec{d}(\overline{x}_n,\overline{x}_{n+1})) \dim \mfg_\alpha \ge 2 \alpha_{new} \zeta_0 d(\overline{x}_n,\overline{x}_{n+1}) \sum_{\alpha \in \Lambda_\tau^+} \dim \mfg_\alpha = 2 \alpha_{new} \zeta_0 c_0 d(\overline{x}_n,\overline{x}_{n+1}).$$
	It follows that the projection $\overline{ \overline{ x}}_{n+1}$ lies at least $2\alpha_{new}\zeta_0c_0(s-2 \delta)$ along the ray $\overline{x}_n\zeta$. 
\end{proof}

In the final step of the proof we used the regularity of the projections to obtain the linear lower-quasigeodesic constant. When the angular radius of $\sigma_{mod}$ with respect to $\zeta$ is strictly less than $\pi/2$, the linear lower-quasigeodesic bound can be chosen independent of the regularity. By \cite[Lemma 5.8]{KL18a}, this happens exactly when $\zeta$ is not contained in a factor of a nontrivial spherical join decomposition of $\sigma_{mod}$. In particular this is always possible when $\X$ is irreducible.

%corrected
\begin{rem}\label{remark straight to morse}
	To provide suitable auxiliary parameters to apply Theorem \ref{straight and spaced sequences are morse}, we may first choose $\epsilon$ small enough to satisfy assumptions \ref{1} and \ref{3} and then choose $s$ large enough to satisfy assumptions \ref{2},\ref{4} and \ref{5}. When we apply Theorem \ref{straight and spaced sequences are morse}	in Section \ref{section applications}, we will choose $\delta=\frac{\zeta_0}{2\kappa_0^2}$ and $\epsilon=\frac{\zeta_0^2}{10\kappa_0^2}$ and then find a large enough $s$ to satisfy the conditions of Theorem \ref{straight and spaced sequences are morse}. 
\end{rem}

\subsection{Morse quasigeodesics have straight and spaced midpoints}

In this section we show that Morse quasigeodesics of sufficiently good quality have straight and spaced midpoint sequences. 

\begin{defn}[{Cf.\ \cite[Definition 7.14]{KLP14}}]
For points $p,q$ in $\X$ we let $\midp(p,q)$ denote the midpoint of the geodesic segment $pq$. A sequence $(p_n)_{n=t_0}^{n=t_{max}}$ in $\X$ satisfies the $(\alpha_0,\tau_{mod},\epsilon,s,k)$-\textit{quadruple condition} if for all $t_1,t_2,t_3,t_4 \in [t_0,t_{max}] \cap \mathbb{Z}$ with $t_2-t_1,t_3-t_2,t_4-t_3 \ge k$ the triple of midpoints 
$$ (\midp(p_1,p_2),\midp(p_2,p_3),\midp(p_3,p_4)) $$
is $(\alpha_0,\tau_{mod},\epsilon)$-straight and $s$-spaced. (Here $p(t_i)=p_i$.)
\end{defn}

Our next theorem says that sufficiently spaced points on  Morse quasigeodesics have straight and spaced midpoint sequences. In an effort to make Theorem \ref{midpoint sequences} readable, we have given up some control over the required spacing. For example, we use only one auxiliary parameter $\alpha_{aux}$ to control the regularity as well as the crude estimate $\sin^{-1}(x) \le \frac{\pi}{2}x$ for $0 \le x \le 1$ (this follows from the fact that $\sin^{-1}$ is convex). The following result is a quantified version of Proposition 7.16 in \cite{KLP14}.

\begin{thm}\label{midpoint sequences}
Assume $k$ is large enough in terms of  $\alpha_{new}<\alpha_0,D,\epsilon,c_1,c_2$ and $s$. To make this precise, we use auxiliary constants $l,\delta,\alpha_{aux}$ and make the following assumptions.
% due to changes in the previous section, some of these hypothoses were updated
\begin{enumerate}  
	\item\label{from k to l} Let $k$ be large enough in terms of the quasigeodesic parameters so that if $\abs{N} \ge k$ then $d(x_n,x_{n+N})\ge 2l$. Precisely, let $k \ge c_1(2l+c_2)$. Our requirements on $k$ will manifest as requirements on $l$;
	\item\label{closetodiamonds} %3.1.2021 corrected
	$$ 1 \le 6 \sinh( \alpha_{aux}(2l-2D))^2, \quad \frac{1}{\alpha_{aux}\zeta_0}\frac{D}{l} \le \frac{\zeta_0^2}{\kappa_0^2}, \quad \text{ and } \quad 5De^{2\kappa_0D-\alpha_0l}  \le \delta $$
	so that midpoints are $\delta$-close to diamonds by Lemma \ref{midpoint projection};	
	\item\label{spacing} We assume $\frac{2 \alpha_{aux}}{\kappa_0} \left( l-\delta-D \right) \ge s$ to ensure that the midpoints are appropriately spaced. 
	\item\label{auxiliaryalpha} We use an auxiliary parameter $\alpha_{aux}$ such that $\alpha_{new}<\alpha_{aux}<\alpha_0,$
	$$ \frac{\alpha_0 \delta + 3 \alpha_0 D+ 2\kappa_0D}{l-\delta-2D} \le \alpha_0-\alpha_{aux}, \quad \text{and} \quad \frac{ 2\kappa_0 \delta (\alpha_{aux}+\kappa_0)}{2 \alpha_{aux}(l-\delta-D)-2 \kappa_0\delta }\le \alpha_{aux}-\alpha_{new} .$$
	so that certain perturbations of regular segments are regular by \ref{regular projections}.
	\item\label{straightness} We assume % This didn't need to be changed.
	$$ \frac{1}{\alpha_{aux}\zeta_0} \frac{D}{l}
	+ \frac{1}{\alpha_{new}\zeta_0} \frac{\kappa_0\delta}{2\alpha_{aux}(l-\delta-D)-\delta \kappa_0} 
	+\frac{1}{2\alpha_{aux} \zeta_0} \frac{\delta}{l-D} 
	+\frac{1}{2\alpha_{new} \zeta_0} \frac{\delta}{l-\delta} + 2\kappa_0 \delta
	\le \frac{\epsilon}{\pi} $$
	to ensure that the midpoint sequence is straight. 
\end{enumerate}
Then every $(\alpha_0,\tau_{mod},D)$-Morse $(c_1,c_2)$-lower-quasigeodesic satisfies the $(\alpha_{new},\tau_{mod},\epsilon,s,k')$-quadruple condition for every $k'\ge k$. 
\end{thm}

Note that in assumption \ref{straightness}, we have in particular assumed $2 \pi \kappa_0 \delta < \epsilon$, so the $\delta$ which appears in the proof is quite small. Our proof follows \cite{KLP14} Proposition 7.16 closely.

\begin{proof}
	Let $(q_n)_{n=t_0}^{n=t_{max}}$ be an $(\alpha_0,\tau_{mod},D)$-Morse quasigeodesic and let $t_1,t_2,t_3,t_4 \in [t_0,t_{max}] \cap \mathbb{Z}$ such that $t_2-t_1,t_3-t_2,t_4-t_3 \ge k$. We abbreviate $p_i \coloneqq q_{t_i}$ and $m_i \coloneqq \midp(p_i,p_{i+1})$. We have $d(p_i,p_{i+1})\ge 2l$, $d(m_i,p_i)\ge l$ and $d(m_i,p_{i+1})\ge l$. 
	
	To show that the midpoint sequence is $(\alpha_{new},\tau_{mod},\epsilon)$-straight it suffices to show that the segment $m_2m_1$ is $(\alpha_{new},\tau_{mod})$-regular and that $\angle_{m_2}^\zeta(p_2,m_1) \le \epsilon/2$ under our assumptions on $k$.
	
	\begin{figure}[h]
		\centering
		\begin{tikzpicture}
		\node[label=below:{$P=P(\tau_-,\tau_+)$}] (f) at (14,0) {};
		\node[circle,scale=0.5,fill=black,label=above:{$p_1$}] (p1) at (2,2.5) {};
		\node[circle,scale=0.5,fill=black,label=above:{$p_2$}] (p2) at (7.5,2.5) {};
		\node[circle,scale=0.5,fill=black,label=above:{$p_3$}] (p3) at (13,2.5) {};
		\node[circle,scale=0.5,fill=black,label=left:{$x_1$}] (x1) at (2,1) {};
		\node[circle,scale=0.5,fill=black,label=above:{$x_3$}] (x3) at (13,1) {};

		\draw (0,0) -- (14,0) -- (15,2) -- (1,2) -- cycle; %draw a parallelogram to represent the parallel set
		\draw (x1) -- (7.5,1.8) -- (x3) -- (7.5,0.2) -- (x1); %draw a diamond
		\draw (p1) .. controls  (4.5,1.2) .. (p2) node[circle,fill,scale=0.5,pos=0.5,label=above:{$m_1$}] (m1) {}; 
		\draw (p2) .. controls  (10.75,1.2) .. (p3) node[circle,fill,scale=0.5,pos=0.5,label=above:{$m_2$}] (m2) {}; 

		\node[circle,scale=0.5,fill=black,label=right:{$\overline{m}_1$}] (bm1) at (m1|-52,0.9) {};
		\node[circle,scale=0.5,fill=black,label=right:{$\overline{p}_2$}] (bp2) at (p2|-52,1.4) {};
		\node[circle,scale=0.5,fill=black,label=right:{$\overline{m}_2$}] (bm2) at (m2|-52,1) {};
		\end{tikzpicture}
		\caption{The projections satisfy $\overline{p}_2 \in V(\overline{m}_1,\st(\tau_+),\alpha_{aux})$ and $\overline{m}_2 \in V(\overline{p}_2,\st(\tau_+),\alpha_{aux})$}
	\end{figure}
		
	By the Morse property there exists a diamond $\diamondsuit_{\alpha_0}(x_1,x_3)$ such that $d(x_1,p_1),d(x_3,p_3)\le D$ and $p_2$ is in the $D$-neighborhood of $\diamondsuit_{\alpha_0}(x_1,x_3)$. The diamond spans a unique parallel set $P=P(\tau_-,\tau_+)$. We denote by $\overline{p}_i$ and $\overline{m}_i$ the projections of $p_i$ and $m_i$ to $P$. 
	
	We first observe that $m_1$ is $\delta$-close to $P$ by the midpoint projection estimate Lemma \ref{midpoint projection}: we have $d(p_1,x_1)\le D,$ $d(p_2,V(x_1,\ost(\tau(x_1x_3)))) \le  d(p_2, \diamondsuit_{\alpha_0}(x_1,x_3)) \le D$ and $p_1p_2$ is $(\alpha_0,\tau_{mod})$-regular with $d(p_1,p_2)\ge 2l$ and $l$ large enough by assumption \ref{closetodiamonds} and assumption \ref{auxiliaryalpha}: %corrected
	$$ d(m_1,P) \le 5De^{2\kappa_0D - \alpha_0l} \le \delta .$$
	
	Next we look at the directions of the segments $\overline{m}_2\overline{m}_1$ and $\overline{m}_2\overline{p}_2$ and show that they have the same $\tau$-direction.
	We have 
	$$ d(\overline{p}_2,V(\overline{p}_1,\st(\tau_+),\alpha_0)) \le d(p_2,\diamondsuit_{\alpha_0}(x_1,x_3))+d(\diamondsuit_{\alpha_0}(x_1,x_3),\overline{p}_1) \le 2D $$
	since projecting to a closed convex subset is distance-non-increasing. If $c_1$ is the geodesic from $p_1$ through $p_2$, the function $t \mapsto d(c_1(t),V(\overline{p}_1,\st(\tau_+),\alpha_0))$ is convex, which implies $\overline{m}_1$ is $2D$-close to $V(\overline{p}_1,\st(\tau_+),\alpha_0)$. We have $d(\overline{m}_1,\overline{p}_1) \ge l-\delta-D$, so by using the point in $V(\overline{p}_1,\st(\tau_+),\alpha_0)$ within $2D$ of $\overline{m}_1$ and Lemma \ref{regular projections} in the presence of assumption \ref{auxiliaryalpha}, we obtain that $\overline{m}_1 \in V(\overline{p}_1,\st(\tau_+),\alpha_{aux})$. Similar arguments show that $\overline{m}_1 \in V(\overline{p}_2,\st(\tau_-),\alpha_{aux})$, or equivalently (by using the geodesic symmetry at $\midp(\overline{p}_1,\overline{p}_2)$) that $\overline{p}_2 \in V(\overline{m}_1,\st(\tau_+),\alpha_{aux})$. By the nestedness of Weyl cones, $\overline{p}_1 \in V(\overline{p}_2,\st(\tau_-),\alpha_{aux})$ and $\overline{p}_2 \in V(\overline{p}_1,\st(\tau_+),\alpha_{aux})$. Similarly, $\overline{m}_2 \in V(\overline{p}_2,\st(\tau_+),\alpha_{aux})$ and $\overline{p}_2 \in V(\overline{m}_2,\st(\tau_-),\alpha_{aux})$. The convexity of Weyl cones implies that also $\overline{m}_1 \in V(\overline{m}_2,\st(\tau_-),\alpha_{aux})$. In particular $\angle_{\overline{m}_2}^\zeta(\overline{p}_2,\overline{m}_1)=0$.
	 
	 It is convenient to show that the midpoint sequence is appropriately spaced at this point in the proof, so that we can use the resulting estimate to control the regularity parameters $\alpha_{aux}$ and $\alpha_{new}$ and the straightness parameter $\epsilon$. The inclusions $\overline{m}_1 \in V(\overline{p}_2,\st(\tau_-),\alpha_{aux})$ and $\overline{m}_2 \in V(\overline{p}_2,\st(\tau_+),\alpha_{aux})$ imply that $d(\overline{m}_1,\overline{m}_2) \ge \frac{\alpha_{aux}}{\kappa_0} \left(  d(\overline{m}_1,\overline{p}_2) + d(\overline{p}_2,\overline{m}_1) \right)$. Therefore by assumption \ref{spacing}, the midpoint sequence is appropriately spaced:
	 $$ d(m_1,m_2) \ge d(\overline{m}_1,\overline{m}_2) \ge \frac{\alpha_{aux}}{\kappa_0} \left(  d(\overline{m}_1,\overline{p}_2) + d(\overline{p}_2,\overline{m}_1) \right) \ge \frac{2 \alpha_{aux}}{\kappa_0} \left( l-\delta-D \right) \ge s. $$	 
	 
	 Using the previous estimate, Lemma \ref{regular projections}, and assumption \ref{auxiliaryalpha}, we see that $m_2m_1$ and $m_2 \overline{m}_1$ are $(\alpha_{new},\tau_{mod})$-regular and $m_2 \overline{p}_2$ is $(\alpha_{aux},\tau_{mod})$-regular.
	 	
	We may now demonstrate the bound $\angle_{m_2}^\zeta(p_2,m_1) \le \epsilon/2$. We have 
	\begin{align*}
	 \angle_{m_2}^\zeta(p_2,m_1) & = \abs{\angle_{m_2}^\zeta(p_2,m_1) - \angle_{\overline{m}_2}^\zeta(\overline{p}_2,\overline{m}_1) }  \\
	 & \le \abs{ \angle_{m_2}^\zeta(p_2,m_1) - \angle_{m_2}^\zeta(\overline{p}_2,\overline{m}_1)} \\ 
	 & + \abs{
	\angle_{m_2}^\zeta(\overline{p}_2,\overline{m}_1) -
	\angle_{m_2}^\zeta \left(\tau(\overline{m}_2\overline{p}_2),\tau(\overline{m}_2\overline{m}_1) \right)} \\ 
	& + \abs{	\angle_{m_2}^\zeta \left( \tau(\overline{m}_2\overline{p}_2),\tau(\overline{m}_2\overline{m}_1) \right) -
	\angle_{\overline{m}_2}^\zeta(\overline{p}_2,\overline{m}_1) }
	\end{align*}
	
	By the triangle inequality for quadruples (on the metric space $(\Flagt,\angle_{m_2}^\zeta)$) we have 
	$$ \abs{ \angle_{m_2}^\zeta(p_2,m_1) - \angle_{m_2}^\zeta(\overline{p}_2,\overline{m}_1)} \le \angle_{m_2}^\zeta(p_2,\overline{p}_2)+\angle_{m_2}^\zeta(m_1,\overline{m}_1) = 2 \sin^{-1}\left(\frac{1}{2} d_{\mfp}(Z_1,Z_2)\right)+2 \sin^{-1} \left(\frac{1}{2} d_{\mfp}(Z_3,Z_4) \right)$$
	where $Z_1,Z_2,Z_3,Z_4$ are the unit vectors at $m_2$ in the directions $\zeta(m_2p_2),\zeta(m_2\overline{p}_2),\zeta(m_2m_1),\zeta(m_2\overline{m}_1)$ respectively. Let $X_1,X_2,X_3,X_4$ be the unit vectors at $m_2$ which in the directions   $p_2,\overline{p}_2,m_1,\overline{m}_1$ respectively. Then by Lemma \ref{zeta projection} and the angle comparison to Euclidean space Lemma \ref{euclidean angle} we have 
	$$ d(Z_1,Z_2) \le \frac{1}{\alpha_{aux}\zeta_0} d\left( X_1,X_2 \right) =  \frac{2}{\alpha_{aux}\zeta_0} \sin \frac{1}{2} \angle_{m_2}(p_2,\overline{p}_2) \le \frac{1}{\alpha_{aux}\zeta_0} \frac{D}{l} .$$
	Similarly, 
	$$ d(Z_3,Z_4) \le \frac{1}{\alpha_{new}\zeta_0} d\left( X_3,X_4 \right) =  \frac{2}{\alpha_{new}\zeta_0} \sin \frac{1}{2} \angle_{m_2}(m_1,\overline{m}_1) \le \frac{1}{\alpha_{new}\zeta_0} \frac{\kappa_0\delta}{2\alpha_{aux}(l-\delta-D)-\delta \kappa_0} .$$
	
	Again by the triangle inequality on $(\Flagt,\angle_{m_2}^\zeta)$, 
	$$ \abs{
		\angle_{m_2}^\zeta(\overline{p}_2,\overline{m}_1) -
		\angle_{m_2}^\zeta \left(\tau(\overline{m}_2\overline{p}_2),\tau(\overline{m}_2\overline{m}_1) \right)}  \le \angle_{m_2}^\zeta(\overline{p}_2,\tau(\overline{m}_2\overline{p}_2)) +
	\angle_{m_2}^\zeta \left(\overline{m}_1,\tau(\overline{m}_2\overline{m}_1) \right) .$$
	
	Asymptotic geodesic rays are bounded by the distance of their tips, so if we let $c_2$ be the geodesic ray from $m_2$ to $\overline{m}_2\overline{p}_2(+\infty)$ we may use Lemma \ref{zeta projection} to obtain 
	$$ \sin \frac{1}{2} \angle_{m_2}^\zeta(\overline{p}_2,\tau(\overline{m}_2\overline{p}_2)) \le \frac{1}{2\alpha_{aux} \zeta_0} \frac{d(\overline{p}_2, \im c_2)}{d(m_2,\overline{p}_2)} \le \frac{1}{2\alpha_{aux} \zeta_0} \frac{\delta}{l-D}.$$
	Similarly by considering the geodesic ray $c_3$ from $\overline{m}_2$ through $\overline{m}_1$,
	$$ \sin \frac{1}{2} \angle_{m_2}^\zeta(\overline{m}_1,\tau(\overline{m}_2\overline{m}_1)) \le \frac{1}{2\alpha_{new} \zeta_0} \frac{d(\overline{m}_1, \im c_3)}{d(m_2,\overline{m}_1)} \le \frac{1}{2\alpha_{new} \zeta_0} \frac{\delta}{l-\delta}. $$	
	
	Write $\tau=\tau(\overline{m}_2\overline{p}_2)$ and $\tau'=\tau(\overline{m}_2\overline{m}_1)$. By the distance-to-angle estimate Corollary \ref{distance to angle},
	\[ \abs{	\angle_{m_2}^\zeta \left( \tau(\overline{m}_2\overline{p}_2),\tau(\overline{m}_2\overline{m}_1) \right) -
		\angle_{\overline{m}_2}^\zeta(\overline{p}_2,\overline{m}_1) } = \abs{	\angle_{m_2}^\zeta(\tau,\tau') -
		\angle_{\overline{m}_2}^\zeta(\tau,\tau')} 
	 \le 4 \sin^{-1} \left( \frac{\kappa_0}{2} d(\overline{m}_2,m_2) \right) \le 4 \sin^{-1} \left( \frac{\kappa_0\delta}{2} \right) 
	\]
	
	Combining these estimates with the fact that $\sin^{-1}(x) \le \frac{\pi}{2}x$ for $0 \le x \le 1$ yields
	$$ \angle_{m_2}^\zeta(p_2,m_1) \le 	\frac{\pi}{2} \bigg[
	\frac{1}{\alpha_{aux}\zeta_0} \frac{D}{l}
	+ \frac{1}{\alpha_{new}\zeta_0} \frac{\kappa_0\delta}{2\alpha_{aux}(l-\delta-D)-\delta \kappa_0} 
	+\frac{1}{2\alpha_{aux} \zeta_0} \frac{\delta}{l-D} 
	+\frac{1}{2\alpha_{new} \zeta_0} \frac{\delta}{l-\delta} + 2 \kappa_0 \delta  \bigg]  
	\le \frac{\epsilon}{2} $$
	by assumption \ref{straightness}. For similar reasons $\angle_{m_2}^\zeta(p_3,m_3) \le \frac{\epsilon}{2}$, so $\angle_{m_2}^\zeta(m_1,m_3) \ge \pi-\epsilon$ as desired. We have already shown that $m_2m_1$ is $(\alpha_{new},\tau_{mod})$-regular and $s$-spaced. For similar reasons the same holds for $m_2m_3$. This concludes the proof. 
\end{proof}

\begin{rem}\label{remark midpoint}
	To provide suitable auxiliary parameters to apply Theorem \ref{midpoint sequences}, we may first choose any $\delta < \frac{\epsilon}{2\pi\kappa_0}$ and any $\alpha_{new}<\alpha_{aux}<\alpha_0$. Then we may choose $l$ large enough to satisfy assumptions \ref{closetodiamonds} through \ref{straightness}, which provides a suitable $k$ via assumption \ref{from k to l}. When we apply Theorem \ref{midpoint sequences} in Section \ref{section applications} we set $\delta = \frac{\epsilon}{20 \pi \kappa_0}$ and $\alpha_{aux} = 0.8 \alpha_0 + 0.2 \alpha_{new}.$
\end{rem}

\subsection{Local-to-global principle for Morse quasigeodesics}

An $L$\textit{-local} $(\alpha_0,\tau_{mod},D)$\textit{-Morse} $(c_1,c_2,c_3,c_4)$\textit{-quasigeodesic} is a sequence $(x_n)_{n=t_0}^{n=t_{max}}$ in $\X$ such that for $t_0 \le t_1 \le t_2 \le t_{max}$ with $t_2-t_1 \le L$, the subsequence $(x_n)_{n=t_1}^{n=t_2}$ is an $(\alpha_0,\tau_{mod},D)$-Morse $(c_1,c_2,c_3,c_4)$-quasigeodesic.

We now come to the main result of the paper. The following result is a quantified local-to-global principle for Morse quasigeodesics. Theorem \ref{local to global} says that for any fixed quality of Morse quasigeodesic, there exists a large enough scale so that a local Morse quasigeodesic of that scale and quality is a global Morse quasigeodesic. It is a quantified version of Theorem 7.18 in \cite{KLP14}, stated as Theorem \ref{main theorem} in the introduction. We will apply Theorem \ref{straight and spaced sequences are morse} and Theorem \ref{midpoint sequences}. While these theorems have cumbersome statements, finding auxiliary parameters which satisfy the required inequalities is easy, as we discussed in Remark \ref{remark straight to morse} and Remark \ref{remark midpoint}, and as we demonstrate in the next section. 

\begin{thm}\label{local to global}
	For any $\alpha_{new}<\alpha_0,D,c_1,c_2,c_3,c_4$, there exists a scale $L$ so that every $L$-local $(\alpha_0,\tau_{mod},D)$-Morse $(c_1,c_2,c_3,c_4)$-quasigeodesic in $\X$ is an $(\alpha_{new},\tau_{mod},D')$-Morse $(c_1',c_2',c_3',c_4')$-quasigeodesic. Precisely, $L=3k$ is large enough if auxiliary parameters $\alpha_{aux},k,\delta,s,\epsilon$ satisfy:
	\begin{enumerate}
		\item\label{apply theorem 4.1} $\epsilon$ is small enough and $s$ is large enough to satisfy the conditions of Theorem \ref{straight and spaced sequences are morse} for $\alpha_{new}<\alpha_{aux},\delta$,
		\item\label{apply theorem 4.4} $k$ is large enough in terms of $\alpha_{aux}<\alpha_0,D,\epsilon,c_1,c_2$ and $s$ to satisfy the conditions of Theorem \ref{midpoint sequences},
	\end{enumerate}
	and the sequence has global Morse parameters
	\begin{enumerate}
		\item $D'=c_3k+\frac{3}{2}c_4+\delta,$
		\item $(c_1')^{-1} = 2 \alpha_{new} \zeta_0 c_0(s-2\delta)k^{-1}$,
		\item $c_2'= 2 \alpha_{new} \zeta_0 c_0 (s-2\delta) + 2\delta + 2 c_3 k + 3 c_4,$ 
		\item $c_3'= c_3+\frac{c_4}{L},$
		\item $c_4'=c_4.$
	\end{enumerate}	
\end{thm}

\begin{proof}
	Let $(x_n)_{n=-\infty}^{n=+\infty}$ be an $L$-local $(\alpha_0,\tau_{mod},D)$-Morse $(c_1,c_2,c_3,c_4)$-quasigeodesic. By Theorem \ref{midpoint sequences} and assumption \ref{apply theorem 4.4}, each subsequence $(x_n)_{n=t_0}^{n=t_0+3k}$ satisfies the $(\alpha_{aux},\tau_{mod},\epsilon,s,k)$-quadruple condition. In particular, the coarse midpoint sequence $m_n=\midp(x_{nk},x_{nk+k})$ is $(\alpha_0,\tau_{mod},\epsilon)$-straight and $s$-spaced. By Theorem \ref{straight and spaced sequences are morse} and assumption \ref{apply theorem 4.1}, the midpoint sequence $(m_n)$ is an $(\alpha_{new},\tau_{mod},\delta)$-Morse $((2\alpha_{new}\zeta_0c_0(s-2\delta))^{-1},2\delta)$-lower quasigeodesic. We now use the midpoint sequence as a coarse approximation of the original sequence to show that $(x_n)$ is a global Morse quasigeodesic.
	
	The subsequences $x_{nk},x_{nk+1},\dots,x_{nk+k-1},x_{nk+k}$ are $(c_3,c_4)$-upper-quasigeodesics (because $L \ge k$), so they lie in uniform neighborhoods of each $m_n$: if $\abs{t-nk} \le \frac{k}{2}$ then
	$$ d(m_n,x_t) \le d(m_n,x_{nk})+d(x_{nk},x_t) \le \frac{d(x_{nk},x_{nk+k})}{2}+d(x_{nk},x_t) \le \frac{c_3}{2}k+\frac{c_4}{2}+c_3 \frac{k}{2}+c_4 = c_3k+\frac{3}{2}c_4 .$$
	In particular, $(x_n)$ is $(\alpha_{new},\tau_{mod},D')$-Morse for $D'=c_3k+\frac{3}{2}c_4+\delta.$ The midpoint sequence is coarsely spaced:
	$$ d(m_n,m_{n+N})\ge 2 \alpha_{new} \zeta_0 c_0(s-2\delta)\abs{N}-2\delta ,$$
	so the original sequence is also coarsely spaced:
	\begin{align*}
		d(x_t,x_{t'}) & \ge d(m_n,m_{n'}) - d(m_n,x_t)-d(m_{n'},x_{t'}) \\
		& \ge 2 \alpha_{new} \zeta_0 c_0 (s-2\delta)\abs{n-n'} -2\delta - 2 c_3 k - 3 c_4 \\
		& \ge 2 \alpha_{new} \zeta_0 c_0 (s-2\delta) k^{-1} \abs{t-t'} - 2 \alpha_{new} \zeta_0 c_0 (s-2\delta) - 2\delta - 2 c_3 k - 3 c_4 .
	\end{align*}
	Finally, if a sequence is $(c_3,c_4)$-coarsely Lipschitz on intervals of length $L$, it then satisfies $d(x_n,x_{n+N})\le \abs{N}(c_3+\frac{c_4}{L})+c_4$ and is $(c_3+\frac{c_4}{L},c_4)$-coarsely Lipschitz.
\end{proof}

\section{Applications of the local-to-global principle}\label{section applications}

In this section we give two applications of the main result, Theorem \ref{local to global}. We describe two explicit neighborhoods of Anosov representations in $\SL(3,\R)$, one for free groups and another for closed surface groups. Each of them is constructed by perturbing a group acting cocompactly on a convex subset of a totally geodesic hyperbolic plane in the associated symmetric space. 

We will need some further estimates in order to quantify these neighborhoods. First we recall a standard proof of the Milnor-Schwarz Lemma so that we may use the explicit quasi-isometry constants it produces. We then give a version of the classical Morse Lemma that will be used in Section \ref{section surface group}. In Section \ref{section matrix estimates} we use elementary linear algebra to control the perturbations of long words in a linear group that results from perturbing the generators. We also relate the Frobenius norm on $d \times d$ matrices to the distance in the symmetric space associated to $\SL(d,\R)$. In the final two sections, we apply the local-to-global principle Theorem \ref{local to global} to describe explicit neighborhoods of Anosov representations.

As one might expect, straightforward applications of Theorem \ref{local to global} as we have done here will yield only very small perturbations. This is partially explained by the following geometric difficulty. The Morse condition implies that the image of each geodesic in the Cayley graph fellow-travels a unique parallel set. After perturbing the representation, one expects the image of the geodesic to fellow-travel a new parallel set. For geodesics through the identity, our techniques merely bound the distance from the perturbed geodesic to its previous parallel set, so for it to fellow-travel for a long time, the perturbation has to be extremely small. If we could identify the new parallel set it fellow-travels and bound the distance to that parallel set, we expect that the perturbation bounds would improve significantly. 

\subsection{Preliminary estimates}\label{section preliminary estimates}

\subsubsection{The Milnor-Schwarz Lemma}\label{section milnor-schwarz}

In this subsection we state and prove a standard result in geometric group theory called the Milnor-Schwarz Lemma. It is a source of concrete quasiisometry parameters for nice enough actions of finitely generated groups, such as those we consider in Sections \ref{section free group} and \ref{section surface group}. The proof given here is taken directly from Sisto's lecture notes \cite{Sis14}.

%I considered altering Sisto's proof, but I couldn't improve the constants is a seemingly worthwhile way. Although I did observe that we can replace the $+2$ with $+1$. this can be further slightly improved by replacing +1 with +[1-2R], so we get (1,0) (optimal) when R is at least 1/2 (often). 

\begin{lem}[Milnor-\v{S}varc Lemma]\label{milnor-schwarz lemma}
	Let $G$ be a group acting properly discontinuously, cocompactly and by isometries on a proper geodesic space $X$. Choose any $p \in X$. Then the group $G$ has a finite generating set $S$ so that the orbit map at $p$ is a quasi-isometry for $G$ with the word metric induced by $S$. In fact,
	$$ wl(g) \le d(p,gp)+1, \quad \text{and} \quad d(p,gp) \le \max_{s \in S} \{d(p,sp)\} wl(g) .$$
\end{lem}

\begin{proof}
	Since the action is cocompact, there exists a constant $R$ so that the $G$-translates of $B_{R}(p)$ cover $X$. Let $S:=\{g \in G \mid d(p,gp) \le 2R+1 \}$. Since $X$ is proper, the closed ball of radius $R+\frac{1}{2}$ centered at $p$ is compact, and since the action is properly discontinuous, $S=\{g \in G \mid B_{R+\frac{1}{2}}(p) \cap B_{R+\frac{1}{2}}(gp) \}$ is finite. Now let $g \in G$. Choose a minimal geodesic from $p$ to $gp$, and subdivide it with points $p_i$ so that $p=p_0,p_1,p_2,\dots,p_{n-1},p_n=gp$ occur monotonically and for $i=0,1,2,\dots,n-2$, we have $d(p_i,p_{i+1})=1$ and $d(p_{n-1},p_n)\le 1$. For each $1 \le i \le n-1$ choose $g_i \in G$ so that $d(g_ip,p_i) \le R$ and set $g_0= \id$ and $g_n=g$. Then for all $0 \le i\le n-1$, we have 
	$$ d(g_ip,g_{i+1}p)\le d(g_ip,p_i)+d(p_i,p_{i+1})+d(p_{i+1},g_{i+1}p) \le 2R+1,$$
	which implies that there exists $s_{i+1} \in S$ so that $g_{i+1}=g_is_{i+1}$. For all $1 \le i \le n$ it follows that $g_i=s_1s_2s_3\cdots s_i$. Therefore $g$ can be written as a product of $n$ elements of $S$, with $n-1 \le d(p,gp)$. It follows that $S$ is a finite generating set for $G$ and the word length of $g$ with respect to $S$ is bounded above by $d(p,gp)+1$.
	
	We have shown that $S$ is a finite generating set for $G$. Write $g=g_1\cdots g_n$ with $g_i \in S$. Then 
	\begin{align*}
		d(p,g_1 g_2 g_3 \cdots g_n p) 
		& \le d(p, g_1 \cdots g_{n-1} p) + d(g_1 \cdots g_{n-1} p, g_1 \cdots g_{n-1} g_n p) \\
		& = d(p, g_1 \cdots g_{n-1} p) + d(p, g_n p) \\
		& \le d(p,g_1p) + \cdots + d(p,g_n p) \\
		& \le \max_{s \in S} \{d(p,sp)\} n,
	\end{align*}
	so the orbit map at $p$ is $\max_{s \in S} \{d(p,sp)\}$-Lipschitz with respect to the generating set $S$. Note that by the definition of $S$, $\max_{s \in S} \{d(p,sp)\} \le 2R+1$.
\end{proof}

The previous lemma provides quasi-isometry constants in terms of only the constant $R$ so that the image of an $R$-ball covers the quotient. In return we give up control over the generating set. In particular, when we apply Lemma \ref{milnor-schwarz lemma} to an action of a closed surface group on the hyperbolic plane in Section \ref{section surface group}, we will give quasiisometry parameters with a nonstandard generating set for the Cayley graph. We will need to control the Frobenious norm of the matrices in our generating set by using Lemma \ref{control frobenius generators}.

 %The right reason is: computing the Morse constant by hand for the standard generating set is a pain, and the Morse Lemma spits out a huge one.

\subsubsection{The classical Morse Lemma}\label{section classical morse lemma}

In Section \ref{section surface group} we will use the following version of the classical Morse Lemma to provide Morse quasiisometry parameters for the orbit map of a surface group acting on a copy of the hyperbolic plane. The following proof is adapted from Bridson-Haefliger \cite{BH99}.

\begin{thm}[Classical Morse Lemma, Cf.\ {\cite{BH99}{Theorem III.H.1.7}}]\label{classical morse lemma}
	Let $D_0$ be an upper bound for 
	$$\{ D \mid D-1 \le \delta \abs{\log_2 ( 2D+2M^2l+6DMl+aM ) } \} $$
	and set $R=D_0+lMD_0+lM^2+\frac{a}{2}$. Then: \\
	
	If $(y_i)_{i=0}^{i=N}$ is a sequence in a $\delta$-hyperbolic geodesic space $\Y$ with 
	$$ d(y_i,y_j) \le M|j-i| \text{ and } |j-i|\le ld(y_i,y_j)+a $$
	then for all $0 \le n \le N$, the distance from $y_n$ to a geodesic segment from $y_0$ to $y_N$ is bounded above by $R$. \\
\end{thm}

\begin{proof} 
	Let $c \colon [0,N]\to \Y$ be the piecewise geodesic curve with $c(i)=y_i$. Let $D$ be minimal so that the closed $D$-neighborhood of $\im c$ covers the geodesic from $p=y_0$ to $q=y_N$. Choose a point $x_0$ on $pq$ realizing $D$, and choose $y,z$ on $pq$ at distance $2D$ from $x_0$ so that $y,x_0,z$ occurs in order (if $x_0$ is too close to $p$, use $p$ for $y$, and likewise for $z$). Choose $y'$ on $\im c$ within $D$ of $y$, and choose $z'$ similarly. Choose $i,j$ so that $y'$ is on $y_iy_{i+1}$ and $z'$ is on $y_{j-1}y_j$. If $c(t)=y'$ and $c(t')=z'$ then the length of $c$ restricted to the $[t,t']$ is at most 
	$$ \length (c|_{[t,t']}) \le \length (c|_{[i,j]}) \le M|j-i| \le M[ld(y_i,y_j)+a] .$$
	Also, 
	$$ d(y_i,y_j) \le d(y_i,y') + d(y',y) + d(y,z) + d(z,z') + d(z',y_j) \le 2M+6D $$
	and it follows that the curve $c'$ formed by following a geodesic segment from $y$ to $y'$ then along $c$ to $z'$ then along a geodesic segment to $z$ has length at most $2D+M[l (2M+6D)+a]$. Proposition III.H.1.6 in \cite{BH99} bounds $D$ in terms of the length of $c'$ and $\delta$. In particular  
	$$ D-1 \le \delta |\log_2 ( 2D+2M^2l+6DMl+aM ) | $$
	which implies an upper bound $D_0$ on $D$.
	
	Now suppose that $(y_n)_{n=a'}^{n=b'}$ is a maximal (consecutive) subsequence outside the $D_0$-neighborhood of $pq$. There exist $s,s'$ such that $0 \le s \le a'$ and $b' \le s' \le N$ within $D_0$ of the same point on $pq$, so $d(c(s),c(s')) \le 2D_0$. As before, by choosing $m,n$ so that $c(s)$ lies on $y_my_{m+1}$ and $c(s')$ lies on $y_ny_{n+1}$ we have that 
	$$ \length( c|_{[s,s']}) \le \length(c|_{[m,n]}) \le M|m-n| \le M(l d(y_m,y_n)+a) $$ 
	and 
	$$ d(y_m,y_n) \le d(y_m,c(s)) + d(c(s),c(s'))+d(c(s'),y_n) \le 2M+2D_0 $$
	so we obtain
	$$ \length c|_{[s,s']} \le M[l(2D_0+2M)+a].$$ It follows that $R= D_0 + M[l(D_0+M)+\frac{a}{2}]$ is an upper bound for the distance from any $y_n$ to $pq$. 
\end{proof}

% order 0 s i j s' N
% length c|[s,s'] \le length c|[i,j] + 2M \le M|j-i| +2M \le M( ld(y_i,y_j)+a +2) \le M( 2lD_0+a +2 +2l) 

\subsubsection{Matrix Estimates} \label{section matrix estimates}

In this subsection we establish a few elementary estimates related to the symmetric space associated to $\SL(d,\R)$. We will control perturbations of long words in a generating set in terms of the Frobenious norm of the generators. As noted above, we use a non-standard generating set for the closed surface group, so we also prepare to control the Frobenious norm of the generators in that case. In Sections \ref{section free group} and \ref{section surface group}, we combine these estimates with the local-to-global principle Theorem \ref{local to global} to guarantee that the Morse subgroups under consideration remain Morse after certain explicit perturbations.

In the rest of the paper, we identify the symmetric space associated to $\SL(d,\R)$ with the space of real, symmetric, positive-definite matrices of determinant $1$. We remind the reader that we take the Riemannian metric to be induced by the Killing form, so at the identity matrix, the Riemannian metric is $2d$ times the Frobenious inner product $\langle X,Y \rangle_{Fr} = \trace (X^T Y)$. 

%This seemingly circuitous estimate is apparently better than the more direct option.

\begin{lem}\label{perturbation of long words}
	Let $\abs{\cdot}$ be any submultiplicative norm on $d \times d$ matrices. Let $w = g_1 g_2 \cdots g_{k-1}g_k$ be a product of $k$ matrices, and let $w'=(g_1+\epsilon_1)(g_2+\epsilon_2)\cdots(g_{k-1}+\epsilon_{k-1})(g_k+\epsilon_k)$ be a product of perturbed matrices. Suppose that for all $1\le i \le k$, $\abs{g_i}\le A$ and $\abs{\epsilon_i} \le \epsilon $. 
	If $k \ge 3$ and $\frac{k-1}{2} \frac{\epsilon}{A} \le 1$ then $\abs{w'-w} \le 2kA^{k-1} \epsilon .$
\end{lem}

\begin{proof}
	We have
	\begin{align*}
		\abs{w' - w} & = \abs{ \prod_{i=1}^k (g_i+\epsilon_i) - \prod_{i=1}^k g_i } \\
		& = \abs{ \sum_{ \substack{ 1 \le i_1 \le \cdots \le i_j \le k \\ j=1}}^{j=k} g_1 g_2 \cdots g_{i_1-1} \epsilon_{i_1} g_{i_1+1} \cdots g_{i_j-1} \epsilon_{i_j} g_{i_j+1} \cdots g_k } \\
		& \le k A^{k-1}\epsilon + \binom{k}{2} A^{k-2}\epsilon^2 + \dots + \binom{k}{j} A^{k-j}\epsilon^j + \dots + \epsilon^k \\
		& = A^k \left[ \left(1+\frac{\epsilon}{A} \right)^k-1 \right] \\
		& \le 2k A^{k-1} \epsilon 
	\end{align*}	
	where the last line follows from the Taylor approximation $(1+\frac{\epsilon}{A})^k-1 \le k \frac{\epsilon}{A} + \frac{k(k-1)}{2}\left( \frac{\epsilon}{A} \right)^2 $, valid when $\frac{\epsilon}{A} \le 1$.
\end{proof}

We next relate the Riemannian distance in $\X$ to the $B_p$-norm on the space of matrices. Recall that when $p$ is the identity matrix, $B_p$ is $2d$ times the Frobenius inner product. We let $B_p$ be defined on all of $\mathfrak{gl}(d,\R)$ as $2d$ times the Frobenius inner product.

\begin{lem}\label{from frob to symm space}
	Let $g \in \SL(d,\R)$ and $p \in \X$ be the identity matrix. Then 
	$$ d_{{\X}}(gp,p) \le \sqrt{d}(d-1) \pnorm{g-1} .$$
\end{lem}

\begin{proof}
	$K=\SO(d)$ acts on $(\mathfrak{gl}(d,\R), B_p)$ by isometries on the left and the right, so $\pnorm{g-1}=\pnorm{ e^A -1 }$ where 
	$$ A = \begin{pmatrix} \lambda_{1} & & \\ & \ddots & \\ & & \lambda_d \end{pmatrix} $$
	is the Cartan projection of $g$. That is, $A$ is the unique diagonal matrix with $\lambda_1 \ge \lambda_2 \ge \cdots \ge \lambda_d$ and $\lambda_1 + \cdots + \lambda_d=0$ such that $g=ke^Ak'$ for some $k,k' \in \SO(d)$. We have $\pnorm{A}=d(gp,p)$. Since $-\lambda_d \le (d-1)\lambda_1$ and $\lambda_1^2 \le (e^{\lambda_1}-1)^2$, 
	$$ d(gp,p)^2 = \pnorm{A}^2 = 2d\sum_{i=1}^d \lambda_i^2 \le 2d^2(d-1)^2 \lambda_1^2 \le 2d^2(d-1)^2 \sum_{i=1}^d (e^{\lambda_i} - 1)^2 = d(d-1)^2 \pnorm{ e^A - 1}^2 = d(d-1)^2 \pnorm{ g-1}^2 .$$
\end{proof}

In the following corollary, we consider a pair of linear representations that map the generating set to nearby generators. We apply a long word to the basepoint using each representation, and bound the resulting distance.

\begin{cor}\label{perturbed representation}
	Let $\Gamma$ be a group with symmetric generating set $S=\{ \gamma_1,\dots,\gamma_n \}$ and let $\rho$ and $\rho'$ be two representations of $\Gamma$ into $\SL(d,\R)$. Assume that
	\begin{enumerate}
		\item For $i \in \{1,\dots,n \},$ $\abs{\rho(\gamma_i)}_{Fr}\le A$ and $\abs{\rho(\gamma_i)-\rho'(\gamma_i)}_{Fr} \le \epsilon$; and
		\item $k \ge 3$ and $ \frac{k-1}{2} \frac{\epsilon}{A} \le 1$.
	\end{enumerate}
	Then for any $\gamma \in \Gamma$ with $d_S(\gamma,1)\le k$, it holds that $d_{{\X}}(\rho'(\gamma)p,\rho(\gamma)p) \le \sqrt{8}d(d-1)kA^{2k-1}\epsilon$.
\end{cor}

\begin{proof}
	Let $g=\rho(\gamma)$ and $g'=\rho'(\gamma)$ for $d_S(\gamma,1)\le k$. Since the Frobenius norm is submultiplicative we have $\abs{g^{-1}}_{Fr} \le A^k$ and moreover because of the assumptions, Lemma \ref{perturbation of long words} applies and we obtain $\abs{g-g'}_{Fr} \le 2kA^{k-1}\epsilon$. We see that 
	$$ \abs{ g^{-1}g'-1}_{Fr} = \abs{ g^{-1}(g'- g)}_{Fr} \le \abs{ g^{-1} }_{Fr} \abs{ g'-g}_{Fr} \le A^k \abs{g' - g}_{Fr} \le 2k A^{2k-1}\epsilon.$$
	Then by applying Lemma \ref{from frob to symm space} to $g^{-1}g'$ we obtain
	$$ d(g'p,gp)= d(g^{-1}g'p,p) \le \sqrt{d}(d-1) \pnorm{g^{-1}g'-1} \le \sqrt{8}d(d-1)kA^{2k-1}\epsilon .$$	
\end{proof}

In the next lemma we give a precise, quantitative version of the following statement: If a representation $\rho$ induces a Morse quasiisometric embedding, then its perturbation $\rho'$ induces a local Morse quasiisometric embedding.

\begin{lem}
	Let $\rho,\rho' \colon \Gamma \to \SL(d,\R)$ be representations and let $S$ be a symmetric generating set for $\Gamma$. If $d( \rho(\gamma)p,\rho'(\gamma)p) \le \epsilon$ for all $d_S(\gamma,1) \le k$ and if the orbit map of $\rho$ at $p$ is an $(\alpha_0,\tau_{mod},D)$-Morse $(c_1,c_2,c_3,c_4)$-quasiisometric embedding then the orbit map of $\rho'$ at $p$ is a $2k$-local $(\alpha_0,\tau_{mod},D+\epsilon)$-Morse $(c_1,c_2+\epsilon,c_3,c_4+\epsilon)$-quasiisometric embedding. 
\end{lem}

\begin{proof}
	If $d( \rho(\gamma)p,\rho'(\gamma)p) \le \epsilon$ for all $d_S(\gamma,1) \le k$, then for every geodesic $(\gamma_n)_{n=-k}^{n=k}$ in $\Gamma$ of length $2k$,
	$$ d(\rho'(\gamma_n)p,\rho'(\gamma_0)p) = d( \rho'(\gamma_0^{-1})\rho'(\gamma_n)p,p) $$
	is within $\epsilon$ of $d(\rho(\gamma_n)p,\rho(\gamma_0)p)$. Additionally, if $(\rho(\gamma_n)p)$ is within $D$ of $\diamondsuit(q,r)$, then $(\rho'(\gamma_n)p)$ is within $D+\epsilon$ of $\diamondsuit( \rho'(\gamma_0)\rho(\gamma_0^{-1})q,\rho'(\gamma_0)\rho(\gamma_0^{-1})r)$. In particular, if $\rho$ induces an $(\alpha_0,\tau_{mod},D)$-Morse $(c_1,c_2,c_3,c_4)$-quasiisometric embedding then $\rho'$ induces a $2k$-local $(\alpha_0,\tau_{mod},D+\epsilon)$-Morse $(c_1,c_2+\epsilon,c_3,c_4+\epsilon)$-quasiisometric embedding. 
\end{proof}

When we apply the Milnor-Schwarz lemma we use the generating set $S=\{ s \in \Gamma \mid d(p,sp) \le 2R+1 \}$, and when we apply Corollary \ref{perturbed representation} we need to bound the size of the generating set. The following Lemma helps us do just that. 

\begin{lem}\label{control frobenius generators}
	Let $p$ be the identity matrix in $\X_d$ and let $g \in \SL(d,\R)$ such that $d(p,gp) \le 2R+1$. Let $\abs{\cdot}_{Fr}$ denote the Frobenius norm. Then 
	$$ \abs{g}_{Fr} \le \exp \left( \frac{2R+1}{\sqrt{2d}} \right) .$$
\end{lem}

\begin{proof}
	Combine
	$$ \abs{g}^2_{Fr} = \abs{gg^T}_{Fr} = \abs{ \exp \log g g^T}_{Fr} \le \exp \abs{\log gg^T}_{Fr} $$
	and
	$$ \sqrt{\frac{d}{2}} \abs{ \log g g^T}_{Fr} = \frac{1}{2} \pnorm{ \log g g^T} = \pnorm{ \log \sqrt{ g g^T} } = d(p,gp) \le 2R+1 $$
	to obtain
	$$ \abs{g}_{Fr} \le \exp \frac{1}{2}\abs{\log g g^T}_{Fr} \le \exp \left( \frac{2R+1}{\sqrt{2d}} \right) .$$
\end{proof}

\subsection{An explicit neighborhood of Anosov free groups}\label{section free group}

In this subsection we obtain an explicit non-empty neighborhood of Anosov free groups. Let $\Gamma_1$ be the subgroup of $\SL(3,\R)$ generated by 
$$ g= \begin{bmatrix} e^t & 0 & 0 \\ 0 & 1 & 0 \\ 0& 0 & e^{-t} \end{bmatrix} , \quad  h= \begin{bmatrix} \cosh t & 0 & \sinh t \\ 0 & 1 & 0 \\ \sinh t & 0 &  \cosh t \end{bmatrix} .$$

As in Section \ref{section matrix estimates} we identify the associated symmetric space with the space of real, symmetric, positive-definite matrices of determinant $1$. Let $p \in \X$ be the identity matrix. $\Gamma_1$ is a subgroup of a copy of $\SL(2,\R)$ preserving a copy of $\HH^2$ containing $p$ of curvature $-\frac{1}{3}$, see Section \ref{curvature}. We will directly estimate the Morse quasiisometry parameters of the orbit map at $p$ on $\Gamma_1$. 

The points $p,gp,hp$ form an isosceles right triangle:
$$ d(p,gp) = \pnorm{\begin{bmatrix} t & 0 & 0 \\ 0 & 0 & 0 \\ 0& 0 & -t \end{bmatrix}} = 2\sqrt{3}t = \pnorm{ \begin{bmatrix} 0 & 0 & t \\ 0 & 0 & 0 \\ t & 0 &  0 \end{bmatrix} } = d(p,hp) .$$

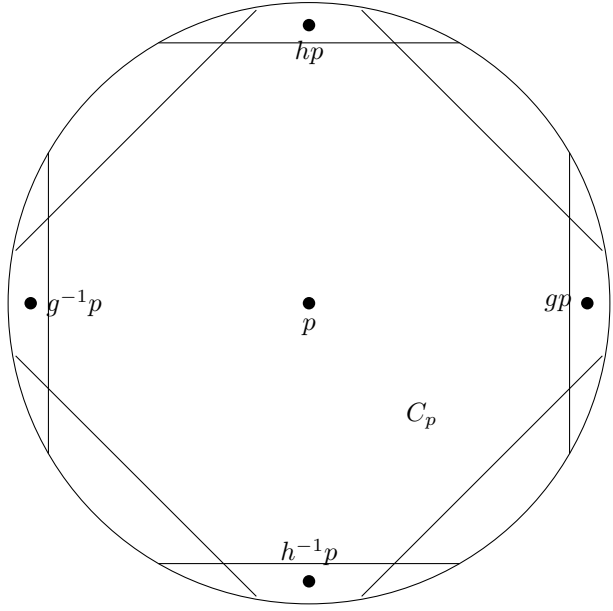
\begin{figure}[h]
	\centering
	\begin{tikzpicture}
		\node[circle,scale=0.5,fill=black,label=below:$p$] (p) at (0,0) {};
		\draw (p) circle (4);
		\coordinate  (x1) at (2,{2*sqrt(3)});
		\coordinate  (x2) at (2,{-2*sqrt(3)});
		\coordinate  (x3) at (-2,{2*sqrt(3)});
		\coordinate  (x4) at (-2,{-2*sqrt(3)});
		\coordinate  (x5) at ({2*sqrt(3)},2);
		\coordinate  (x6) at ({2*sqrt(3)},-2);
		\coordinate  (x7) at ({-2*sqrt(3)},2);
		\coordinate  (x8) at ({-2*sqrt(3)},-2);
		\coordinate  (y1) at (0.7,3.9);
		\coordinate  (y2) at (0.7,-3.9);
		\coordinate  (y3) at (-0.7,3.9);
		\coordinate  (y4) at (-0.7,-3.9);
		\coordinate  (y5) at (3.9,0.7);
		\coordinate  (y6) at (3.9,-0.7);
		\coordinate  (y7) at (-3.9,0.7);
		\coordinate  (y8) at (-3.9,-0.7);
		\draw (x1.center) -- (x3.center);
		\draw (x2.center) -- (x4.center);
		\draw (x5.center) -- (x6.center);
		\draw (x7.center) -- (x8.center);
		\draw (y1.center) -- (y5.center);
		\draw (y2.center) -- (y6.center);
		\draw (y3.center) -- (y7.center);
		\draw (y4.center) -- (y8.center);
		\node[circle,scale=0.5,fill=black,label=left:$gp$] (gp) at (3.7,0) {};
		\node[circle,scale=0.5,fill=black,label=right:$g^{-1}p$] (gip) at (-3.7,0) {};
		\node[circle,scale=0.5,fill=black,label=below:$hp$] (hp) at (0,3.7) {};
		\node[circle,scale=0.5,fill=black,label=above:$h^{-1}p$] (hip) at (0,-3.7) {};
		\node (Cp) at (1.5,-1.5) {$C_p$};
	\end{tikzpicture}
	\caption{The Dirichlet domain $C_p$ in the projective model for $\HH^2$}
\end{figure}

Write $T=\tanh(t)$. If $\sqrt{2}T >1$, then $\Gamma_1$ acts cocompactly on a closed convex subset $C$ of $\HH^2$, with a Dirichlet domain $C_p$. The domain $C_p$ is an octagon with geodesic boundary and neighbors $gC_p,g^{-1}C_p,hC_p,h^{-1}C_p$ in $C$. Since $C$ is convex, the minimum distance between any pair of neighbors is bounded below by the length of an arc in $C_p$ joining non-adjacent edges. This has lower bound
$$ c_1^{-1} = \sqrt{3} \min  \left\{  t, \frac{1}{2} \log \left(\frac{T^2+\sqrt{2 T^2-1}}{T^2-\sqrt{2 T^2-1}}\right) , \frac{1}{2} \log \left(\frac{1+ 2T \sqrt{1-T^2}}{1-2 T \sqrt{1-T^2}}\right) \right\} .$$
We also set $c_3 = 2\sqrt{3}t$. The orbit map is a $(c_1,0,c_3,0)$ quasi-isometry. Set $R = \sqrt{3}\tanh^{-1} \left(\sqrt{T^{-2}-2+2T^2} \right)$. Then $C$ is within the $R$-neighborhood of $\Gamma_1 \cdot p$ and the diameter of $C_p$ is $2R$. The orbit map is $R$-Morse. \\

We are now in position to prove Theorem \ref{main theorem demonstration 1}.

\begin{manualtheorem}{1.2}
	Let $\Gamma_1$ be the subgroup of $\SL(3,\R)$ generated by 
	$$ g= \begin{bmatrix} e^t & 0 & 0 \\ 0 & 1 & 0 \\ 0& 0 & e^{-t} \end{bmatrix} , \quad  h= \begin{bmatrix} \cosh t & 0 & \sinh t \\ 0 & 1 & 0 \\ \sinh t & 0 &  \cosh t \end{bmatrix} ,$$
	with $\tanh t=0.75$. If $\Gamma_1'$ is generated by $g',h'$ where $\max \{ \abs{g-g'}_{Fr},\abs{h-h'}_{Fr} \} \le  10^{-15,309} $, then $\Gamma_1'$ is Anosov.
\end{manualtheorem}

Before proceeding to the proof, we discuss how to choose suitable parameters in the application of Theorem \ref{local to global}. There are a number of auxiliary parameters appearing in Theorems \ref{straight and spaced sequences are morse} and \ref{midpoint sequences}. We will choose these auxiliary parameters in the same way in Section \ref{section surface group}. Because of the large number of auxiliary parameters, it is not clear how to obtain optimal estimates, even when treating Theorems \ref{straight and spaced sequences are morse}, \ref{midpoint sequences} and \ref{local to global} as black boxes. The choices we make here are simply the result of selecting auxiliary parameters in a few different ways and choosing the best result (smallest $k$) we achieved. We used a Mathematica notebook to verify the system of inequalities for each theorem. 

%corrected.
First we choose auxiliary parameters $\delta=\frac{\zeta_0}{2\kappa_0^2}$ and $\alpha_{aux}:= 0.5 \alpha_0 + 0.5 \alpha_{new}$. We apply Theorem \ref{straight and spaced sequences are morse} with $\alpha_{aux} < \alpha_0$ and $\delta = \frac{\zeta_0}{2\kappa_0^2}$ by setting $\epsilon = \frac{\zeta_0^2}{10\kappa_0^2}$ and then choosing $s$ large enough to satisfy the assumptions of the theorem. In Theorem \ref{midpoint sequences}, for any choice of auxiliary parameters $\delta_{aux} < \frac{\epsilon}{2\pi \kappa_0}$ and any $\alpha_{aux} < \alpha_{aux}' < \alpha_0$, there is a large enough auxiliary parameter $l$ to satisfy the assumptions. We select $\delta_{aux}:= 0.1 \frac{\epsilon}{2\pi\kappa_0}=0.1 \frac{\zeta_0^2}{20\pi\kappa_0^3}$ and $\alpha_{aux}':= 0.8 \alpha_0 + 0.2 \alpha_{aux}$.

\begin{proof}[Proof of Theorem 1.2]
As discussed earlier in this section, the orbit map of $\Gamma_1$ is a $(\zeta_0,\sigma_{mod},3.18)$-Morse $((1.28)^{-1},0,3.38,0)$-quasigeodesic embedding. We relax the parameters, asking the perturbation to induce a $33,602$-local $(\zeta_0,\sigma_{mod}, 3.28)$-Morse $(1,0.1,3.38,0.1)$-quasiisometric embedding. By Theorem \ref{local to global}, such an orbit map is a global $(0.95 \zeta_0; \sigma_{mod};37,858)$-Morse $(91;75,838;3.38;0)$-quasiisometric embedding. 
%corrected 3.8.2021

If $g',h' \in \SL(3,\R)$ satisfy $\abs{g-g'}_{Fr},\abs{h-h'}_{Fr}\le 10^{-15,309}$, then for $d_{\Gamma_1}(w,1) \le k=16,801$ we have $d(\rho(w)p,\rho'(w)p)\le 0.1$ by Corollary \ref{perturbed representation}, so $\rho'$ also induces a $33,602$-local $(\zeta_0,\sigma_{mod}, 3.28)$-Morse $(1,0.1,3.38,0.1)$-quasiisometric embedding and therefore its orbit map at $p$ is a (global) Morse quasiisometric embedding. In particular, $g',h'$ generate an Anosov subgroup of $\SL(3,\R)$ and our proof of Theorem \ref{main theorem demonstration 1} is complete. 
\end{proof}

\subsection{An explicit neighborhood of Anosov surface groups} \label{section surface group}

Let $\Gamma_2$ be the subgroup of $\SL(3,\R)$ generated by 
$$ S = \left\{ \begin{bmatrix} \cos \theta & 0 & \sin \theta \\ 0 & 1 & 0 \\ -\sin \theta & 0 & \cos \theta \end{bmatrix}  \begin{bmatrix} \lambda & 0 & 0 \\ 0 & 1 & 0 \\ 0 & 0 & \lambda^{-1} \end{bmatrix} 
\begin{bmatrix} \cos \theta & 0 & -\sin \theta \\ 0 & 1 & 0 \\ \sin \theta & 0 & \cos \theta \end{bmatrix} \bigg\vert \, \theta \in \left\{0,\frac{\pi}{8},\frac{\pi}{4},\frac{3\pi}{8}\right\} \right\} $$
for $\log \lambda = \cosh^{-1} (\cot \frac{\pi}{8})$. This group acts cocompactly on a complete, totally geodesic submanifold of $\X$ of constant curvature $-\frac{1}{3}$, see Section \ref{curvature}, with quotient a closed surface of genus $2$. A fundamental domain for this action is given by a regular octagon in $\HH^2$ with center $p$, the identity matrix in $\X$. This octagon decomposes into $16$ triangles with vertices at the center, the vertices of the octagon, and the midpoints of the edges. These triangles are isosceles with angles  $\frac{\pi}{2},\frac{\pi}{8},\frac{\pi}{8}$. By the hyperbolic law of cosines (for curvature $-\frac{1}{3}$),
$$ \cos \gamma = -\cos \alpha \cos \beta + \sin \alpha \sin \beta \cosh (\frac{1}{\sqrt{3}} c), $$
we see that the distance from the center $p$ to the vertex is $R=\sqrt{3}\cosh^{-1} \left( \cot^2 \frac{\pi}{8} \right) $. The $\Gamma_2$ translates of $B_R(p)$ cover $\HH^2$, so by the Milnor-Schwarz Lemma the orbit map $\orb_p \colon \Gamma_2 \to \HH^2$ is a $(1,1,2R+1,0)$-quasi-isometric embedding. One checks that $2R+1 \le 9.5$. Here we use the symmetric generating set $S'=\{ \gamma \in \Gamma_2 \mid d(p,\gamma p) \le 9.5 \}$. Note that the $S'$ here agrees with the one in the introduction because $d(p,\gamma p) = \sqrt{6}\abs{\log \gamma}_{Fr}$. Every geodesic in this copy of $\HH^2$ is $(\frac{1}{2\sqrt{3}},\sigma_{mod})$-regular in $\X$. Representations of this form were studied by Barbot in \cite{Bar10}. 

We may now prove

\begin{manualtheorem}{1.3}
	If $\rho \colon \Gamma_2 \to \SL(3,\R)$ is a representation satisfying $\abs{ \rho(s) - s}_{Fr} \le 10^{-3,698,433}$ for all $s \in S'$, then $\rho$ is Anosov. 
\end{manualtheorem}

\begin{proof}
	 From the classical Morse Lemma (Theorem \ref{classical morse lemma}), we get a Morse constant of $D=163$. Thus the orbit map at $p$ is a $(\frac{1}{2\sqrt{3}},\sigma_{mod},163)$-Morse $(1,1,9.5,0)$-quasiisometric embedding. We relax the additive parameters by $10$ and ask a perturbation to be a $(2.2 \times 10^6)$-local $(\frac{1}{2 \sqrt{3}},\sigma_{mod},173)$-Morse $(1,11,9.5,10)$-quasiisometric embedding. By Theorem \ref{local to global}, such an orbit map is a global $(\frac{1}{4\sqrt{3}},\sigma_{mod},6.8 \times 10^6)$-Morse $(108,214;1.4 \times 10^7;9.5;0)$-quasiisometric embedding.
	%corrected 3.5.2021
	
	If $\rho \colon \Gamma_2 \to \SL(3,\R)$ is another representation such that $\abs{\rho(s)-s}_{Fr} \le 10^{-3,698,433}$ then for $d_{S'}(w,1) \le k = 1.1 \times 10^6$ we have $d_{{\X}}(\rho(w)p,wp) \le 10$ by Corollary \ref{perturbed representation} so $\rho$ also induces a $(2.2 \times 10^6)$-local $(\frac{1}{2 \sqrt{3}},\sigma_{mod},173)$-Morse $(1,11,9.5,10)$-quasiisometric embedding and therefore a global $(\frac{1}{4\sqrt{3}},\sigma_{mod},6.8 \times 10^6)$-Morse $(108,214;1.4 \times 10^7;9.5;0)$-quasiisometric embedding. In particular, $\rho$ is Anosov and our proof of Theorem \ref{main theorem demonstration 2} is complete.
\end{proof}

% it would be nice to do this math and add this section.
%\subsection{Morse-Schottky groups}

\newpage

\printbibliography

\end{document}